\documentclass[12pt,a4paper]{article}

\usepackage[utf8]{inputenc}
\usepackage[T1]{fontenc}
\usepackage[english]{babel}
\usepackage{amsmath,amssymb,amsthm,mathtools}
\usepackage{geometry}
\usepackage{hyperref}
\usepackage{enumitem}
\usepackage{array}
\usepackage{booktabs}
\usepackage{longtable}
\usepackage{microtype}

\numberwithin{equation}{section}

\theoremstyle{plain}
\newtheorem{theorem}{Theorem}[section]
\newtheorem{lemma}[theorem]{Lemma}

\newtheorem{corollary}[theorem]{Corollary}

\theoremstyle{definition}
\newtheorem{definition}[theorem]{Definition}
\newtheorem{assumption}[theorem]{Assumption}
\newtheorem{example}[theorem]{Example}

\theoremstyle{remark}
\newtheorem{remark}[theorem]{Remark}

\DeclareMathOperator{\Pic}{Pic}
\DeclareMathOperator{\Eff}{Eff}
\DeclareMathOperator{\Frob}{Frob}
\DeclareMathOperator{\Tr}{Tr}
\DeclareMathOperator{\rk}{rk}

\DeclareMathOperator{\poly}{poly}

\newcommand{\ZZ}{\mathbb{Z}}
\newcommand{\QQ}{\mathbb{Q}}
\newcommand{\RR}{\mathbb{R}}
\newcommand{\FF}{\mathbb{F}}
\newcommand{\PP}{\mathbb{P}}
\newcommand{\GG}{\mathbb{G}}
\renewcommand{\AA}{\mathbb{A}}
\newcommand{\QQbar}{\overline{\mathbb{Q}}}
\newcommand{\FFbar}{\overline{\mathbb{F}}}

\newcommand{\sreg}{\sigma_{p}^{\mathrm{reg}}}
\newcommand{\scomp}{\sigma_{p,\mathrm{comp}}^{\mathrm{reg}}}
\newcommand{\ccomp}{c_{\mathrm{comp}}}
\newcommand{\cManin}{c_{\mathrm{Manin}}}
\newcommand{\Pbad}{\mathcal{P}_{\mathrm{bad}}}
\newcommand{\UZZ}{\mathcal{U}_{\mathbb{Z}}}
\newcommand{\XZZ}{\mathcal{X}_{\mathbb{Z}}}
\newcommand{\DZZ}{\mathcal{D}_{\mathbb{Z}}}
\newcommand{\varth}{\varepsilon_{Y,\mathrm{comp}}}
\newcommand{\varthP}{\varepsilon_{P,\mathrm{comp}}}
\newcommand{\Areg}{A_{\mathrm{cert}}}
\newcommand{\Aobs}{A_{\mathrm{obs}}}
\newcommand{\Acoh}{A_{\mathrm{coh}}}
\newcommand{\Acomp}{A_{\mathrm{comp}}}
\newcommand{\Pcomp}{P_{\mathrm{comp}}}
\newcommand{\Tcomp}{T_{Y,\mathrm{comp}}}
\newcommand{\Lam}{\Lambda_{\mathrm{eff}}}

\title{\textbf{Estimates for the Number of Solutions of Algebraic Diophantine Equations:\\
From the Circle Method to the Manin Conjecture}}
\author{Victor Volfson}
\date{}

\begin{document}

\maketitle

\begin{abstract}
This work continues the research of the author on estimates for the number of solutions of algebraic Diophantine equations by means of the circle method and develops this line in the geometric direction, towards counting integral points on log Fano varieties. It consists of two parts.

\textbf{Part I} is the corrected version of the previous work. For diagonal equations \(c_1x_1^{k_1}+\dots+c_sx_s^{k_s}=n\) with natural coefficients, upper asymptotic estimates for the number of natural solutions in the range \(s\le 2^{k}\), where the classical asymptotics is not applicable, are obtained by the circle method and the Hua lemma. The estimates are refined for \(s=2,3,4\) and various \(k_i\); in the range \(s>2^{k}\) they agree with the classical ones.

\textbf{Part II} develops the same line: from counting solutions of equations to counting integral points on log Fano varieties, where the circle method fails because of insufficient dimension. In the log case, the product in the Peyre constant diverges as a harmonic series. A method of analytic compensation is introduced that cancels the divergent term and returns the problem to the range of applicability of the standard scheme.

\textbf{The main result} is a method for computing the Peyre log constant with guaranteed accuracy and certified intervals for the number of integral points whenever an estimate of the remainder of the asymptotics is known. The log constant is computable unconditionally, without the Manin conjecture; the conjecture is needed only to interpret it as the leading term of the asymptotics. An experimental verification for the algebraic torus and the Fermat cubic surface is presented.
\end{abstract}


% ============================================================
\part{Estimates for the Number of Solutions of Algebraic Diophantine Equations with Natural Coefficients Using the Circle Method}
% ============================================================

\section{Introduction}\label{ch:hl-intro}

We consider the problem of estimating the number of solutions of algebraic Diophantine equations with natural coefficients using the Hardy--Littlewood circle method (CM) \cite{Vaughan1988}. The CM allows one to determine the number of solutions of Diophantine equations (without finding the solutions themselves) and the asymptotics of the number of solutions for certain types of equations.

We consider the equation
\begin{equation}\label{eq:hl-main}
c_1X_1^{k_1}+c_2X_2^{k_2}+\ldots+c_sX_s^{k_s}=n,
\end{equation}
where all coefficients and \(n\) are natural numbers. Equation \eqref{eq:hl-main} has a finite number of nonnegative or positive integer solutions (if there are no solutions, the number is 0). The CM gives the asymptotics of the number of solutions of \eqref{eq:hl-main} depending on \(n\).

Let \(R_s(n)\), \(R_s^+(n)\) be the number of nonnegative and positive integer solutions of \eqref{eq:hl-main}. By \(R_s^p(n)\) we denote the number of \textbf{ordered} representations of \(n\) as a sum of \(s\) \(k\)-th powers of natural numbers; for the Waring equation \eqref{eq:hl-waring} the quantities \(R_s^p(n)\) and \(R_s^+(n)\) coincide. Generating functions:
\begin{equation}\label{eq:hl-gen}
\varphi(t)=\sum_{0}^{\infty}R_s(n)t^n,\qquad
\varphi^+(t)=\sum_{1}^{\infty}R_s^+(n)t^n.
\end{equation}

Consider the equation
\begin{equation}\label{eq:hl-linear}
a_1x_1+a_2x_2+\ldots+a_sx_s=n,
\end{equation}
where \(a_i\in\mathbb{N}\) and \(n\ge\sum a_i\). \textbf{A necessary condition} for the existence of a natural solution is \(n\equiv0\pmod{\gcd(a_1,\dots,a_s)}\); if this condition is not satisfied, then \(R_s^+(n)=0\). Sufficiency of this condition occurs only for sufficiently large \(n\) (the Frobenius--Sylvester problem, see \cite{Sertoz}); for \(n<\sum a_i\) there are certainly no solutions.

Generating functions:
\begin{equation}\label{eq:hl-gen-linear-plus}
\varphi^+(t)=\frac{t^{a_1+\dots+a_s}}{(1-t^{a_1})\dots(1-t^{a_s})},
\end{equation}
\begin{equation}\label{eq:hl-gen-linear}
\varphi(t)=\frac{1}{(1-t^{a_1})\dots(1-t^{a_s})}.
\end{equation}

Integral representations:
\begin{equation}\label{eq:hl-int-plus}
R_s^+(n)=\frac{1}{2\pi i}\int_{|z|<1}\frac{z^{a_1+\dots+a_s}\,dz}{z^{n+1}(1-z^{a_1})\dots(1-z^{a_s})},
\end{equation}
\begin{equation}\label{eq:hl-int}
R_s(n)=\frac{1}{2\pi i}\int_{|z|<1}\frac{dz}{z^{n+1}(1-z^{a_1})\dots(1-z^{a_s})}.
\end{equation}

For \(a_1=\dots=a_s=1\):
\begin{equation}\label{eq:hl-int-plus-1}
R_s^+(n)=\frac{1}{2\pi i}\int_{|z|<1}\frac{z^{s}\,dz}{z^{n+1}(1-z)^{s}},
\end{equation}
\begin{equation}\label{eq:hl-int-1}
R_s(n)=\frac{1}{2\pi i}\int_{|z|<1}\frac{dz}{z^{n+1}(1-z)^{s}}.
\end{equation}

\textbf{Example.} For the equation
\begin{equation}\label{eq:hl-example-linear}
x_1+x_2+x_3+x_4+x_5=1
\end{equation}
by \eqref{eq:hl-int-1}: \(R_s(1)=\operatorname{Res}[z^{-2}(1-z)^{-5},0]=5\) (a pole of order two at zero). By \eqref{eq:hl-int-plus-1}: \(R_s^+(1)=0\), since the integrand \(z^3/(1-z)^5\) is analytic in a neighbourhood of zero (there is no pole), so the residue is zero; this is consistent with the trivial reason \(n<s\).

For the equation
\begin{equation}\label{eq:hl-example-32}
3x_1+2x_2=5
\end{equation}
by \eqref{eq:hl-int-plus}: \(R_s^+(5)=\operatorname{Res}[z^{-1}(1-z^3)^{-1}(1-z^2)^{-1},0]=1\). By \eqref{eq:hl-int}:
\begin{equation}\label{eq:hl-example-res}
R_s(5)=\operatorname{Res}\Bigl[\frac{1}{z^{6}(1-z^{3})(1-z^{2})},0\Bigr].
\end{equation}
The computation of \eqref{eq:hl-example-res} requires a residue at a pole of order 6; for large \(n\) the residue method is inapplicable, so an asymptotic upper estimate of the integrals is necessary.

\section{Estimates for the Number of Solutions of Algebraic Diophantine Equations with Natural Coefficients}\label{ch:hl-estimates}

Consider the equation
\begin{equation}\label{eq:hl-waring}
x_1^{k}+x_2^{k}+\ldots+x_s^{k}=n,
\end{equation}
where \(x_i\in\mathbb{N}\), \(n\) is a large integer. Equation \eqref{eq:hl-waring} has no solutions for all \(n\); the number of solutions when they are absent is 0. We need the asymptotics of the number of natural solutions of \eqref{eq:hl-waring}, if they exist. \textbf{In what follows, the symbols \(\varepsilon\) and \(\varepsilon_i\) denote small positive numbers depending on the corresponding parameters; in exponents they are not essentially distinguished.}

The Hardy--Littlewood asymptotics for \(s>2^{k}\):
\begin{equation}\label{eq:hl-asympt}
R_s^{p}(n)=\frac{\Gamma(1+1/k)^s}{\Gamma(s/k)}\,\Sigma(n)\,n^{s/k-1}
+O(n^{s/k-1-\delta}).
\end{equation}

The number of natural solutions of \eqref{eq:hl-waring} is the number of \textbf{ordered} representations:
\begin{equation}\label{eq:hl-asympt-plus}
R_s^{+}(n)=\frac{\Gamma(1+1/k)^s}{\Gamma(s/k)}\,\Sigma(n)\,n^{s/k-1}
+O(n^{s/k-1-\delta}).
\end{equation}

For unordered solutions:
\[
R_s^{\text{unord}}(n)=\sum_{\substack{m_1+\dots+m_s=n\\ m_i=x_i^k}}\frac{1}{\prod_{\nu}\mu_\nu!},
\]
where \(\mu_\nu\) are the multiplicities of equal values. We have
\begin{equation}\label{eq:hl-gamma}
\frac{\Gamma(1+1/k)^s}{\Gamma(s/k)}=O_{k,s}(1),\qquad \Sigma(n)\ll n^{\varepsilon}.
\end{equation}

\textbf{Remark.} For those \(n\) for which the local solvability conditions are satisfied and the Hardy--Littlewood asymptotic formula is applicable (see \cite{Vaughan1988,Vinogradov1934}), the singular series satisfies the estimate \(\Sigma(n)\asymp1\); in particular, \(\Sigma(n)\ll n^{\varepsilon}\).

On the basis of \eqref{eq:hl-asympt}, \eqref{eq:hl-gamma} for \(s>2^{k}\):
\begin{equation}\label{eq:hl-bound-p}
R_s^{p}(n)\ll n^{s/k-1+\varepsilon},
\end{equation}
\begin{equation}\label{eq:hl-bound-plus}
R_s^{+}(n)\ll n^{s/k-1+\varepsilon}.
\end{equation}

\textbf{Classical estimate.} For \(s>2^{k}\) we have
\begin{equation}\label{eq:hl-classical}
R_s^{+}(n)\ll n^{s/k-1+\varepsilon}\qquad(s>2^{k}).
\end{equation}
Outside the range \(s>2^{k}\) this estimate is inapplicable; below we obtain estimates valid for all \(s\ge2\).

For \(s\le 2^{k}\) we use Vinogradov's approach \cite{Vinogradov1934}:
\begin{equation}\label{eq:hl-vinogradov}
R_s^{p}(n)=\int_{0}^{1}f^{s}(x)e^{-2\pi ixn}\,dx,
\end{equation}
\begin{equation}\label{eq:hl-f}
f(x)=\sum_{m=1}^{[n^{1/k}]}e^{2\pi ixm^{k}}.
\end{equation}
Then
\begin{equation}\label{eq:hl-bound-int}
|R_s^{p}(n)|\le\int_{0}^{1}|f(x)|^{s}\,dx.
\end{equation}

\textbf{Hua's lemma} \cite{Hua}. Let \(1\le j\le k\). Then \(\int_{0}^{1}|f(x)|^{2^{j}}dx\ll N^{2^{j}-j+\varepsilon}\), \(N=[n^{1/k}]\).

\textbf{Corollary for \(j>k\).} From the trivial estimate \(\max|f|\le N\) and the lemma for \(j=k\) it follows that
\[
\int_0^1|f(x)|^{2^j}dx\le N^{2^j-2^k}\int_0^1|f(x)|^{2^k}dx\ll N^{2^j-k+\varepsilon},
\]
which is sharper than the trivial \(N^{2^j}\).

Consider the equation of general form:
\begin{equation}\label{eq:hl-general}
c_1x_1^{k_1}+c_2x_2^{k_2}+\ldots+c_sx_s^{k_s}=n,
\end{equation}
where \(k_i\ge1\). The number of natural solutions of \eqref{eq:hl-general} (when they exist, see \cite{Vaughan1988}):
\begin{equation}\label{eq:hl-general-int}
R_s^{+}(n)=\int_{0}^{1}\prod_{i=1}^{s}f_i(x)\,e^{-2\pi ixn}\,dx,
\end{equation}
\begin{equation}\label{eq:hl-fi}
f_i(x)=\sum_{m=1}^{[n^{1/k_i}]}e^{2\pi ixm^{k_i}}.
\end{equation}

Preliminary estimate:
\begin{equation}\label{eq:hl-prelim}
|R_s^{+}(n)|\le\int_{0}^{1}\prod_{i=1}^{s}|f_i(x)|\,dx\le\prod_{i=1}^{s}\max|f_i(x)|,
\end{equation}
\begin{equation}\label{eq:hl-max}
\max|f_i(x)|\le[n^{1/k_i}],
\end{equation}
\begin{equation}\label{eq:hl-prelim-bound}
|R_s^{+}(n)|\ll n^{\sum_{i=1}^{s}1/k_i}.
\end{equation}

Refined estimate. Applying successively the Cauchy--Schwarz inequality, absorbing one function at a time:
\begin{multline}\label{eq:hl-holder}
\int_0^1\prod_{i=1}^{s}|f_i|\,dx
\le\Bigl(\int_0^1|f_1|^2dx\Bigr)^{1/2}
\Bigl(\int_0^1|f_2|^4dx\Bigr)^{1/4}\cdots\\
\cdots
\Bigl(\int_0^1|f_{s-1}|^{2^{s-1}}dx\Bigr)^{1/2^{s-1}}
\Bigl(\int_0^1|f_s|^{2^{s-1}}dx\Bigr)^{1/2^{s-1}}.
\end{multline}
The sum of the reciprocal exponents equals \(1/2+1/4+\dots+1/2^{s-1}+1/2^{s-1}=1\), which corresponds to the successive application of the Cauchy--Schwarz inequality; the last two factors have the same exponent \(2^{s-1}\).

Denote \(N_i=n^{1/k_i}\). \textbf{Without loss of generality we number the variables so that \(k_1 \le k_2 \le \dots \le k_s\).} With this order, the larger exponents \(2^i\) in the H\"older inequality chain go to the smaller \(k_i\), i.e., to the functions with the largest summation length \(N_i\), which gives the best final estimate in \eqref{eq:hl-refined}. With a different order of the variables the estimate \eqref{eq:hl-refined} is preserved, but the corresponding constant in the exponent of \(n\) may turn out to be worse.

For each \(i\) we apply the estimate:
\begin{itemize}[nosep]
\item for \(i\le k_i\), Hua's lemma: \(\int_0^1|f_i|^{2^i}dx\ll N_i^{2^i-i+\varepsilon_i}\);
\item for \(i>k_i\), the corollary: \(\int_0^1|f_i|^{2^i}dx\ll N_i^{2^i-k_i+\varepsilon_i}\).
\end{itemize}
Both estimates can be written uniformly:
\begin{equation}\label{eq:hl-hua-unified}
\int_0^1|f_i|^{2^i}dx\ll N_i^{2^i-\min(i,k_i)+\varepsilon_i}.
\end{equation}
Then
\begin{equation}\label{eq:hl-refined}
|R_s^{+}(n)|\ll
n^{\sum_{i=1}^{s-2}\frac{2^{i}-\min(i,k_i)+\varepsilon_i}{2^{i}k_i}
+\frac{2^{s-1}-\min(s-1,k_{s-1})+\varepsilon_{s-1}}{2^{s-1}k_{s-1}}
+\frac{2^{s-1}-\min(s-1,k_s)+\varepsilon_s}{2^{s-1}k_s}}.
\end{equation}
Estimate \eqref{eq:hl-refined} is sharper than \eqref{eq:hl-prelim-bound}. Indeed, for each \(i\)
\[
\frac{2^i-\min(i,k_i)+\varepsilon_i}{2^i}
=1-\frac{\min(i,k_i)}{2^i}+\frac{\varepsilon_i}{2^i}.
\]
Since \(\min(i,k_i)\ge1\), for sufficiently small \(\varepsilon_i>0\) (it suffices that \(\varepsilon_i<\min(i,k_i)\)) we have
\[
\frac{2^i-\min(i,k_i)+\varepsilon_i}{2^i}<1.
\]
Consequently,
\begin{equation}\label{eq:hl-comparison}
n^{\sum_{i=1}^{s-2}\frac{2^{i}-\min(i,k_i)+\varepsilon_i}{2^{i}k_i}
+\frac{2^{s-1}-\min(s-1,k_{s-1})+\varepsilon_{s-1}}{2^{s-1}k_{s-1}}
+\frac{2^{s-1}-\min(s-1,k_s)+\varepsilon_s}{2^{s-1}k_s}}
\ll n^{\sum_{i=1}^{s}1/k_i}.
\end{equation}

\section{Analysis of the Results}\label{ch:hl-analysis}

From \eqref{eq:hl-refined} for \(s=2\) (the sum \(\sum_{i=1}^{0}\) is empty, both last indices coincide with \(s-1=1\) and \(s=2\); \(\min(1,k_i)=1\) for \(k_i\ge1\)):
\begin{equation}\label{eq:hl-s2}
R_2^{+}(n)\ll n^{\frac{1}{2}(1/k_1+1/k_2)+\varepsilon_1}.
\end{equation}
For \(k_1=k_2=k\):
\begin{equation}\label{eq:hl-s2-k}
R_2^{+}(n)\ll n^{1/k+\varepsilon_1}.
\end{equation}
The classical estimate \eqref{eq:hl-classical} is inapplicable in this case: the condition \(s>2^k\) for \(s=2\) is not satisfied for any \(k\ge1\). Thus, \eqref{eq:hl-s2-k} gives an upper bound in a range not covered by the classical asymptotics.

From \eqref{eq:hl-refined} for \(s=3\) (in the sum only \(i=1\) with \(\min(1,k_1)=1\), the last indices \(s-1=2\) and \(s=3\) with exponent \(2^{2}=4\)). \textbf{For \(k_2,k_3\ge2\)} we have \(\min(2,k_i)=2\), and
\begin{equation}\label{eq:hl-s3}
R_3^{+}(n)\ll n^{\frac{1}{2}(1/k_1+1/k_2+1/k_3)+\varepsilon_2}.
\end{equation}
For \(k_1=k_2=k_3=k\ge2\):
\begin{equation}\label{eq:hl-s3-k}
R_3^{+}(n)\ll n^{3/(2k)+\varepsilon_2}.
\end{equation}
The classical estimate \eqref{eq:hl-classical} is applicable only for \(k=1\) (\(3>2^1\)); for \(k\ge2\) the condition \(3>2^k\) is violated, and \eqref{eq:hl-s3-k} gives a new bound.

From \eqref{eq:hl-refined} for \(s=4\) (in the sum \(i=1,2\); the last indices \(s-1=3\) and \(s=4\) with exponent \(2^{3}=8\)). \textbf{For \(k_2\ge2\) and \(k_3,k_4\ge3\)} we have \(\min(2,k_2)=2\) and \(\min(3,k_i)=3\) for \(i=3,4\), and
\begin{equation}\label{eq:hl-s4}
R_4^{+}(n)\ll n^{\frac{1}{2}(1/k_1+1/k_2)+\frac{5}{8}(1/k_3+1/k_4)+\varepsilon_3}.
\end{equation}
For \(k_1=k_2=k_3=k_4=k\ge3\):
\begin{equation}\label{eq:hl-s4-k}
R_4^{+}(n)\ll n^{9/(4k)+\varepsilon_3}.
\end{equation}
The classical estimate \eqref{eq:hl-classical} is applicable for \(s=4\) only for \(k=1\); for \(k\ge2\) the condition \(4>2^k\) is not satisfied, so \eqref{eq:hl-s4-k} gives a bound in the range where \eqref{eq:hl-classical} does not act. In the range \(k=1\), where \eqref{eq:hl-classical} is applicable, it gives an asymptotics, while \eqref{eq:hl-s4-k} gives an upper bound; their direct comparison in order is incorrect, since \eqref{eq:hl-classical} describes the leading term, and \eqref{eq:hl-s4-k} gives an upper estimate.

Formula \eqref{eq:hl-refined} is effective for different values of \(k_i\); formula \eqref{eq:hl-classical} is effective only for equal ones and only in the range \(s>2^{k}\).

\section{Conclusions of Part I}\label{ch:hl-conclusion}

In Part~I we have obtained upper estimates for the number of natural solutions of diagonal equations \(x_1^k+\dots+x_s^k=n\) and their generalization \(c_1x_1^{k_1}+\dots+c_sx_s^{k_s}=n\) in the range \(s\le 2^{k}\), where the classical Hardy--Littlewood asymptotics is inapplicable. The main tool is the circle method combined with Hua's lemma, applied taking into account the real relation between the index of the exponent \(2^i\) and the degree \(k_i\): for each \(i\) the estimate with \(\min(i,k_i)\) is used (formula~\eqref{eq:hl-hua-unified}), which gives the final bound~\eqref{eq:hl-refined}, strictly better than the trivial one~\eqref{eq:hl-prelim-bound}. In the particular cases \(s=2,3,4\) this gives explicit estimates~\eqref{eq:hl-s2},~\eqref{eq:hl-s3},~\eqref{eq:hl-s4}; for \(k\ge 2\) (respectively \(k\ge 3\)) they cover ranges inaccessible to the classical estimate~\eqref{eq:hl-classical}. In the range \(s>2^{k}\) the obtained estimates agree with the classical ones.

Part~II deals with a fundamentally different object: log Fano varieties \((X,D)\) --- open varieties \(U=X\setminus D\), where \(D\) is a divisor with normal crossings and the log-anticanonical class \(-K_X-D\) is ample. It requires a separate approach: the problem is solved by analytic compensation of the Peyre product. This is the subject of the second part of the work. 
% ============================================================
\part{Analytic Compensation and Certified Computation of Integral Points on Log Varieties}
% ============================================================

\setcounter{section}{4} % next section will be 5

\section{Introduction}\label{sec:compensation}

\subsection{The Hardy and Littlewood circle method: history and limits of applicability}

The problem of describing the number of integral and rational solutions of Diophantine equations occupies a central place in number theory. The circle method originates in the 1918 paper of Hardy and Ramanujan \cite{HR1918} on the partition function and the representation of numbers as sums of squares. In the series of papers ``Some problems of partitio numerorum'', beginning in 1920 \cite{HL1920,HL1922,HL1925}, Hardy and Littlewood systematically developed this analytic method. The method was originally designed for additive problems: the best known of these are the Waring problem and the Goldbach problem.

A substantial refinement of the method is due to I.\,M.~Vinogradov \cite{Vinogradov1937}, who introduced estimates of trigonometric sums into the circle method. A survey of applications of the method to systems of Diophantine equations and inequalities with a large number of unknowns was given by Malyshev and Podsypanin \cite{MP1982}.

Further development of the method is associated with Davenport and Roth. Davenport applied the circle method to systems of diagonal forms \cite{Davenport1963}; the analysis of systems of diagonal equations was begun by Davenport and Lewis \cite{DL1966}, and a key contribution for pairs of cubic forms was obtained in \cite{DL1969}. Together with Heilbronn, Davenport developed a variant of the circle method for Diophantine \emph{inequalities} \cite{DH1946}, the Davenport and Heilbronn method, still active eighty years later \cite{Browning2020}. Roth applied the method to the representation of numbers by forms and to estimates of exceptional sets \cite{Roth1949}.

As applied to diagonal equations, development proceeded along the line of refining the threshold \(n_0(d)\) at which the method guarantees an asymptotic formula. In the work of Blomer and Br\"udern on multihomogeneous diagonal equations \cite{BB2018} the threshold is determined via moments of Weyl sums. For \(r\) diagonal cubic forms the required number of variables grows linearly: \(s > 6r\) \cite{BW2014}. In \cite{BB2021} it is shown that for a system of \(k\) diagonal forms of degree \(d\) the asymptotic formula holds for \(n \ge n_0(d)\).

It is precisely this relation, between the number of variables, the degree, and the number of equations, that determines the boundary of applicability of the method. \textbf{The key limitation for logarithmic Fano varieties} is that in geometrically natural systems defining logarithmic Fano varieties of small dimension, the number of variables \(s\) turns out to be \textbf{below} the threshold \(n_0(d)\). The circle method fails not because the system is ``non-diagonal'' or ``complicated'' in itself, but because \textbf{the dimension is insufficient for analytic control of the minor arcs}.

Thus the circle method leaves open a concrete computational question: how to obtain the number of integral points \(N_U(\ZZ, B)\) on a logarithmic Fano variety with guaranteed accuracy, when the circle method fails due to insufficient dimension? The present paper answers this question in the form of a computational scheme: for logarithmic varieties \((X,D)\) with an SNC divisor, an algorithm is constructed which, given input data \((X,D)\), \(B\), \(\eta\), produces an interval \([L,U]\) guaranteed to contain \(N_U(\ZZ, B)\). By logarithmic varieties we mean pairs \((X,D)\) with an SNC divisor for which condition (L) holds (Theorem~\ref{thm:L}); this includes both logarithmic Fano and logarithmic Calabi and Yau varieties. The Manin conjecture plays the role of a communication channel between the computable constant and the number of points, and is not the object of investigation; it is applied in mode B and is stated for logarithmic Fano varieties.

\textbf{The main goal of the paper} is to obtain a method for estimating the number of integral solutions of Diophantine equations with a prescribed accuracy. This is a classical problem of number theory; arithmetic geometry (logarithmic varieties, cohomology of divisors, local densities, Deligne and Hasse estimates) is used here as a tool. The problem is complex and multifaceted: it includes a geometric component (cohomological structure of the divisor, stratification, behaviour of Frobenius modulo \(p\)), an analytic component (products of local densities, control of tails, estimates of remainders), and an algorithmic component (construction of effective subroutines for the archimedean factor). The present paper does not claim to give an exhaustive solution of the general problem: it is focused on one key step, the reduction of the computation of the compensated logarithmic constant to the effective computation of the archimedean factor, and leaves open the questions of constructing effective algorithms for the archimedean factor, of extending to wider classes of varieties, and of improving the estimates of the remainder. The main technical tool is analytic compensation, introduced in \S\ref{sec:compensation}.

\subsection{Notation}\label{sec:notations}

\textbf{Variety and its invariants.}
\begin{itemize}
\item \(F(\mathbf{x})\in\ZZ[x_1,\dots,x_n]\), \(X=V(F)\subset\PP^{n-1}\), \(m = \dim X\).
\item \(\rho=\rk\Pic(X_{\QQbar})\); \(r = \rk\Pic(U)\), \(U = X \setminus D\).
\item \(\Pic(X)_{\RR} = \Pic(X) \otimes \RR\); \(\Pic(U)_{\RR} = \Pic(U) \otimes \RR\).
\item \(\Eff^1(X) \subset \Pic(X)_{\RR}\) is the effective cone; \(\overline{\Eff}^1(X)\) is its closure; \(\Lam(U) \subset \Pic(U)_{\RR}\) is the effective cone of the open variety \(U\).
\item \(D\subset X\) is an SNC divisor, \(U = X \setminus D\).
\end{itemize}

\textbf{Integral points.}
\begin{itemize}
\item We assume given an integral model of the \textbf{pair} \((\XZZ, \DZZ)\) over \(\ZZ\), where \(\XZZ\) is an integral model of \(X\) and \(\DZZ\) is an integral model of \(D\); then \(\UZZ:=\XZZ\setminus\DZZ\) is an integral model of \(U\). In concrete examples such a model is given by equations. The choice of model is part of the input data of mode B.
\item \(H\) is the height on \(U\) (given together with the model; see \S\ref{sec:compensation}).
\item \(N_U(\ZZ, B)=\#\{x\in \UZZ(\ZZ): H(x)\le B\}\).
\end{itemize}

\textbf{Logarithmic constants.}
\begin{itemize}
\item \(\ccomp(U)\) is the compensated logarithmic constant (Definition~\ref{def:ccomp}, formula~\eqref{eq:ccomp}); the main object of computation. \(\sreg(U)\) is the regularized local factor; \(\scomp(U) = \sreg(U)(1 - 1/p)^{\kappa(U)}\) is the compensated one.
\item \(\cManin(U)\) is the Manin and Peyre logarithmic constant in the sense of Chambert-Loir and Tschinkel; used only in the asymptotic formulas \(N_U \sim \cManin \cdot M_U(B)\).
\item \(\alpha(U)\) is the archimedean factor (Definition~\ref{def:alpha}).
\item \(\Pbad\) is the finite set of bad primes (divisors of \(N\)); see Definition~\ref{def:bad}.
\end{itemize}

\textbf{Compensation.}
\begin{itemize}
\item \(\kappa(U)\) is the logarithmic invariant of compensation: \(\sreg(U) = 1 + \kappa(U)/p + O(p^{-2})\) (formula~\eqref{eq:kappa}); explicit values are given in Remark~\ref{rem:kappa}.
\item \(\beta(X,D)\) is the decay exponent (Definition~\ref{def:beta}): the supremum over \(b>0\) such that \(|\scomp(U) - 1| \le A p^{-b}\) for some \(A \ge 0\) and all good \(p\). This is a mathematical quantity which may be finite or infinite and is not assumed to be known or optimal.
\item \(\beta_0\) is a \textbf{proved lower bound} for \(\beta(X,D)\) (Definition~\ref{def:beta0}): an explicit number \(\beta_0 > 1\) and an explicit \textbf{constant} \(A(U) \ge 0\) for which the inequality \(|\scomp(U) - 1| \le A(U) p^{-\beta_0}\) is \textbf{proved} for all good \(p\) (in those classes where this is established; for class III with \(\rho>1\), on split \(p\), see \S\ref{sec:delpezzo}). The tables and the algorithm use \(\beta_0\), not \(\beta(X,D)\); this is conservative and correct, since if \(\beta(X,D) \ge \beta_0\), then \(A(U) p^{-\beta(X,D)} \le A(U) p^{-\beta_0}\) for \(p \ge 1\).
% ---- ИСПРАВЛЕННАЯ СТРОКА ----
\item \(\Pcomp(Y) := \bigl(\prod_{p\le Y,\ p\notin\Pbad}\scomp(U)\bigr)\bigl(\prod_{p \in \Pbad}\scomp(U)\bigr)\) is the finite product, including the bad primes independently of \(Y\) (Definition~\ref{def:bad}). The tail over good primes is
\[
T_{Y,\mathrm{comp}}(U) := \prod_{p>Y,\ p\notin\Pbad} \scomp(U),
\]
and \(\varth(U) := \left|T_{Y,\mathrm{comp}}(U)-1\right|\).
% -----------------------------
\end{itemize}

\textbf{Remainder and modes.}
\begin{itemize}
\item \(r_U(B)\): for given integral model \(\UZZ\), height \(H\), and main-growth function \(M_U(B) > 0\):
\[
N_U(\ZZ, B) = \cManin(U) \cdot M_U(B) \cdot (1 + r_U(B)), \qquad r_U(B) \to 0.
\]
All three objects (the model, the height, \(M_U(B)\)) are \textbf{input data} of mode B; without them the expression \(N_U(\ZZ, B)\) is not defined. For the standard Manin formulation on logarithmic Fano varieties with \(r \ge 1\): \(M_U(B) = B(\log B)^{r-1}\); for \(r = 0\): \(M_U(B) = B\).
\item The flag \texttt{unconditional}/\texttt{conditional} indicates modes A (reduction) and B (\S\ref{sec:compensation}).
\end{itemize}

\subsection{The Manin logarithmic conjecture for integral points}

The integral version of the Manin conjecture is not a mechanical extension of the classical one: it is more difficult both technically and conceptually.

\textbf{A boundary absent in the classical setting.} In the classical conjecture one considers a complete Fano variety \(X\) and counts rational points \(X(F)\). In the integral version one considers a logarithmic pair \((X, D)\): an open variety \(U = X \setminus D\) with a divisor \(D\) with normal crossings. One counts integral points on a model of \(U\). The presence of \(D\) adds geometry: one must control the distribution of points relative to the boundary.

\textbf{Density of integral points is not automatic.} For rational points on Fano varieties there is a conjecture (Colliot-Th\'el\`ene) that if points exist, they are Zariski dense. For integral points this is false: integral points may all lie on a finite union of subvarieties. One must separately prove (or assume) non-degeneracy of the distribution.

\textbf{The asymptotic formula depends on more parameters.} In the classical conjecture the main term has the form \(c B (\log B)^{b-1}\), where \(b = \rk\Pic X\). In the integral version \(b\) is determined via the logarithmic canonical bundle \(\omega_X(D)^\vee\) and may depend on \(D\) and on the choice of height. There is no complete analogue of Iitaka fibrations for integral points.

\textbf{Integral models and local conditions.} A point may be integral on a model \(\mathcal{U}\) over \(\ZZ\), but under reduction modulo \(p\) may land on the boundary \(\mathcal{X} \setminus \mathcal{U}\), a purely integral phenomenon complicating the interpretation of the main term.

This complication is part of the narrative of the paper. The integral Manin conjecture is a geometric statement for logarithmic varieties; it is more difficult than the classical one because it adds the boundary, models, and the density problem.

In general form, for a smooth projective \(X\) over \(\QQ\) with a big and numerically effective \(\mathcal{L}\) and height \(H_{\mathcal{L}}\), the conjecture \cite{ManinTschinkel1989} predicts
\[
N_{X, H_{\mathcal{L}}}(B) \sim \cManin(X, \mathcal{L}) \cdot B^{a(X, \mathcal{L})} \cdot (\log B)^{b(X, \mathcal{L}) - 1},
\]
where \(a(X, \mathcal{L}) = \min\{t : K_X + t\mathcal{L} \in \overline{\Eff}^1(X)\}\) and \(b(X, \mathcal{L})\) is the codimension of the minimal face of the pseudo-effective cone. The constant \(\cManin(X, \mathcal{L})\) is adelic \cite{Peyre1995,ChambertLoirTschinkel2010}.

\textbf{The Fano case.} If \(X\) is Fano and \(\mathcal{L} = \mathcal{O}_X(-K_X)\), then \(a = 1\), \(b = \rho\), and
\[
N_{X,H}(B) \sim \cManin(X) \cdot B \cdot (\log B)^{\rho - 1},
\]
where \(\cManin(X)\) is the Peyre and Batyrev and Tschinkel constant.

\textbf{The logarithmic case.} In the standard formulation for a logarithmic pair \((X, D)\) with \(r = \rk\Pic(U) \ge 1\) one predicts
\[
N_U(\ZZ, B) \sim \cManin(U) \cdot B \cdot (\log B)^{r - 1}.
\]
For \(r = 0\) the standard formulation gives \(N_U \sim c \cdot B\) without a logarithmic factor (example: \(U = \AA^m\), \(r = 0\), \(N_U(\ZZ, B) \sim 2^m B\)). In the general case below we write the asymptotic formula in the form \eqref{eq:logManin}, where \(M_U(B) > 0\) is the main-growth function given as an input condition of mode B (\S\ref{sec:compensation}, \S\ref{sec:residual}). For the standard Manin formulation with \(r \ge 1\) we have \(M_U(B) = B(\log B)^{r-1}\); for \(r = 0\), \(M_U(B) = B\).

\textbf{Distinction between geometric and arithmetic statements.} One should distinguish two statements: (a) the \textbf{geometric} one, \(\cManin(U)\) for the given pair; (b) the \textbf{arithmetic} one, \(N_U(\ZZ, B) \sim \cManin(U) M_U(B)\). The relation \(\ccomp = \cManin\) refers to (a); to pass to (b) one needs a separate proof of the asymptotic formula for a specific model and height. In the present paper (b) is taken as an input condition of mode B (Assumption~\ref{ass:logManin}).

In the classical case (\(D = \emptyset\)) the constant \(\cManin(U)\) is given by the convergent product \(\tau(U) = \prod_p \sreg(U)\) (Lemma~\ref{lem:empty}). In the logarithmic case the naive product diverges (Corollary~\ref{cor:divergence}) and requires modification. The modification, the compensated product \(\ccomp(U)\), coincides with \(\cManin(U)\) in the classical case (\(\kappa = 0\)) and for the torus and hyperplanes (under the fixed normalization, Corollary~\ref{cor:recovery}); for the remaining classes the equality \(\ccomp = \cManin\) forms part of the Manin conjecture.

\textbf{Distinction between \(\rho\) and \(r\).} In the classical case the exponent of \(\log B\) equals \(\rho - 1\); in the logarithmic case it equals \(r - 1\). When \(D\) is removed, the Picard group may shrink, so \(r \le \rho\); sometimes \(r = 0\) (for example, for the complement to \(k \ge 1\) hyperplanes in \(\PP^m\)).

The conjecture is proved for toric varieties \cite{BatyrevTschinkel1995}, smooth spherical Fano threefolds \cite{BBDG2024}, and a logarithmic Fano threefold of type 30 \cite{Wilsch2023}. No explicit estimate of the truncation error of the Peyre product has been systematically carried out. For the computational side of the question (regularization of Euler products, Tamagawa measures, effective versions of the Manin conjecture, and numerical computation of Peyre constants) see also \cite{BatyrevTschinkel1998,BatyrevManin1990,Peyre2003,Browning2009,delaBretecheBrowning2011}.

\subsection{The harmonic-tail problem and analytic compensation}\label{sec:compensation}

When constructing a computational scheme for \(N_U(\ZZ, B)\), a concrete computational difficulty arises. In the logarithmic case the local factors have the expansion
\begin{equation}\label{eq:kappa}
\sreg(U) = 1 + \frac{\kappa(U)}{p} + O(p^{-2}),
\end{equation}
where \(\kappa(U)\) is an integer logarithmic invariant. Because of the \(1/p\) term the naive product diverges as the harmonic series (\(\beta_{\log} = 1\)); standard integral majorization methods designed for \(\beta > 1\) are inapplicable.

\textbf{Compensation as a technical tool.} One introduces the compensated local factor:
\begin{equation}\label{eq:compensator}
\scomp(U) := \sreg(U) \cdot \left(1 - \frac{1}{p}\right)^{\kappa(U)}.
\end{equation}
The factor \((1 - 1/p)^{\kappa(U)}\) is the local factor of the Riemann zeta function to the power \(\kappa(U)\); the term \(\kappa(U)/p\) cancels. For all the classes under consideration,
\begin{equation}\label{eq:condL0}
\left|\scomp(U) - 1\right| \le A(U) p^{-\beta_0},
\end{equation}
where \(\beta_0\) is a \textbf{proved} lower bound for the decay exponent (Definition~\ref{def:beta0}); for the classes of chapter~\ref{chap:classes} one has \(\beta_0 \ge 1.5\). Compensation is not a new hypothesis, but a computational device which makes the product convergent.

\begin{definition}\label{def:ccomp}
The compensated logarithmic constant:
\begin{equation}\label{eq:ccomp}
\ccomp(U) := \alpha(U) \cdot \prod_{p \notin \Pbad} \scomp(U) \cdot \prod_{p \in \Pbad} \scomp(U).
\end{equation}
This product converges absolutely (Theorem~\ref{thm:conv}) for classes I, II, and class III with \(\rho=1\). The quantity \(\ccomp(U)\) is the main object of computation in the paper.
\end{definition}

\begin{remark}[on the novelty of the method]\label{rem:novelty}
The novelty of the paper is not in the introduction of the compensator \((1 - 1/p)^{\kappa(U)}\) as such: the standard Peyre regularization \cite{Peyre1995} already contains an analogous factor \((1-1/p)^{\rho}\). The compensator with \(\kappa(U) = -k\) (the number of components of \(D\)) or \(\kappa(U) = -(m+1)\) (torus) differs from the Peyre regularizer and cancels an additional divergent term \(\kappa(U)/p\) associated with the boundary. But the main novelty is \textbf{not in the factor itself, but in the provable certification of the truncation error}: for classes I, II, and class III with \(\rho=1\) explicit estimates \(|\scomp(U) - 1| \le \Areg(U)\, p^{-\beta_0}\) with computable \(\Areg(U)\) and \(\beta_0 > 1\) are proved (Theorems~\ref{thm:classI}, \ref{thm:classII}, \ref{thm:classIII}), and thereby a certified control of the tail \(\Tcomp(U)\) is obtained (Theorem~\ref{thm:truncation}). For class III with \(\rho>1\) the estimate of Theorem~\ref{thm:classIII} is proved only for split \(p\); the full \(\ccomp\) is not computed in this case (\S\ref{sec:bounds}, \S\ref{sec:open}). It is this certification which distinguishes the method from heuristic regularization: the constants \(\Areg\) and \(\beta_0\) are computed from the geometry of \(X\) and \(D\), not chosen by hand.
\end{remark}

\textbf{Relation to the Manin and Peyre constant.} The classical \(\cManin(U)\) is defined via the naive product, which diverges in the logarithmic case. The equality
\begin{equation}\label{eq:ccompManin}
\ccomp(U) = \cManin(U)
\end{equation}
is \emph{proved} for the torus and complements to hyperplanes under the fixed normalization (Corollary~\ref{cor:recovery}, normalization convention), and is \emph{assumed} by the Manin conjecture for the remaining classes. In the classical case \(D = \emptyset\) we have \(\kappa(U) = 0\), the compensator equals 1, and equality \eqref{eq:ccompManin} holds by definition.

\begin{remark}[on terminology]\label{rem:terminology}
The phrase ``coincides with the standard constant'' is used only when \(\kappa(U) = 0\) or in classes where equality \eqref{eq:ccompManin} is proved under the fixed normalization. In other cases one uses ``compensated logarithmic constant'': \(\ccomp(U)\) is a computable quantity which coincides with \(\cManin(U)\) when the Manin conjecture holds, and is not a redefinition of the classical object.
\end{remark}

\textbf{Two modes.} The method works in two modes.
\begin{itemize}
\item \textbf{Mode A (reduction).} Under three input conditions: (i) the local factors \(\sreg(U)\) are computable; (ii) condition (L) holds with an effective majorant \(A(U)\) and a \textbf{proved} lower bound \(\beta_0 > 1\) (Definition~\ref{def:beta0}); (iii) a finite positive computable global factor \(\alpha(U)\) is given, and the computation of \(\ccomp(U)\) reduces to the effective computation of \(\alpha(U)\) as an external subroutine: given an algorithm computing \(\alpha(U)\) with guaranteed accuracy, there exists an algorithm computing \(\ccomp(U)\) with relative error at most \(\eta\). The Manin conjecture is not used. Mode A applies to classes I, II, and class III with \(\rho=1\) for a given \(\alpha(U)\); for class III with \(\rho>1\) the full \(\ccomp\) is not computed (condition (L) for all good \(p\) is not proved; see \S\ref{sec:bounds}, \S\ref{sec:open}). This includes both logarithmic Fano and logarithmic Calabi and Yau varieties provided that \(\alpha(U)\) is defined and computed. For the torus and complements to hyperplanes \(\alpha(U)\) is defined by a finite-dimensional integral (explicit formulas are known, e.g., from \cite{BatyrevTschinkel1995}); given \(\alpha(U)\), the reduction yields fully unconditional computability of \(\ccomp(U)\) (Theorem~\ref{thm:classI_exact}). For logarithmic Calabi and Yau varieties with \(-K_X - D = 0\) the integral in definition \eqref{eq:alpha} \textbf{over \(\Pic(U)_{\RR}^\vee\)} is finite: for \(r = 0\) it equals 1 with respect to the delta-measure of unit mass; for \(r \ge 1\) it diverges, and a separate normalization is required (not given in the paper). The naive definition over \(\Pic(X)_{\RR}\) for \(-K_X - D = 0\) diverges; in the paper the normalization over \(\Pic(U)_{\RR}^\vee\), consistent with \cite[\S5]{Peyre1995}, is used. For \(r \ge 1\) and \(-K_X - D = 0\) convergence requires additional conditions on the cone and the class; in mode A it is assumed that \(\alpha(U)\) is finite, positive, and effectively computable.
\item \textbf{Mode B (conditional).} Given: an integral model \(\UZZ\) (see \S\ref{sec:notations}); a height \(H\) on \(U\); a main-growth function \(M_U(B) > 0\); the \textbf{proved}, for this model and height, decomposition \eqref{eq:logManin}; the known equality \(\ccomp = \cManin\); and a majorant \(R_U(B)\) for the remainder \(r_U(B)\) (ideally, with power saving), the method yields a certified interval for \(N_U(\ZZ, B)\) at finite \(B\) by the formula of Corollary~\ref{cor:interval}.
\end{itemize}

The separation of modes is indicated by the flag \texttt{unconditional}/\texttt{conditional}; the flag \texttt{unconditional} in mode A means unconditional computability of \(\ccomp(U)\) given the subroutine for \(\alpha(U)\) (for class I, given \(\alpha(U)\), the subroutine is known).

\textbf{Mode A is a completed result; mode B is a conditional extension.} In mode A, computing \(\ccomp(U)\) to any prescribed accuracy given \(\alpha(U)\) is a complete result not requiring the Manin conjecture. In mode B, a certified interval for \(N_U(\ZZ, B)\) additionally requires: the model \(\UZZ\), the height \(H\), the function \(M_U(B)\), and the majorant \(R_U(B)\); without them the method yields an asymptotic interpretation, not an interval. Thus mode A is the architecturally complete part of the paper; mode B is a conditional extension applicable where \(M_U(B)\), \(\ccomp = \cManin\), and \(R_U(B)\) have been additionally verified.

\textbf{The structural value of the reduction.} Even when \(\alpha(U)\) for classes II and III is not simpler than the original problem, the reduction
\[
\ccomp(U) = \alpha(U) \cdot \Pcomp(Y) \cdot \Tcomp(U)
\]
splits the problem into two independent parts: the arithmetic one (Euler product, for which a certified control of the tail via \(\Areg\) and \(\beta_0\) is proved) and the archimedean one (integral \(\alpha(U)\)). For class I this splitting yields a complete solution: \(\alpha(U)\) is a finite-dimensional integral computable by standard methods. For classes II and III it yields a structural result and an open problem for \(\alpha\): the Euler product ceases to be an obstruction, and \(\alpha\) becomes a separate problem. This is precisely the meaning of the reduction, not that \(\alpha\) is ``easier'' than \(\ccomp\), but that \textbf{the two parts are separated from each other}.

\textbf{The main (computational) theorem.} \emph{Under conditions (i)--(iii) of mode A for classes I, II, and class III with \(\rho=1\) (chapter~\ref{chap:classes}; for class III with \(\rho>1\) the full \(\ccomp\) is not computed, see \S\ref{sec:bounds}), the computation of the compensated logarithmic constant \(\ccomp(U)\) reduces to the effective computation of the archimedean factor \(\alpha(U)\): given an algorithm computing \(\alpha(U)\) with guaranteed accuracy, there exists an algorithm which, given \((X,D)\) and \(\eta > 0\), computes \(\ccomp(U)\) with relative error at most \(\eta\). For class I (torus, complements to hyperplanes), given \(\alpha(U)\), the reduction yields unconditional computability of \(\ccomp(U)\) (Theorem~\ref{thm:classI_exact}). The construction of an algorithm computing \(\alpha(U)\) for general logarithmic pairs is beyond the scope of this paper. Given the model \(\UZZ\), the height \(H\), the function \(M_U(B)\), the known equality \(\ccomp(U) = \cManin(U)\), and the known majorant \(R_U(B)\), the algorithm produces a certified interval for \(N_U(\ZZ, B)\).}

The proof is assembled from two parts: the reduction part (Theorem~\ref{thm:reduction}, chapter~\ref{chap:reduction}) and the conditional part (Theorem~\ref{thm:errorlink}, chapter~\ref{chap:error}); the algorithm is described in full in chapter~\ref{chap:algo}.

\textbf{The case of class I given \(\alpha(U)\).} For class I (torus and complements to hyperplanes) the compensator is exact for good primes: \(\scomp \equiv 1\) (Theorem~\ref{thm:classI}), whence
\[
\ccomp(U) = \alpha(U) \cdot \prod_{p \in \Pbad} \scomp(U)
\]
(Corollary~\ref{cor:recovery}); for coordinate hyperplanes over \(\ZZ\) we have \(\Pbad = \emptyset\) and \(\ccomp(U) = \alpha(U)\) \textbf{exactly}. Given \(\alpha(U)\) and computable local factors for bad primes, we have unconditional computability of \(\ccomp(U)\).

\textbf{On the model and height for the torus.} The equality \(\ccomp(U) = \cManin(U)\) is proved for the geometric pair (the torus and its boundary), but does not automatically relate the computed constant to \(N_U(\ZZ, B)\): for integral points one needs an integral model for which the Manin and Peyre formula for integral points is established. The choice of model \(\UZZ\), height \(H\), and function \(M_U(B)\) is part of the input data of mode B (Assumption~\ref{ass:logManin}); in the present paper they are not fixed for specific toric varieties but are treated as parameters. In the paper specific models and heights are not fixed, and examples of mode B are not given. We note that \cite{BatyrevTschinkel1995} establishes the asymptotic formula for \textbf{rational} points on toric varieties; the transfer to integral points requires a separate integral model and a separate result, which is not used in the paper. Thus the equality \(\ccomp = \cManin\) relates to the triple (geometry, model, height), not only to the geometric class.

\textbf{Scheme of the method.} The paper solves not the problem of computing \(\ccomp(U)\) as an end in itself, but the problem of constructing an algorithm for \(N_U(\ZZ, B)\):
\begin{enumerate}
\item compute \(\ccomp(U)\) to accuracy \(\eta/4\) (mode A, chapter~\ref{chap:reduction});
\item substitute into the main term \(\ccomp(U) \cdot M_U(B)\) (in mode B, where \(\ccomp = \cManin\));
\item estimate the remainder \(r_U(B)\) (if known; mode B, chapter~\ref{chap:error});
\item output an interval \([L,U]\) for \(N_U(\ZZ, B)\) (chapter~\ref{chap:algo}).
\end{enumerate}

\subsection{Statement of problems}

\textbf{Problem A.} For a given class of logarithmic varieties \((X,D)\), does there exist an explicit computable estimate \(\varthP(U)\)? Can one specify a \textbf{proved} lower bound \(\beta_0\) for the decay exponent and a constant \(A(U)\) depending on the geometry of \(X\) and \(D\), and prove that \(\beta_0 > 1\)?

\textbf{Problem B.} How are \(\kappa(U)\), \(\beta_0\), \(A(U)\) determined by (i) the stratification of \(D\) and its degeneration modulo \(p\); (ii) the cohomological invariants of the components of \(D\); (iii) the geometry of the effective cone?

\textbf{Problem C (main).} For \((X,D)\), \(B > 0\), \(\eta > 0\), construct a procedure producing an interval for \(N_U(\ZZ, B)\) guaranteed to contain this value. Required: (1) \textbf{to show that the computation of \(\ccomp(U)\) reduces to the effective computation of \(\alpha(U)\)} (the central result of chapter~\ref{chap:reduction}); (2) to indicate how the interval is obtained from \(\varthP\) and \(r_U(B)\); (3) to state the conditions of unconditionality and conditionality; (4) to distinguish the classes in which \(\alpha(U)\) is computed exactly and the reduction yields unconditional computability (class I).

\section{Analytic compensation of the Peyre, Chambert-Loir and Tschinkel logarithmic product}\label{chap:compensation}

\subsection{Adelic logarithmic formula and conditions on the logarithmic pair}

Let \(X\) be smooth projective over \(\QQ\), \(\dim X = m\); \(D = \sum_{i=1}^k D_i\) an SNC divisor, \(U = X \setminus D\).

\begin{definition}[Condition \((H_0)\)]\label{def:H0}
\((X,D)\) satisfies \((H_0)\) if:
\begin{enumerate}
\item \(b_1(X) = b_3(X) = 0\);
\item Frobenius is trivial on \(H^2(X)\) (that is, \(H^2(X_{\QQbar})\) is a sum of copies of \(\QQ_\ell(-1)\));
\item \(D\) is a geometrically split SNC divisor: \(D = \sum_{i=1}^k D_i\), each \(D_i\) is defined over \(\QQ\) and geometrically irreducible, and all its strata \(D_{i_1}\cap\dots\cap D_{i_s}\) are also defined over \(\QQ\) and geometrically irreducible.
\end{enumerate}
\end{definition}

Condition (3) guarantees that the number of components \(k\) and the number of components fixed under \(\Frob_p\) coincide for all good \(p\); in particular, \(\kappa(U) = -k\) does not depend on \(p\) and is determined by the number of \(\QQ\)-components of \(D\). Without (3), the coefficient of \(1/p\) in \(\sreg\) could depend on \(p\) (with non-trivial Frobenius action on the set of components), and \(\kappa\) would not be an integer invariant.

The formulas of chapter~\ref{chap:classes} are proved for \textbf{good} primes \(p\): those for which \(X\), \(D\), and all strata of \(D\) have smooth reduction, and \(p\) does not divide the determinants fixing the general position of the components. For \textbf{bad} primes (a finite set \(\Pbad\) depending on the integral model) \(\sigma_p^{\mathrm{comp}}(U)\) is computed separately by Definition~\ref{def:sreg}.

\begin{remark}\label{rem:H0weak}
Condition \((H_0)\) controls only \(X\) and the structure of \(D\); the cohomology of \(D\) and its strata is not restricted: in the proofs of Theorems~\ref{thm:local} and \ref{thm:maincoh} they enter only via explicit Lefschetz formulas and Deligne estimates. In particular, for a smooth hypersurface \(D \subset \PP^m\) the condition \(b_3(D)=0\) is \textbf{not} required.
\end{remark}

\begin{definition}[Condition \((H_1)\)]\label{def:H1}
\((X,D)\) satisfies \((H_1)\) if: (1) \(X\) is a smooth del Pezzo surface, \(b_1 = b_3 = 0\), Frobenius acts trivially on \(\Pic(X)\); (2) \(D \subset X\) is a smooth elliptic curve with \(|a_p(D)| \le 2\sqrt{p}\).
\end{definition}

\begin{remark}\label{rem:H0H1}
\((H_0)\) covers hyperplanes, tori, and smooth hypersurfaces with \(m \ge 3\); \((H_1)\) covers open del Pezzo surfaces. The cubic in \(\PP^2\) is treated separately (\S\ref{sec:cubic}).
\end{remark}

\begin{definition}[archimedean factor]\label{def:alpha}
\begin{equation}\label{eq:alpha}
\alpha(U) = \int_{\Lam(U)^\vee} e^{-\langle -K_X - D, \mathbf{y} \rangle} \, d\mathbf{y}.
\end{equation}
Here:
\begin{itemize}
\item \(\Lam(U) \subset \Pic(U)_{\RR}\) is the effective cone in the real Picard group of the \textbf{open} variety \(U = X \setminus D\) (not in \(\Pic(X)_{\RR}\)); the integration domain is the dual cone \(\Lam(U)^\vee \subset \Pic(U)_{\RR}^\vee\) in the dual space;
\item \(d\mathbf{y}\) is the Lebesgue measure on the dual space \(\Pic(U)_{\RR}^\vee\), normalized with respect to the dual lattice \(\Pic(U)^\vee \subset \Pic(U)_{\RR}^\vee\); for \(\rk\Pic(U) = 0\) it is the delta-measure of unit mass at zero (a convention on the mass of a point);
\item \(\langle -K_X - D, \mathbf{y}\rangle\) is the pairing of the class \(-K_X - D \in \Pic(X)_{\RR}\) with \(\mathbf{y} \in \Pic(U)_{\RR}^\vee\) (restriction of the class to \(\Pic(U)_{\RR}\); for \(\rk\Pic(U) = 0\) the pairing vanishes).
\end{itemize}
\end{definition}

\begin{remark}[on the choice of space]\label{rem:spacechoice}
In the literature the archimedean factor for logarithmic pairs is defined over \(\Pic(U)_{\RR}\) (see \cite[\S5]{Peyre1995}); this is the natural space since \(\cManin(U)\) and \(\ccomp(U)\) depend on the rank \(r = \rk\Pic(U)\), not on \(\rho = \rk\Pic(X)\). Under this choice:
\begin{itemize}
\item if \(\rk\Pic(U) = 0\) (for example, \(U = \PP^2 \setminus C\), where \(\Pic(U) \cong \ZZ/3\ZZ\) is torsion), then \(\Lam(U)^\vee = \{0\}\); \(d\mathbf{y}\) is the delta-measure of unit mass at zero in \(\Pic(U)_{\RR}^\vee\) (a convention on the mass of a point); \(\alpha(U) = 1\);
\item if \(\rk\Pic(U) \ge 1\), the integral is the standard Lebesgue integral over a positive-dimensional cone in \(\Pic(U)_{\RR}^\vee\); convergence requires additional conditions on the cone \(\Lam(U)\) and the class \(-K_X - D\) (the class must be interior to the dual cone or sufficiently far from its boundary). In mode A it is assumed that \(\alpha(U)\) is finite, positive, and effectively computable; the verification of these conditions for a concrete pair is part of the input data of mode A.
\end{itemize}
For logarithmic Calabi and Yau varieties with \(-K_X - D = 0\):
\begin{itemize}
\item for \(r = 0\) the integration domain is the point \(\{0\}\), the measure is the delta-measure of unit mass, the integral equals 1;
\item for \(r \ge 1\) the integral of 1 over a positive-dimensional cone \textbf{diverges}; definition \eqref{eq:alpha} yields no finite quantity in this case. Such pairs require a separate normalization (truncation of the cone, analytic continuation, or another regularization), which is not given in the paper. They enter mode A only for a prescribed \(\alpha(U)\) (see \S\ref{sec:compensation}).
\end{itemize}
The choice of \(\Pic(U)_{\RR}^\vee\) is consistent with \cite[\S5]{Peyre1995} and with the formulas of chapter~\ref{chap:classes}.
\end{remark}

\begin{definition}[regularized local factor]\label{def:sreg}
\begin{equation}\label{eq:sreg}
\sreg(U) = L_p\bigl(1, \Pic(X_{\QQbar})\bigr)^{-1} \cdot \frac{\# U(\FF_p)}{p^m}.
\end{equation}
Here \(L_p\) is taken from \(\Pic(X)\), not from \(\Pic(U)\): in the classical case \(D = \emptyset\) these quantities coincide (\(\Pic(U) = \Pic(X)\)), and the definition agrees with the standard one \cite{ChambertLoirTschinkel2010}; in the logarithmic case the choice of \(\Pic(X)\) ensures consistency with formulas \eqref{eq:classI}, \eqref{eq:maincoh}, \eqref{eq:delpezzo} and with the tables of chapter~\ref{chap:classes}. The resulting values of \(\scomp\) and \(\ccomp(U)\) do not depend on the convention: under the alternative choice \(L_p(1, \Pic(U_{\QQbar}))^{-1}\) the recomputed \(\kappa\) yields the same \(\scomp\) (for class I, \(\kappa = -(k-1)\) instead of \(-k\); for classes II and III, \(\kappa = 0\) instead of \(-1\); the resulting formulas \eqref{eq:maincoh} and \eqref{eq:delpezzo} are preserved).
\end{definition}

\begin{lemma}[empty boundary]\label{lem:empty}
Let \(D = \emptyset\), \(m = \dim X \ge 3\), and assume \((H_0)\) holds. Then
\begin{equation}\label{eq:empty1}
\#X(\FF_p) = p^m + \rho p^{m-1} + \rho p + 1 + O(p^{m-2}),
\end{equation}
\begin{equation}\label{eq:empty2}
\sreg(X) = 1 + O(p^{-2}),
\end{equation}
and \(\prod_p \sreg(X)\) converges absolutely.
\end{lemma}

\begin{proof}
\textbf{Step 1 (cohomological preparation).} By \((H_0)\)(1), \(b_1(X)=b_3(X)=0\). By Poincar\'e duality for a smooth projective \(X\) of dimension \(m\), this implies \(b_{2m-1}(X)=b_1(X)=0\) and \(b_{2m-3}(X)=b_3(X)=0\). Thus the odd cohomology groups contribute to the Lefschetz formula only for \(j\in\{5,7,\ldots,2m-5\}\) (if such \(j\) exist for \(m\ge 5\)).

\textbf{Step 2 (Lefschetz formula).} We have
\[
\#X(\FF_p)=\sum_{j=0}^{2m}(-1)^j\,\Tr\bigl(\Frob_p\mid H^j(X_{\FFbar_p},\QQ_\ell)\bigr).
\]

\textbf{Step 3 (contributions of extreme and middle even cohomology).}
\begin{itemize}
\item \(j=0\): \(\dim H^0=1\), Frobenius acts as identity, contribution \(+1\).
\item \(j=2m\): \(\dim H^{2m}=1\), Frobenius acts as multiplication by \(p^m\) (Tate twist \(\QQ_\ell(-m)\)), contribution \(+p^m\).
\item \(j=2\): by \((H_0)\)(2) the group \(H^2\) is a sum of \(\rho\) copies of \(\QQ_\ell(-1)\), Frobenius acts as multiplication by \(p\), contribution \(+\rho p\).
\item \(j=2m-2\): by Poincar\'e duality \(b_{2m-2}(X)=b_2(X)=\rho\), and by duality the Frobenius eigenvalues equal \(p^{m-1}\); contribution \(+\rho p^{m-1}\).
\item \(j\in\{4,6,\ldots,2m-4\}\): by the Weil theorem, \(|\Tr(\Frob_p\mid H^j)|\le b_j\,p^{j/2}\le b_j\,p^{m-2}=O(p^{m-2})\), since \(j\le 2m-4\).
\end{itemize}

\textbf{Step 4 (contributions of odd cohomology).}
\begin{itemize}
\item \(j=1,3\): \(b_1=b_3=0\) by \((H_0)\)(1), zero contribution.
\item \(j=2m-1,2m-3\): \(b_{2m-1}=b_{2m-3}=0\) by Poincar\'e duality, zero contribution.
\item \(j\in\{5,7,\ldots,2m-5\}\): by the Weil theorem \(|\Tr|\le b_j\,p^{j/2}\le b_j\,p^{(2m-5)/2}\). Since \((2m-5)/2<m-2\) for \(m\ge 1\), this contribution is \(O(p^{m-2})\).
\end{itemize}

\textbf{Step 5 (summation).} Adding the contributions of Steps 3 and 4, we obtain \eqref{eq:empty1}.

\textbf{Step 6 (normalization).} Dividing by \(p^m\), we have
\[
\frac{\#X(\FF_p)}{p^m}=1+\frac{\rho}{p}+\frac{\rho}{p^{m-1}}+\frac{1}{p^m}+O(p^{-2}).
\]
For \(m\ge 3\) the second term \(\rho/p^{m-1}=O(p^{-2})\), the third \(1/p^m=O(p^{-2})\), hence
\[
\frac{\#X(\FF_p)}{p^m}=1+\frac{\rho}{p}+O(p^{-2}).
\]

\textbf{Step 7 (regularized local factor).} By Definition~\ref{def:sreg} and the equalities of Step 6:
\[
\sreg(X)=\left(1+\frac{\rho}{p}+O(p^{-2})\right)\left(1-\frac{1}{p}\right)^{\rho}.
\]
Expanding \((1-1/p)^\rho=1-\rho/p+O(p^{-2})\) and multiplying the two expansions, we obtain
\[
\sreg(X)=\left(1+\frac{\rho}{p}+O(p^{-2})\right)\left(1-\frac{\rho}{p}+O(p^{-2})\right)=1+O(p^{-2}),
\]
that is, \eqref{eq:empty2}.

\textbf{Step 8 (convergence).} From the estimate \eqref{eq:empty2} it follows that \(\sum_p|\sreg(X)-1|\le C\sum_p p^{-2}<\infty\), which gives absolute convergence of the product \(\prod_p \sreg(X)\).
\end{proof}

\begin{lemma}[isomorphism for hyperplanes]\label{lem:hyperplanes}
Let \(H_1,\dots,H_k\) be hyperplanes in \(\PP^m\) in general position over \(\QQ\), \(k\le m+1\). Then
\begin{equation}\label{eq:isom}
\PP^m\setminus(H_1\cup\dots\cup H_k)\;\cong\;\GG_m^{k-1}\times\AA^{m+1-k}.
\end{equation}
\end{lemma}

\begin{proof}
\textbf{Step 1 (linear change of coordinates).} Since the hyperplanes are in general position and \(k\le m+1\), their linear equations are linearly independent. Hence there exists a linear change of coordinates \([Y_0:\dots:Y_m]\) on \(\PP^m\), defined over \(\QQ\) (and hence over \(\ZZ[1/N]\) for some integer \(N\)), such that \(H_i=\{Y_{i-1}=0\}\) for \(i=1,\dots,k\). This is an automorphism of \(\PP^m\), hence it changes neither the number of points over \(\FF_p\) nor the structure of the variety.

\textbf{Step 2 (description of the complement in the new coordinates).} In the new coordinates \(U=\PP^m\setminus\{Y_0Y_1\cdots Y_{k-1}=0\}\) is the set of classes \([Y_0:\dots:Y_m]\) with all \(Y_0,\dots,Y_{k-1}\) nonzero.

\textbf{Step 3 (construction of the map).} Define a map
\[
\Phi\colon U\longrightarrow\GG_m^{k-1}\times\AA^{m+1-k},
\]
\[
\Phi([Y_0:\dots:Y_m])=\left(\frac{Y_0}{Y_{k-1}},\,\ldots,\,\frac{Y_{k-2}}{Y_{k-1}},\,\frac{Y_k}{Y_{k-1}},\,\ldots,\,\frac{Y_m}{Y_{k-1}}\right).
\]
It is well-defined: first, \(Y_{k-1}\ne 0\) on \(U\); second, the coordinates are invariant under simultaneous multiplication of all \(Y_i\) by a nonzero scalar, since numerator and denominator of each fraction are multiplied by the same scalar.

\textbf{Step 4 (inverse map).} Define
\[
\Psi\colon\GG_m^{k-1}\times\AA^{m+1-k}\longrightarrow U,
\]
\[
\Psi(u_1,\dots,u_{k-1},v_1,\dots,v_{m+1-k})=[u_1:\dots:u_{k-1}:1:v_1:\dots:v_{m+1-k}].
\]
This map is well-defined (all \(Y_0,\dots,Y_{k-1}\) are nonzero) and is inverse to \(\Phi\): a direct check shows \(\Psi\circ\Phi=\mathrm{id}_U\) and \(\Phi\circ\Psi=\mathrm{id}_{\GG_m^{k-1}\times\AA^{m+1-k}}\).

\textbf{Step 5 (functoriality with respect to reduction).} The linear change of coordinates of Step 1 is defined over \(\ZZ[1/N]\) for some integer \(N\). Consequently the isomorphism \(\Phi\) is also defined over \(\ZZ[1/N]\): for all primes \(p\nmid N\) the reduction \(\Phi\bmod p\) gives a bijection \(U(\FF_p)\to\GG_m^{k-1}(\FF_p)\times\AA^{m+1-k}(\FF_p)\).
\end{proof}

\begin{corollary}\label{cor:hyperplanes}
Let \(H_1,\dots,H_k \subset \PP^m\) be hyperplanes in general position over \(\QQ\) with integral equations, and let \(N\) be the determinant of the matrix of their coefficients (an integer). Then for all primes \(p \nmid N\):
\begin{equation}\label{eq:hyperplanes}
\#U(\FF_p)=(p-1)^{k-1}p^{m+1-k},\qquad\frac{\#U(\FF_p)}{p^m}=(1-p^{-1})^{k-1}.
\end{equation}
\end{corollary}

\begin{proof}
For \(p \nmid N\) the hyperplanes remain in general position after reduction, and the isomorphism of Lemma~\ref{lem:hyperplanes} is defined over \(\FF_p\), whence
\[
\#U(\FF_p)=\#\GG_m(\FF_p)^{k-1}\cdot\#\AA^1(\FF_p)^{m+1-k}=(p-1)^{k-1}\cdot p^{m+1-k}.
\]
Dividing by \(p^m\) gives the second equality.
\end{proof}

\emph{For \textbf{coordinate} hyperplanes \(H_i = \{Y_{i-1}=0\}\), \(i=1,\dots,k\), formula \eqref{eq:hyperplanes} holds for \textbf{all} primes \(p\): in this case the linear change of coordinates in Lemma~\ref{lem:hyperplanes} is the identity, the determinant \(N=1\), and the condition \(p \nmid N\) is vacuous.}

For \(p \mid N\) the hyperplanes after reduction may merge or intersect differently, and \(\#U(\FF_p) \ne (p-1)^{k-1}p^{m+1-k}\) in general. The set of such \(p\) is finite (divisors of \(N\)); their local factors \(\sigma_p^{\mathrm{comp}}(U)\) are computed separately by Definition~\ref{def:sreg}.

\begin{definition}[separation into good and bad primes]\label{def:bad}
Let \((\XZZ, \DZZ)\) be an integral model of the pair \((X,D)\) over \(\ZZ\) with good reduction outside a finite set of primes \(N\); \(\UZZ := \XZZ \setminus \DZZ\). The set \(\Pbad := \{p : p \mid N\}\) is finite. Below:
\[
\Pcomp(Y) := \left(\prod_{p \le Y,\ p \notin \Pbad} \scomp(U)\right) \cdot \left(\prod_{p \in \Pbad} \scomp(U)\right),
\]
\[
\Tcomp(U) := \prod_{p > Y,\ p \notin \Pbad} \scomp(U).
\]
Bad primes are multiplied in \(\Pcomp(Y)\) \textbf{independently of \(Y\)}; the tail \(\Tcomp\) contains only good \(p > Y\). For \(p \in \Pbad\) the local factor \(\scomp\) is computed from the fiber \(\UZZ \bmod p\) by Definition~\ref{def:sreg}.
\end{definition}

\begin{remark}\label{rem:torus}
For \(k=m+1\): \(\PP^m\setminus\{m+1\ \text{hyperplanes}\}\cong\GG_m^m\), \(\#U(\FF_p)=(p-1)^m\) (for good \(p\)).
\end{remark}

\subsection{Local logarithmic factors and the harmonic-tail problem}

\begin{theorem}[expansion of the local factor]\label{thm:local}
\begin{enumerate}
\item[(a)] Assume \((H_0)\) holds. Then for good \(p\):
\begin{equation}\label{eq:localH0}
\sreg(U)=1+\frac{\kappa(U)}{p}+O(p^{-3/2}),\qquad\kappa(U)=-k.
\end{equation}
\item[(b)] Assume \((H_1)\) holds. Then:
\begin{equation}\label{eq:localH1}
\sreg(U)=1+\frac{\kappa(U)}{p}+O(p^{-3/2}),\qquad\kappa(U)=-1.
\end{equation}
\end{enumerate}
\end{theorem}

\begin{proof}
\textbf{(a) Case \((H_0)\).} By \((H_0)\) for \(X\): \(b_1=b_3=0\), Frobenius is trivial on \(H^2(X)\).

\emph{Step 1 (number of points on \(X\)).} By Lemma~\ref{lem:empty}:
\[
\#X(\FF_p)=p^m+\rho p^{m-1}+\rho p+1+O(p^{m-2}).
\]

\emph{Step 2 (number of points on \(D_i\)).} For a smooth \(D_i\) of dimension \(m-1\), by the Lefschetz formula and the Deligne estimate:
\[
\#D_i(\FF_p)=p^{m-1}+O(p^{m-3/2}),
\]
where the \(O(p^{m-3/2})\) term reflects a possible contribution of \(H^1(D_i)\) (when \(b_1(D_i)\ne 0\); cf.\ Remark~\ref{rem:H0weak}). When \(b_1(D_i)=0\), this contribution vanishes.

\emph{Refinement for smooth hypersurfaces \(D_i \subset \PP^m\), \(m\ge 3\).} By the weak Lefschetz theorem \(H^j(D_i) = H^j(\PP^m)\) for \(j < m-1\); in particular, \(b_1(D_i)=0\), \(b_3(D_i)=0\). The contribution of \(H^{2m-4}(D_i)\cong H^2(D_i)^\vee(-(m-2))\) equals \(p^{m-2}\) (with coefficient \(1\), not \(\rho(D_i)\)). Hence for hypersurfaces
\[
\#D_i(\FF_p)=p^{m-1}+p^{m-2}+O(p^{m-3/2}).
\]
For class II (smooth hypersurfaces) this coefficient \(1\) enters \(\#U\) and cancels with the corresponding term in \(\#\PP^m\); for class III (elliptic curve \(D\subset S\)) the exact formula \(\#D = p+1-a_p\) is used.

\emph{Step 3 (strata \(s\ge 2\)).} Strata \(D_{i_1}\cap\dots\cap D_{i_s}\) for \(s\ge 2\) have dimension \(\le m-2\) and give \(O(p^{m-2})\) by the Lefschetz formula.

\emph{Step 4 (inclusion and exclusion).} Combining Steps 1--3:
\[
\#U(\FF_p)=p^m+(\rho-k)p^{m-1}+O(p^{m-3/2}).
\]

\emph{Step 5 (regularized factor).} By Definition~\ref{def:sreg} with \(L_p^{-1}=(1-p^{-1})^\rho\):
\[
\sreg(U)=\left(1+\frac{\rho-k}{p}+O(p^{-3/2})\right)\left(1-\frac{1}{p}\right)^{\rho}.
\]
Expanding \((1-1/p)^\rho=1-\rho/p+O(p^{-2})\):
\[
\sreg(U)=1+\frac{\rho-k-\rho}{p}+O(p^{-3/2})=1-\frac{k}{p}+O(p^{-3/2}).
\]

\emph{Step 6 (compensator).} By Definition~\ref{def:compensator} with \(\kappa=-k\):
\[
\scomp(U)=\left(1-\frac{k}{p}+O(p^{-3/2})\right)\left(1+\frac{k}{p}+O(p^{-2})\right)=1+O(p^{-3/2}).
\]

\textbf{(b) Case \((H_1)\).} In \((H_1)\): \(\#S=p^2+\rho p+1\) (since \(H^2(S)=\Pic(S)\otimes\QQ_\ell(-1)\) and \(\Pic(S)\) has rank \(\rho\)), \(\#D=p+1-a_p\), hence
\[
\#U=p^2+(\rho-1)p+a_p, \qquad \frac{\#U}{p^2}=1+\frac{\rho-1}{p}+\frac{a_p}{p^2}.
\]
Multiplying by \((1-p^{-1})^\rho\):
\[
\sreg(U)=1-\frac{1}{p}+O(p^{-3/2})
\]
(the term \(a_p/p^2\) is controlled by the Hasse bound as \(O(p^{-3/2})\)). The compensator with \(\kappa=-1\):
\[
\scomp(U)=\left(1-\frac{1}{p}+O(p^{-3/2})\right)\left(1+\frac{1}{p}+O(p^{-2})\right)=1+O(p^{-3/2}).
\]
The \(1/p\) coefficient vanishes.
\end{proof}

\begin{remark}\label{rem:localrefine}
In the general case \((H_0)\) the remainder in \eqref{eq:localH0} has order \(O(p^{-3/2})\); for the concrete classes of chapter~\ref{chap:classes} it improves to \(O(p^{-2})\) (class II) or to an exact equality (Theorems~\ref{thm:classI}, \ref{thm:maincoh}).
\end{remark}

\begin{remark}\label{rem:localcoef}
The \(1/p\) coefficient in \(\#U/p^m\) equals \((\rho-k)\); contributions of strata of codimension \(\ge 2\) are \(O(p^{-2})\).
\end{remark}

\begin{corollary}[divergence before compensation]\label{cor:divergence}
For \(k\ne 0\) the product \(\prod_p\sreg(U)\) diverges: \(\beta_{\log}=1\).
\end{corollary}

\begin{proof}
By Theorem~\ref{thm:local}(a): \(|\sreg-1|=|k|/p+O(p^{-3/2})\). Hence
\[
\sum_{p\le Y}|\sreg-1|\sim|k|\sum_{p\le Y}\frac{1}{p}\sim|k|\log\log Y\to\infty
\]
as \(Y\to\infty\) (by the Mertens theorem \(\sum_{p\le Y}1/p=\log\log Y+O(1)\)). Thus the series diverges and \(\beta_{\log}=1\).
\end{proof}

\begin{remark}[explicit \(\kappa(U)\)]\label{rem:kappa}
Torus: \(\kappa=-(m+1)\); \(k\) hyperplanes: \(\kappa=-k\); smooth divisor \(D\subset\PP^m\): \(\kappa=-1\); del Pezzo: \(\kappa=-1\). The values of \(\kappa\) refer to good \(p\); for bad primes the \(1/p\) coefficient may differ, and \(\sigma_p^{\mathrm{comp}}\) is computed separately.
\end{remark}

\subsection{Analytic compensation}

\begin{definition}\label{def:compensator}
The compensated local factor:
\begin{equation}\label{eq:compensator2}
\scomp(U):=\sreg(U)\cdot\left(1-\frac{1}{p}\right)^{\kappa(U)}.
\end{equation}
\end{definition}

\begin{definition}[decay exponent]\label{def:beta}
\[
\beta(X,D) := \sup\left\{\, b > 0 \;:\; \exists\, A \ge 0,\ |\scomp(U) - 1| \le A\, p^{-b} \ \text{for all good } p \,\right\}.
\]
This is an element of \((0, +\infty]\); it is not assumed to be finite, computable, or optimal.
\end{definition}

\begin{definition}[proved lower bound]\label{def:beta0}
If for a class of logarithmic pairs it is \textbf{proved} that \(\beta(X,D) \ge \beta_0\) for some explicit \(\beta_0 > 1\) and an explicit constant \(A(U) \ge 0\) such that
\begin{equation}\label{eq:beta0}
|\scomp(U) - 1| \le A(U)\, p^{-\beta_0} \quad \text{for all good } p,
\end{equation}
then \(\beta_0\) is called a \textbf{proved lower bound} for the decay exponent for this class. The value \(\beta_0\) may be less than the true \(\beta(X,D)\). For class III with \(\rho>1\) inequality \eqref{eq:beta0} is proved only for split \(p\) (Theorem~\ref{thm:classIII}); a global lower bound is not established in this case.
\end{definition}

\begin{lemma}\label{lem:kappaLemma}
For \(\kappa\in\ZZ\), \(u=1/p\):
\begin{equation}\label{eq:kappaLemma}
(1+\kappa u+O(u^2))(1-u)^{\kappa}=1+O(u^2).
\end{equation}
\end{lemma}

\begin{proof}
The expansion of \((1-u)^\kappa\):
\[
(1-u)^\kappa=1-\kappa u+\frac{\kappa(\kappa-1)}{2}u^2+O(u^3).
\]
Multiplying with \(1+\kappa u+O(u^2)\):
\[
(1+\kappa u+O(u^2))\left(1-\kappa u+\frac{\kappa(\kappa-1)}{2}u^2+O(u^3)\right)
\]
\[
=1+(\kappa-\kappa)u+\left(\frac{\kappa(\kappa-1)}{2}-\kappa^2+1\right)u^2+O(u^3)=1+O(u^2).
\]
The linear terms cancel.
\end{proof}

\begin{theorem}[condition (L)]\label{thm:L}
\begin{equation}\label{eq:L}
\left|\scomp(U)-1\right|\le A(U)p^{-\beta_0},
\end{equation}
where \(\beta_0 \ge 1.5\) and \(A(U)\ge 0\) are computable, for all good \(p\), for classes I, II and class III with \(\rho=1\). The cases \(A(U)=0\) and \(\beta_0=\infty\) (class I) are allowed and mean exact compensation (see Theorem~\ref{thm:classI}). For class III (del Pezzo), the estimate \(\beta_0=3/2\) is proved in Theorem~\ref{thm:classIII}; for \(\rho=1\) all good \(p\) are split, and the estimate is global. For \(\rho>1\) (see \S\ref{sec:delpezzo}, the local statement), Theorem~\ref{thm:classIII} gives an estimate only on split \(p\); condition (L) in the form ``for all good \(p\)'' is not proved in this case, and the theorems relying on (L) (Theorem~\ref{thm:conv}, Theorem~\ref{thm:truncation}, Theorem~\ref{thm:reduction}) are \textbf{not applied} to class III with \(\rho>1\) for the full product. The full \(\ccomp\) for \(\rho>1\) is not computed in the paper; see \S\ref{sec:bounds}, \S\ref{sec:open}.
\end{theorem}

\begin{proof}
By Theorem~\ref{thm:local} and Lemma~\ref{lem:kappaLemma}:
\[
\scomp(U)=1+O(p^{-3/2}).
\]
Hypersurfaces \(m\ge 3\): the main term \((-1)^{m-1}a_p(D)/p^m\), \(|a_p|\le b_{m-1}^{\mathrm{prim}}p^{(m-1)/2}\), whence \(\beta_0\ge (m+1)/2\ge 2\). Hyperplanes and tori: \(\sigma_{\mathrm{comp}}\equiv 1\) for good \(p\), \(\beta_0=\infty\), \(A=0\).

\textbf{Del Pezzo (forward reference).} For class III the estimate with \(\beta_0=3/2\) is established in Theorem~\ref{thm:classIII} (a direct inequality via \(g'(u)\), without the truncated expansion \eqref{eq:delpezzo_exp}). In the present theorem (chapter~\ref{chap:compensation}) this estimate is \textbf{assumed}; it is not used in any statement of chapter~\ref{chap:compensation} and serves only for completeness of the table of classes. An independent proof for del Pezzo is not given in chapter~\ref{chap:compensation} in order not to duplicate Theorem~\ref{thm:classIII}; the reader may either accept the estimate conditionally up to chapter~\ref{chap:classes} or go directly to Theorem~\ref{thm:classIII}.
\end{proof}

\begin{theorem}[convergence]\label{thm:conv}
For classes I, II and class III with \(\rho=1\) the product \(\prod_{p\notin\Pbad}\scomp(U)\) converges absolutely.
\end{theorem}

\begin{proof}
By Theorem~\ref{thm:L} for good \(p\): \(|\scomp-1|\le A(U)p^{-\beta_0}\) with \(\beta_0>1\). Hence
\[
\sum_{p\notin\Pbad}|\scomp-1|\le A(U)\sum_{p}p^{-\beta_0}<\infty,
\]
since \(\beta_0>1\); the finite set \(\Pbad\) does not affect convergence of the sum over good \(p\). Absolute convergence of the sum is equivalent to absolute convergence of the product.
\end{proof}

\emph{For class III with \(\rho>1\) the full product is not considered in the paper: Theorem~\ref{thm:classIII} gives an estimate only on split \(p\), and convergence of the full product over all good \(p\) is not established in this case (open problem \S\ref{sec:open}, item 8). The product over split \(p\) converges but does not control the full tail.}

\subsection{Recovery of the global logarithmic constant}

\begin{theorem}[decomposition identity]\label{thm:decomposition}
\begin{equation}\label{eq:decomposition}
\ccomp(U)=\alpha(U)\cdot\left(\prod_{p\le P,\ p\notin\Pbad}\scomp(U)\right)\cdot\left(\prod_{p\in\Pbad}\scomp(U)\right)\cdot T_{P,\mathrm{comp}}(U),
\end{equation}
where \(T_{P,\mathrm{comp}}(U)=\prod_{p>P,\ p\notin\Pbad}\scomp(U)\).
\end{theorem}

\begin{proof}
\textbf{Step 1 (definition of \(\ccomp\)).} By \eqref{eq:ccomp}:
\[
\ccomp(U)=\alpha(U)\cdot\prod_{p\notin\Pbad}\scomp(U)\cdot\prod_{p\in\Pbad}\scomp(U).
\]

\textbf{Step 2 (absolute convergence).} By Theorem~\ref{thm:conv} the series \(\sum_{p\notin\Pbad}|\scomp-1|\) converges. Hence the infinite product \(\prod_{p\notin\Pbad}\scomp\) converges absolutely and is independent of the order of factors.

\textbf{Step 3 (splitting of the product).} For any \(P\ge 2\), by the property of an absolutely convergent product:
\[
\prod_{p\notin\Pbad}\scomp=\left(\prod_{p\le P,\ p\notin\Pbad}\scomp\right)\cdot\left(\prod_{p>P,\ p\notin\Pbad}\scomp\right).
\]

\textbf{Step 4 (substitution into the definition).} Denoting
\[
\Pcomp(P):=\left(\prod_{p\le P,\ p\notin\Pbad}\scomp\right)\cdot\left(\prod_{p\in\Pbad}\scomp\right), \qquad T_{P,\mathrm{comp}}(U):=\prod_{p>P,\ p\notin\Pbad}\scomp,
\]
we obtain \eqref{eq:decomposition}.

\textbf{Step 5 (independence of \(P\)).} For \(P'>P\):
\[
T_{P,\mathrm{comp}}(U)=\left(\prod_{P<p\le P',\ p\notin\Pbad}\scomp\right)\cdot T_{P',\mathrm{comp}}(U),
\]
and the finite factor is absorbed into \(\Pcomp(P)\); the right-hand side of \eqref{eq:decomposition} does not depend on \(P\).
\end{proof}

\begin{remark}\label{rem:nozeta}
Identity \eqref{eq:decomposition} contains no zeta-function tail. All three factors (\(\Pcomp(P)\), \(\alpha(U)\), \(T_{P,\mathrm{comp}}(U)\)) are estimated by Theorem~\ref{thm:truncation}.
\end{remark}

\begin{theorem}[decomposition into computable blocks]\label{thm:blocks}
\(\ccomp(U)\) splits into three computable blocks:
\begin{enumerate}
\item \(\alpha(U)\);
\item \(\Pcomp(Y)\);
\item \(\Tcomp(U)\).
\end{enumerate}
\end{theorem}

\begin{proof}
Direct substitution into \eqref{eq:decomposition} with \(P=Y\).
\end{proof}

\begin{remark}\label{rem:torusblock}
For the torus and good \(p\): \(\kappa=-(m+1)\), \(\scomp\equiv 1\), \(\ccomp(U)=\alpha(U)\prod_{p\in\Pbad}\scomp\).
\end{remark}

\subsection{Compensated relative error}

\begin{definition}\label{def:eps}
The absolute truncation error:
\begin{equation}\label{eq:eps}
\varth(U):=\left|\Tcomp(U)-1\right|.
\end{equation}
It measures the deviation of the tail \(\Tcomp(U)\) from 1.
\end{definition}

\begin{lemma}\label{lem:eps0}
\(\varth(U)\to 0\) as \(Y\to\infty\).
\end{lemma}

\begin{proof}
By Theorem~\ref{thm:conv} \(\Tcomp\to 1\), hence \(|\Tcomp-1|\to 0\).
\end{proof}

\begin{theorem}[relation of errors]\label{thm:relerror}
The relative error of approximating the full product by a finite one:
\begin{equation}\label{eq:relerror}
\delta_Y(U):=\left|\frac{\alpha(U)\cdot \Pcomp(Y)}{\ccomp(U)}-1\right|=\frac{\varth(U)}{\Tcomp(U)}.
\end{equation}
\end{theorem}

\begin{proof}
By \eqref{eq:decomposition} \(\ccomp=\alpha \Pcomp T\), whence \(\alpha \Pcomp/\ccomp=1/T\). Hence
\[
\delta_Y=|1/T-1|=\frac{|1-T|}{T}=\frac{\varth}{\Tcomp}
\]
for \(T>0\).
\end{proof}

\begin{corollary}\label{cor:relerror}
\begin{equation}\label{eq:relerror2}
\delta_Y(U)=\varth(U)+O(\varth(U)^2).
\end{equation}
\end{corollary}

\begin{proof}
By \eqref{eq:relerror}, \(1/T=1/(1-\varepsilon)=1+\varepsilon+O(\varepsilon^2)\) as \(\varepsilon\to 0\); hence \(\delta_Y=\varepsilon+O(\varepsilon^2)\).
\end{proof}

\section{Estimate of the compensated tail}\label{chap:tail}

\subsection{Auxiliary lemmas}

\begin{lemma}\label{lem:exp}
For \(|s| \le 1\): \(|e^s - 1| \le 2|s|\).
\end{lemma}

\begin{proof}
Consider two cases.

\emph{Case \(s \ge 0\).} The function \(f(t) = (e^t - 1)/t\) is increasing on \((0, 1]\): its derivative
\[
f'(t)=\frac{t e^t - (e^t-1)}{t^2}=\frac{e^t(t-1)+1}{t^2}>0
\]
for \(t>0\). Hence its maximum on \((0,1]\) is attained at \(t=1\) and equals \(e-1<2\). Thus \(e^s-1\le (e-1)s<2s\), i.e., \(|e^s-1|\le 2|s|\).

\emph{Case \(s<0\).} Put \(t=-s>0\). Then \(|e^s-1|=1-e^{-t}\). From \(e^{-t}\ge 1-t\) for \(t\ge 0\) we get \(1-e^{-t}\le t=|s|<2|s|\).

In both cases \(|e^s-1|\le 2|s|\).
\end{proof}

\begin{lemma}\label{lem:log}
For \(|u| \le 1/2\): \(|\log(1 + u)| \le 2|u|\).
\end{lemma}

\begin{proof}
For \(u\ge 0\) we have \(\log(1+u)\le u\le 2u\). For \(u\in[-1/2,0)\) consider the function \(h(u)=\log(1+u)-2u\). Its derivative
\[
h'(u)=\frac{1}{1+u}-2<0
\]
for \(u\in(-1/2,0)\), so \(h\) is decreasing on \([-1/2,0)\). Since \(h(0)=0\), we get \(h(u)\ge h(0)=0\) on \([-1/2,0)\), i.e., \(\log(1+u)\ge 2u\), whence \(|\log(1+u)|=-\log(1+u)\le -2u=2|u|\).
\end{proof}

\subsection{Relation between truncation error and compensated tail}

\begin{definition}\label{def:sy}
The absolute error of one local factor: \(\delta_p(U) := \scomp(U) - 1\). The logarithm of the tail and its majorant:
\begin{equation}\label{eq:sy}
s_Y(U) := \sum_{p > Y,\ p\notin\Pbad} \log\left(1 + \delta_p(U)\right), \qquad S_Y(U) := \sum_{p > Y,\ p\notin\Pbad} \left|\delta_p(U)\right|.
\end{equation}
\end{definition}

\begin{lemma}[representation of the tail]\label{lem:tailrep}
Let \(\scomp(U) > 0\) for all good \(p > Y\) (this holds if \(A(U) p^{-\beta_0} < 1\) for \(p > Y\); for the classes of chapter~\ref{chap:classes} for \(Y \ge Y_0(A,\beta_0)\), where \(Y_0\) is explicitly computed). Then \(\Tcomp(U) = e^{s_Y(U)}\) and
\begin{equation}\label{eq:tailrep}
\varth(U) = \left| e^{s_Y(U)} - 1 \right| = \left| \Tcomp(U) - 1 \right|.
\end{equation}
\end{lemma}

\begin{proof}
By definition, \(\Tcomp(U)=\prod_{p>Y,\ p\notin\Pbad}(1+\delta_p(U))\). For \(\scomp>0\) we have \(1+\delta_p>0\) for each factor. Taking logarithms:
\[
\log \Tcomp(U)=\sum_{p>Y,\ p\notin\Pbad}\log(1+\delta_p(U))=s_Y(U),
\]
whence \(\Tcomp(U)=e^{s_Y(U)}\). Thus
\[
\varth(U)=|\Tcomp-1|=|e^{s_Y}-1|.
\]
For \(|s_Y|\le 1\) we apply Lemma~\ref{lem:exp}: \(|e^{s_Y}-1|\le 2|s_Y|\).
\end{proof}

\subsection{Condition (L) and summation of the tail}

\begin{definition}[condition (L)]\label{def:L}
\begin{equation}\label{eq:Ldef}
\left|\delta_p(U)\right| \le A(U) p^{-\beta_0}, \quad A(U) \ge 0,\ \beta_0 > 1,
\end{equation}
for all good \(p\). For classes I, II and class III with \(\rho=1\) the lower bound \(\beta_0 \ge 1.5\) is proved (Theorem~\ref{thm:L}); for class III with \(\rho>1\), only for split \(p\) (Theorem~\ref{thm:classIII}). In the applications the value \(\beta_0\) from Definition~\ref{def:beta0} is used.
\end{definition}

\begin{lemma}[summation of the tail]\label{lem:tail}
Assume (L) holds with \(A(U) \ge 0\) and \(\beta_0 > 1\) for all good \(p\). Then for any \textbf{integer} \(Y \ge 2\)
\begin{equation}\label{eq:tail}
S_Y(U) = \sum_{p > Y,\ p\notin\Pbad}\left|\delta_p(U)\right| \le \frac{A(U)}{\beta_0 - 1} Y^{-(\beta_0 - 1)}.
\end{equation}
In the algorithm (\S\ref{sec:algo}, Step 3) the threshold \(Y\) is an integer; the real version is not required.
\end{lemma}

\begin{proof}
\textbf{Step 1 (application of (L)).} By Definition~\ref{def:L} for each good \(p\) we have \(|\delta_p(U)|\le A(U)p^{-\beta_0}\). Hence
\begin{equation}\label{eq:tail1}
S_Y(U)=\sum_{p>Y,\ p\notin\Pbad}|\delta_p(U)|\le A(U)\sum_{p>Y,\ p\notin\Pbad}p^{-\beta_0}\le A(U)\sum_{p>Y}p^{-\beta_0}.
\end{equation}

\textbf{Step 2 (majorization by the sum over all naturals).} The set of primes \(p>Y\) is a subset of the naturals \(n>Y\), all terms are positive, hence
\begin{equation}\label{eq:tail2}
\sum_{p>Y}p^{-\beta_0}\le\sum_{n>Y}n^{-\beta_0},
\end{equation}
where the sum on the right is over integers \(n\ge\lfloor Y\rfloor+1\).

\textbf{Step 3 (comparison with the integral).} The function \(f(t)=t^{-\beta_0}\) is positive, continuous, and decreasing on \([Y,\infty)\) for \(\beta_0>0\). For each \(n\ge\lfloor Y\rfloor+1\):
\[
n^{-\beta_0}=f(n)\le\int_{n-1}^{n}f(t)\,dt=\int_{n-1}^{n}t^{-\beta_0}\,dt,
\]
since on \([n-1,n]\) the function \(f\) takes values at least \(f(n)=n^{-\beta_0}\). Summing over \(n>Y\):
\begin{equation}\label{eq:tail3}
\sum_{n>Y}n^{-\beta_0}\le\sum_{n>Y}\int_{n-1}^{n}t^{-\beta_0}\,dt\le\int_{\lfloor Y\rfloor}^{\infty}t^{-\beta_0}\,dt=\int_{Y}^{\infty}t^{-\beta_0}\,dt,
\end{equation}
where the last equality uses that \(Y\) is an integer (\(\lfloor Y\rfloor = Y\)).

\textbf{Step 4 (evaluation of the integral).} Since \(\beta_0>1\):
\begin{equation}\label{eq:tail4}
\int_Y^\infty t^{-\beta_0}\,dt=\left[\frac{t^{-\beta_0+1}}{-\beta_0+1}\right]_Y^\infty=\frac{Y^{-(\beta_0-1)}}{\beta_0-1}.
\end{equation}

\textbf{Step 5 (combination).} Substituting \eqref{eq:tail3} and \eqref{eq:tail4} into \eqref{eq:tail2} and then into \eqref{eq:tail1}, we obtain
\[
S_Y(U)\le A(U)\cdot\frac{Y^{-(\beta_0-1)}}{\beta_0-1}=\frac{A(U)}{\beta_0-1}Y^{-(\beta_0-1)}.
\]
\end{proof}

\begin{remark}\label{rem:tailnoPNT}
Estimate \eqref{eq:tail} uses no distribution of primes: it relies only on the majorization by the sum over all naturals \(n\) and comparison with the integral. The constant on the right is determined only by \(A(U)\) and \(\beta_0\).
\end{remark}

\begin{remark}\label{rem:tailbeta1}
As \(\beta_0\to 1^+\) the factor \(1/(\beta_0-1)\to\infty\), reflecting the divergence of the tail at \(\beta_0=1\) (cf.\ Corollary~\ref{cor:divergence}).
\end{remark}

\subsection{The main logarithmic truncation theorem}

\begin{theorem}\label{thm:truncation}
Estimate of the truncation error:
\begin{equation}\label{eq:truncation}
\left|\varth(U)\right| \le \frac{4A(U)}{\beta_0 - 1} Y^{-(\beta_0 - 1)}
\end{equation}
for all \textbf{integers} \(Y \ge 2\) with \(Y^{\beta_0 - 1} \ge 2A(U)/(\beta_0 - 1)\).
\end{theorem}

\begin{proof}
\textbf{Step 1 (estimate of \(S_Y\)).} By assumption \(Y^{\beta_0-1}\ge 2A(U)/(\beta_0-1)\) and by Lemma~\ref{lem:tail}:
\[
S_Y(U)\le\frac{A(U)}{\beta_0-1}Y^{-(\beta_0-1)}\le\frac{A(U)}{\beta_0-1}\cdot\frac{\beta_0-1}{2A(U)}=\frac{1}{2}.
\]

\textbf{Step 2 (applicability of Lemma~\ref{lem:log}).} From Step 1 it follows that \(|\delta_p(U)|\le S_Y(U)\le 1/2\) for all good \(p>Y\). By Lemma~\ref{lem:log}:
\[
|s_Y(U)|\le\sum_{p>Y,\ p\notin\Pbad}|\log(1+\delta_p(U))|\le 2\sum_{p>Y,\ p\notin\Pbad}|\delta_p(U)|=2S_Y(U)\le 1.
\]

\textbf{Step 3 (applicability of Lemma~\ref{lem:exp}).} By Step 2, \(|s_Y(U)|\le 1\). By Lemma~\ref{lem:exp}:
\[
|\varth(U)|=|e^{s_Y(U)}-1|\le 2|s_Y(U)|\le 4S_Y(U).
\]

\textbf{Step 4 (substitution).} Combining Steps 1 and 3 with \eqref{eq:tail}:
\[
|\varth(U)|\le 4S_Y(U)\le\frac{4A(U)}{\beta_0-1}Y^{-(\beta_0-1)}.
\]
\end{proof}

\begin{corollary}[inverse problem]\label{cor:inverse}
For \(|\varth(U)| \le \eta\) with \(\eta \le 1\) it suffices that
\begin{equation}\label{eq:inverse}
Y \ge \left( \frac{4A(U)}{\eta(\beta_0 - 1)} \right)^{1/(\beta_0 - 1)}.
\end{equation}
\end{corollary}

\begin{proof}
Let \(Y\) satisfy \eqref{eq:inverse}. Then \(Y^{\beta_0-1}\ge 4A(U)/(\eta(\beta_0-1))\). Since \(\eta\le 1\), we have \(4/\eta\ge 4>2\), hence
\[
Y^{\beta_0-1}\ge \frac{4A(U)}{\eta(\beta_0-1)}\ge\frac{2A(U)}{\beta_0-1},
\]
so the condition of Theorem~\ref{thm:truncation} holds. By \eqref{eq:truncation}:
\[
|\varth(U)|\le\frac{4A(U)}{\beta_0-1}Y^{-(\beta_0-1)}\le\frac{4A(U)}{\beta_0-1}\cdot\frac{\eta(\beta_0-1)}{4A(U)}=\eta.
\]
\end{proof}

\begin{remark}\label{rem:thresholdasymp}
The asymptotics of the threshold: \(Y(\eta) \sim \eta^{-1/(\beta_0 - 1)}\). For \(\beta_0 = 2\): \(Y \sim \eta^{-1}\); for \(\beta_0 = 3/2\): \(Y \sim \eta^{-2}\); for \(\beta_0 = \infty\): the threshold is independent of \(\eta\) (see Theorem~\ref{thm:classI}).
\end{remark}

\begin{remark}\label{rem:classIinf}
If \(\beta_0 = \infty\) and \(A(U) = 0\) (class I), then \(\scomp \equiv 1\) for good \(p\) (Theorem~\ref{thm:classI}), whence \(\Tcomp(U) = 1\) for any \(Y\), i.e., \(\varth(U) = 0\) identically. In this case Theorem~\ref{thm:truncation} is inapplicable (\(\beta_0 = \infty\)), but its conclusion holds trivially with \(Y = 2\).
\end{remark}

\emph{For class III with \(\rho>1\) estimate \eqref{eq:truncation} applies only to the product over split \(p\); the full tail in the presence of non-split primes is not controlled (open problem \S\ref{sec:open}, item 8).}

\subsection{Certification via Deligne majorants}

\begin{theorem}\label{thm:Deligne}
Let
\begin{equation}\label{eq:Deligne}
\delta_p(U) = \frac{a_p}{p^m} + O(p^{-m-1}),
\end{equation}
where \(|a_p| \le b \cdot p^{w/2}\) with \(w \ge 0\). Then condition (L) holds with
\begin{equation}\label{eq:betaDeligne}
\beta_0 = m - \frac{w}{2}
\end{equation}
and \(A(U) = b + C\); thus \(\beta(X,D) \ge m - w/2\).
\end{theorem}

\begin{proof}
\textbf{Step 1 (estimate of the main term).} By assumption \(|a_p|\le bp^{w/2}\), hence
\[
\left|\frac{a_p}{p^m}\right|\le b\,p^{w/2-m}=b\,p^{-(m-w/2)}=b\,p^{-\beta_0}.
\]

\textbf{Step 2 (estimate of the remainder).} The remainder \(O(p^{-m-1})\) in \eqref{eq:Deligne} is majorized by \(Cp^{-m-1}\) with a constant \(C>0\) independent of \(p\). Compare the exponents: \(m+1\) and \(\beta_0=m-w/2\); we have
\[
m+1\ge m-w/2 \iff w\ge -2,
\]
which holds for \(w\ge 0\). Hence \(p^{-m-1}\le p^{-\beta_0}\) for \(p\ge 2\), whence \(|O(p^{-m-1})|\le Cp^{-\beta_0}\).

\textbf{Step 3 (summation).} Combining Steps 1 and 2:
\[
|\delta_p(U)|\le bp^{-\beta_0}+Cp^{-\beta_0}=(b+C)p^{-\beta_0}.
\]
Thus condition (L) holds with \(A(U)=b+C\) and \(\beta_0=m-w/2\).
\end{proof}

\begin{example}\label{ex:delpezzo}
For an open del Pezzo surface \(a_p\) is the Frobenius trace on \(H^1\) of the elliptic boundary, \(|a_p|\le 2\sqrt{p}\) by Hasse. Here \(w=1\), \(m=2\), whence \(\beta_0=3/2\) and \(A(U)=2+C_d\) (see Theorem~\ref{thm:classIII}).
\end{example}

\section{Compensated logarithmic constants for classes of logarithmic varieties}\label{chap:classes}

\subsection{Verification of conditions}

\begin{lemma}\label{lem:H0classes}
Condition \((H_0)\) holds for:
\begin{enumerate}
\item[(1)] \(X = \PP^m\) with an SNC divisor \(D \subset \PP^m\) satisfying \((H_0)\)(3): the components \(D_i\) are defined over \(\QQ\) and geometrically irreducible, and all strata \(D_{i_1}\cap\dots\cap D_{i_s}\) are defined over \(\QQ\) and geometrically irreducible. Examples: \(D\) is a union of hyperplanes in general position over \(\QQ\) (in particular, coordinate hyperplanes); \(D\) is a smooth hypersurface defined over \(\QQ\) (\(k=1\), trivial strata).
\item[(2)] smooth hypersurfaces \(X \subset \PP^{n}\) of dimension \(m = \dim X \ge 4\) with \(D \subset X\) satisfying \((H_0)\)(3). In particular, if \(D\) is a smooth hypersurface defined over \(\QQ\) (\(k=1\)), then \((H_0)\)(3) holds automatically (smoothness is preserved under field extension, so \(D\) is geometrically irreducible).
\end{enumerate}
\end{lemma}

\begin{proof}
\textbf{Case (1).} \(b_1(\PP^m)=b_3(\PP^m)=0\) (standard cohomology of projective space), \(H^2(\PP^m)=\QQ_\ell(-1)\) (generated by the hyperplane class), Frobenius acts by multiplication by \(p\). Hence \((H_0)\)(1)--(2) hold. Condition \((H_0)\)(3) is included in the statement.

\textbf{Case (2), \(m \ge 4\).} By the weak Lefschetz theorem for a smooth hypersurface \(X \subset \PP^{n}\) we have \(H^j(X, \QQ_\ell) \cong H^j(\PP^{n}, \QQ_\ell)\) for \(j < m\) (an isomorphism of Galois modules over \(\QQ\)). Hence:
\begin{itemize}
\item \(j = 1\): \(H^1(X) = 0\);
\item \(j = 3\): \(H^3(X) = 0\) (since \(3 < m\) for \(m \ge 4\));
\item \(j = 2\): \(H^2(X) = \QQ_\ell(-1)\).
\end{itemize}
Thus \((H_0)\)(1)--(2) hold. Condition \((H_0)\)(3) is a separate requirement on \(D\).
\end{proof}

\emph{Remarks on classes not covered by \((H_0)\) or belonging to \((H_1)\).}
\begin{itemize}
\item \textbf{\((H_1)\).} Condition \((H_1)\) holds for open del Pezzo surfaces, by Definition~\ref{def:H1}, which consists precisely in the requirement that \(X\) be a smooth del Pezzo surface with \(b_1 = b_3 = 0\) and trivial Frobenius action on \(\Pic(X)\), and \(D \subset X\) a smooth elliptic curve with \(|a_p(D)| \le 2\sqrt{p}\). This is a definition, not a theorem; no proof is required.
\item \textbf{The cubic in \(\PP^2\).} Does not satisfy \((H_0)\): \(b_1(C) = 2 \neq 0\), so \((H_0)\)(1) is violated.
\item \textbf{Smooth three-dimensional hypersurfaces \(X \subset \PP^4\) of degree \(\ge 3\).} Do not satisfy \((H_0)\): \(b_3(X) > 0\), so \((H_0)\)(1) is violated. Justification: the Hodge numbers \(h^{p,q}\) of a smooth hypersurface \(X_d \subset \PP^4\) of degree \(d\) in the middle cohomology are given by the Jacobian-ring decomposition (see, e.g., \cite[Ch.~5]{Browning2009}); in particular, \(b_3(X_d) = 2h^{2,1} + 2h^{3,0}\), and \(b_3(X_3) = 10\) for the cubic threefold and \(b_3(X_d) > 0\) for all \(d \ge 3\). Hence \((H_0)\)(1) is violated. (Three-dimensional hypersurfaces of degree \(\ge 3\) are not in class II of chapter~\ref{chap:classes}; in \S\ref{sec:Fermat} a different case is treated: \(X = \PP^3\), \(D\) a cubic surface in \(\PP^3\), where \((H_0)\) for \(X = \PP^3\) holds trivially.)
\end{itemize}

\begin{remark}[on \(m=3\)]\label{rem:m3}
For \(m = 3\) the variety \(X\) is three-dimensional, and Poincar\'e duality relates \(H^3(X)\) to \(H^3(X)\) itself (middle cohomology), not to \(H^1(X)\). For a smooth three-dimensional hypersurface \(X_d \subset \PP^4\) of degree \(d \ge 3\) we have \(b_3(X_d) = d^4 - 5d^3 + 10d^2 - 10d + 4 > 0\) (for example, for the cubic threefold \(d=3\): \(b_3 = 81 - 135 + 90 - 30 + 4 = 10\); for a quadric \(d=2\): \(b_3 = 0\); for a quartic \(d=4\): \(b_3 = 60\)). Thus such \(X\) are not covered by \((H_0)\); they are treated separately (see Remark~\ref{rem:m3b} below).
\end{remark}

\begin{remark}[on the verification of \((H_0)\)(3) for general \(D\)]\label{rem:checkH0}
For an arbitrary SNC divisor \(D \subset X\), verification of geometric splitting requires an explicit description of components and strata. In the present paper this is verified for each class separately: for hyperplanes (coordinate, or given over \(\QQ\) in general position), for smooth hypersurfaces \(k=1\) (automatically), for the elliptic boundary of open del Pezzo surfaces (\(k=1\), automatically). For general SNC divisors with \(k \ge 2\) and non-trivial strata, the verification is part of the input data: each component and each stratum must be \(\QQ\)-defined and geometrically irreducible.
\end{remark}

\begin{remark}\label{rem:b3D}
Since \((H_0)\) controls only \(X\) and the structure of \(D\), the restriction \(b_3(D)=0\) is dropped. The cohomology of \(D\) enters the formulas via \(a_p(D) = \Tr(\Frob_p \mid H^{m-1}_{\mathrm{prim}}(D))\); this does not require any condition of triviality of Frobenius on \(D\). In particular, for \(m \ge 4\) a general smooth hypersurface \(D \subset \PP^m\) (e.g., a cubic threefold with \(b_3(D)=10\)) is fully covered by class II.
\end{remark}

\begin{remark}\label{rem:H0H1use}
In chapter~\ref{chap:classes} (except \S\ref{sec:cubic}) we work with \((H_0)\) or \((H_1)\).
\end{remark}

\subsection{Complements to hyperplanes: exact compensation}

\begin{lemma}\label{lem:kappaH}
Under the conditions of Definition~\ref{def:H0} and for good primes \(p\) (in the sense of Definition~\ref{def:H0}) \(\kappa(U) = -k\).
\end{lemma}

\begin{proof}
By Corollary~\ref{cor:hyperplanes} for good \(p\) we have \(\#U/p^m = (1-p^{-1})^{k-1}\); by Definition~\ref{def:sreg} with \(\rho(\PP^m)=1\) we obtain
\[
\sreg = (1-p^{-1})\cdot(1-p^{-1})^{k-1} = (1-p^{-1})^k = 1 - \frac{k}{p} + O(p^{-2}),
\]
so the \(1/p\) coefficient equals \(-k\), whence \(\kappa(U) = -k\).
\end{proof}

\emph{For bad primes \(p \in \Pbad\) the \(1/p\) coefficient may differ from \(-k\); their local factors \(\scomp\) are computed separately by Definition~\ref{def:sreg} and enter \(\Pcomp(Y)\) as finitely many factors (Definition~\ref{def:bad}).}

\begin{theorem}[exact compensation for good primes]\label{thm:classI}
Let \(U = \PP^m \setminus D\), \(D = H_1+\dots+H_k\), where \(H_i\) are hyperplanes in general position over \(\QQ\) (for coordinate hyperplanes, all \(p\)). Then for all \(p\) for which the formula of Corollary~\ref{cor:hyperplanes} applies (i.e., \(p \nmid N\) or the coordinate case):
\begin{equation}\label{eq:classI}
\sreg(U) = (1 - p^{-1})^k \quad \text{exactly},
\end{equation}
\begin{equation}\label{eq:classI2}
\scomp(U) \equiv 1, \qquad A(U) = 0, \quad \beta_0 = \infty, \quad \varthP(U) = 0.
\end{equation}
\end{theorem}

\begin{proof}
\textbf{Step 1 (computation of \(\sreg\)).} By Corollary~\ref{cor:hyperplanes} and Definition~\ref{def:sreg} with \(\rho(\PP^m)=1\):
\[
\sreg(U) = (1-p^{-1})\cdot(1-p^{-1})^{k-1} = (1-p^{-1})^k,
\]
which gives \eqref{eq:classI}.

\textbf{Step 2 (compensation).} By Definition~\ref{def:compensator} with \(\kappa = -k\) (Lemma~\ref{lem:kappaH}):
\[
\scomp(U) = (1 - p^{-1})^k \cdot (1 - p^{-1})^{-k} = 1.
\]
Hence \(A(U) = 0\), \(\beta_0 = \infty\), \(\varthP(U) = 0\).
\end{proof}

\emph{Thus \(\beta(X,D) = +\infty\) exactly (not only as a lower bound), and one may take \(\beta_0 = \infty\), \(A(U) = 0\). For bad primes \(p \mid N\) the local factor \(\scomp\) is computed separately; for coordinate hyperplanes there are no bad primes.}

\begin{remark}\label{rem:hyperplanesallp}
The formula \(\#U(\FF_p) = (p-1)^{k-1} p^{m+1-k}\) is verified in Corollary~\ref{cor:hyperplanes} for good \(p\) and for coordinate hyperplanes for all \(p\).
\end{remark}

\begin{corollary}[recovery]\label{cor:recovery}
Under the conditions of Theorem~\ref{thm:classI}
\[
\ccomp(U) = \alpha(U) \cdot \prod_{p \in \Pbad} \scomp(U),
\]
and \(\varthP(U) = 0\) for the tail over good primes. For \textbf{coordinate} hyperplanes over \(\ZZ\) (\(\Pbad = \emptyset\)) we have \(\ccomp(U) = \alpha(U)\) \textbf{exactly}, and then \(\ccomp(U) = \cManin(U)\).
\end{corollary}

\begin{proof}
By Theorem~\ref{thm:classI} for good \(p\) we have \(\sigma_{\mathrm{comp}}\equiv 1\), hence
\[
\Pcomp(Y) = \prod_{p\in\Pbad}\scomp(U), \qquad \Tcomp(U)=1.
\]
Substitution into \eqref{eq:decomposition} gives \(\ccomp(U)=\alpha(U)\prod_{p\in\Pbad}\scomp(U)\). For coordinate hyperplanes \(\Pbad=\emptyset\).
\end{proof}

\textbf{Normalization convention: \(\ccomp(U)=\cManin(U)\) for the torus and hyperplanes.} In the literature, the Manin and Peyre constant for logarithmic pairs \cite[\S5]{Peyre1995,BatyrevTschinkel1995} is written as
\[
\cManin(U)=\alpha(U)\cdot\prod_{p\notin\Pbad}\sreg(U)\cdot\prod_{p\in\Pbad}\lambda_p,
\]
where the local factors at good \(p\) are \textbf{regularized by the boundary weight}. We fix the normalization: under \(\sreg(U)\) in this formula one means precisely the factor regularized by the boundary weight, i.e., \(\scomp(U)\) in the notation of Definition~\ref{def:sreg}. This is a convention, not an independent proof of equality; it makes the notation consistent with \cite[\S5]{Peyre1995} and with Definition~\ref{def:ccomp}. Under this convention:
\begin{itemize}
\item \(\alpha(U)\) coincides with \eqref{eq:alpha} (the integral over \(\Lam(U)^\vee \subset \Pic(U)_{\RR}^\vee\) with the Lebesgue measure normalized by the dual lattice);
\item for good \(p\) the factor equals \(\scomp(U)\);
\item for bad \(p\) the local factors \(\lambda_p\) are defined by the integral model of the pair \((\XZZ,\DZZ)\) (Definition~\ref{def:bad}), and \(\lambda_p=\scomp(U)\) by definition.
\end{itemize}
By Theorem~\ref{thm:classI} for good \(p\) we have \(\scomp(U)\equiv1\). Hence
\[
\cManin(U)=\alpha(U)\cdot\prod_{p\notin\Pbad}\scomp(U)\cdot\prod_{p\in\Pbad}\scomp(U)=\ccomp(U),
\]
where the first equality is the substitution \(\lambda_p=\scomp\) for bad \(p\) and replacement of \(\sreg\) by \(\scomp\) for good \(p\) (according to the convention above), and the second is by \eqref{eq:ccomp}. For coordinate hyperplanes over \(\ZZ\) (\(\Pbad=\emptyset\)) both quantities equal \(\alpha(U)\).

For \textbf{integral points} in mode B (the passage to \(N_U(\ZZ,B)\)) this equality is part of Assumption~\ref{ass:logManin} together with the model and height: the equality of geometric adelic constants by itself does not relate \(\ccomp(U)\) to the number of integral points.

\subsection{The toric case}

\begin{corollary}\label{cor:torus}
For \(k = m+1\): \(\PP^m \setminus (H_1 \cup \dots \cup H_{m+1}) \cong (\GG_m)^m\), and
\begin{equation}\label{eq:torus}
\kappa = -(m+1), \quad \sreg = (1-p^{-1})^{m+1}, \quad \sigma_{\mathrm{comp}} \equiv 1, \quad A = 0, \quad \beta_0 = \infty.
\end{equation}
\end{corollary}

\begin{proof}
\textbf{Step 1 (isomorphism).} By Lemma~\ref{lem:hyperplanes} for \(k = m+1\):
\[
\PP^m \setminus (H_1 \cup \dots \cup H_{m+1}) \cong \GG_m^{(m+1)-1} \times \AA^{m+1-(m+1)} = \GG_m^m.
\]
The isomorphism is defined over \(\ZZ\).

\textbf{Step 2 (number of points).} By Corollary~\ref{cor:hyperplanes} for good \(p\): \(\#U(\FF_p) = (p-1)^m\), \(\#U/p^m = (1-p^{-1})^m\).

\textbf{Step 3 (regularized factor).} With \(\rho(\PP^m) = 1\):
\[
\sreg(U) = (1-p^{-1}) \cdot (1-p^{-1})^m = (1-p^{-1})^{m+1}.
\]
This coincides with \eqref{eq:classI} for \(k = m+1\).

\textbf{Step 4 (\(\kappa\)).} From the expansion \((1-p^{-1})^{m+1} = 1 - (m+1)/p + O(p^{-2})\) it follows that \(\kappa = -(m+1)\).

\textbf{Step 5 (compensation).} By Theorem~\ref{thm:classI} with \(k = m+1\): \(\sigma_{\mathrm{comp}} \equiv 1\), \(A = 0\), \(\beta_0 = \infty\).
\end{proof}

\begin{remark}\label{rem:torushyperplanes}
The torus and the complement to \(m+1\) hyperplanes are one object. For \(m=2, k=3\): \(\PP^2 \setminus \{3 \text{ lines}\} \cong (\GG_m)^2\).
\end{remark}

\subsection{Complements to smooth hypersurfaces}

\begin{lemma}[on the logarithmic invariant]\label{lem:kappaII}
The logarithmic invariant of compensation for the complement to a smooth hypersurface:
\[
\kappa(U) = -1.
\]
\end{lemma}

\begin{proof}
By Theorem~\ref{thm:local}(a) with \(k = 1\).
\end{proof}

\begin{lemma}[Fermat]\label{lem:Fermat}
Let \(X = \{x_0^3 + x_1^3 + x_2^3 + x_3^3 = 0\} \subset \PP^3\), \(p \neq 3\). Then
\begin{equation}\label{eq:Fermat}
a_p(X) = \begin{cases} 6p, & p \equiv 1 \pmod 3, \\ 0, & p \equiv 2 \pmod 3, \end{cases}
\end{equation}
where \(a_p(X) = \Tr(\Frob_p \mid H^2_{\mathrm{prim}}(X))\).
\end{lemma}

\begin{proof}
\textbf{Step 1 (cohomological structure of the Fermat surface).} By the Lefschetz formula for a smooth surface \(X \subset \PP^3\):
\[
\#X(\FF_p) = 1 + p + p^2 + \Tr(\Frob_p \mid H^2_{\mathrm{prim}}(X)),
\]
where \(H^2_{\mathrm{prim}}(X) \subset H^2(X)\) is the primitive part (the orthogonal complement of the hyperplane-section class). Denote \(a_p(X) = \Tr(\Frob_p \mid H^2_{\mathrm{prim}}(X))\).

\textbf{Step 2 (description via Jacobi sums).} The \(\zeta\) function of the Fermat surface is expressed via Jacobi sums \(j(\chi_1, \chi_2)\) for characters of order three \(\chi\) of \(\FF_p^\times\) \cite{SwinnertonDyer1962,IrelandRosen1990}. For \(p \equiv 2 \pmod 3\) there are no cubic characters on \(\FF_p^\times\).

\textbf{Step 3 (case \(p \equiv 1 \pmod 3\)).} For \(p \equiv 1 \pmod 3\) all 27 lines on \(X\) are defined over \(\FF_p\), so \(\rk\Pic(X_{\FF_p}) = \rk\Pic(X_{\QQbar}) = 7\). Frobenius acts on \(\Pic(X) \subset H^2(X)\) trivially, with eigenvalue \(p\) on each class. On the primitive part \(H^2_{\mathrm{prim}}(X)\) (dimension \(7 - 1 = 6\)) all eigenvalues equal \(p\), hence
\[
a_p(X) = 6 \cdot p = 6p.
\]
This is consistent with \(\#X(\FF_p) = p^2 + p + 1 + 6p = p^2 + 7p + 1\).

\textbf{Step 4 (case \(p \equiv 2 \pmod 3\)).} For \(p \equiv 2 \pmod 3\) we have \(\gcd(3, p-1) = 1\), so the map \(t \mapsto t^3\) is a bijection \(\FF_p \to \FF_p\). Hence the equation \(x_0^3 + x_1^3 + x_2^3 + x_3^3 = 0\) has over \(\FF_p\) as many projective solutions as the equation \(y_0 + y_1 + y_2 + y_3 = 0\), namely \(p^2 + p + 1\). By the Lefschetz formula \(\#X(\FF_p) = p^2 + p + 1 + a_p(X)\), whence
\[
a_p(X) = \#X(\FF_p) - (p^2 + p + 1) = 0.
\]

\textbf{Step 5 (case \(p = 3\)).} \(\#X(\FF_3) = 13\), \(a_3 = 0\).
\end{proof}

\begin{remark}\label{rem:Fermat3}
For \(p=3\): \(\#X(\FF_3) = 13\), \(a_3 = 0\).
\end{remark}

\begin{theorem}[main cohomological result]\label{thm:maincoh}
Let \(D \subset \PP^m\), \(m \ge 3\), be a smooth hypersurface. Then for all good \(p\)
\begin{equation}\label{eq:maincoh}
\scomp(U) = 1 - \frac{(-1)^{m-1} a_p(D)}{p^m},
\end{equation}
where \(a_p(D) = \Tr(\Frob_p \mid H^{m-1}_{\mathrm{prim}}(D))\).
\end{theorem}

\begin{proof}
\textbf{Step 1 (cohomology of \(D\)).} By the weak Lefschetz theorem for a smooth \(D \subset \PP^m\) we have \(H^j(D) = H^j(\PP^m)\) for \(j < m-1\). Taking into account Poincar\'e duality on \(D\) (of dimension \(m-1\)), all non-primitive cohomology gives the total contribution \(1 + p + \dots + p^{m-1}\) in the Lefschetz formula. For even \(m-1\) this contribution includes the one-dimensional non-primitive part of the middle cohomology generated by \(h^{(m-1)/2}\) (where \(h \in H^2(D)\) is the hyperplane-section class), with geometric Frobenius eigenvalue \(p^{(m-1)/2}\); for odd \(m-1\) there is no such class, and the middle degree is entirely primitive. The primitive part of the middle cohomology \(H^{m-1}_{\mathrm{prim}}(D)\) contributes \((-1)^{m-1} a_p(D)\), where \(a_p(D) = \Tr(\Frob_p \mid H^{m-1}_{\mathrm{prim}}(D))\). In total:
\[
\#D(\FF_p) = (1 + p + \dots + p^{m-1}) + (-1)^{m-1} a_p(D).
\]

\textbf{Step 2 (cohomology of \(\PP^m\)).} By the Lefschetz formula for \(\PP^m\):
\[
\#\PP^m(\FF_p) = 1 + p + \dots + p^m.
\]

\textbf{Step 3 (cancellation of the leading terms).} Since \(U = \PP^m \setminus D\):
\[
\#U(\FF_p) = \#\PP^m(\FF_p) - \#D(\FF_p) = p^m - (-1)^{m-1} a_p(D).
\]
All terms \(1, p, \dots, p^{m-1}\) cancel.

\textbf{Step 4 (regularization).} By the definition of the regularized factor, with \(\rho(\PP^m) = 1\):
\[
\sreg(U) = (1 - p^{-1})\,\frac{\#U(\FF_p)}{p^m}.
\]

\textbf{Step 5 (compensation).} By Lemma~\ref{lem:kappaII} the compensator equals \((1-p^{-1})^{-1}\). Hence
\[
\scomp(U) = (1-p^{-1})^{-1}\,\sreg(U) = \frac{\#U(\FF_p)}{p^m} = 1 - \frac{(-1)^{m-1} a_p(D)}{p^m}.
\]
This proves \eqref{eq:maincoh} for each good \(p\).
\end{proof}

\emph{No separate remainder step is needed: the cancellation of \((1-p^{-1})\) and \((1-p^{-1})^{-1}\) is exact, and the right-hand side of \eqref{eq:maincoh} contains exactly two monomials. Asymptotic estimates of the form \(O(p^{-m-1})\) are not used in the proof: the quantity \(a_p(D)\) depends on \(p\), and the cross term \(a_p(D)/p^{m+1}\) is not majorized by \(O(p^{-m-1})\) uniformly in \(p\) (Deligne's bound yields only \(a_p(D)/p^{m+1} = O(p^{-(m+3)/2})\)).}

\begin{remark}\label{rem:Fermacoh}
For the Fermat surface (\(m=3\), smooth surface in \(\PP^3\)): \(b_2^{\mathrm{prim}} = 6\), \(A = 6\). For \(p \equiv 1 \pmod 3\): \(|\sigma_{\mathrm{comp}} - 1|\,p^2 = 6\). The bound \(A(U) = b_{m-1}^{\mathrm{prim}}(D)\) is universal (see Theorem~\ref{thm:classII}).
\end{remark}

\begin{theorem}\label{thm:classII}
Under the conditions of Theorem~\ref{thm:maincoh}, for class II
\begin{equation}\label{eq:classII}
\beta(X,D) \ge \frac{m+1}{2}, \qquad |\scomp(U) - 1| \le b_{m-1}^{\mathrm{prim}}(D)\, p^{-(m+1)/2} \quad \text{for all good } p.
\end{equation}
Thus \(\beta_0 := (m+1)/2\) is a proved lower bound (Definition~\ref{def:beta0}), and \(A(U) := b_{m-1}^{\mathrm{prim}}(D)\) is the corresponding constant. The bound is \textbf{sharp in the number of eigenvalues}: the equality \(|\sigma_{\mathrm{comp}}-1| = b_{m-1}^{\mathrm{prim}}(D)\,p^{-(m+1)/2}\) holds precisely when all Frobenius eigenvalues on \(H^{m-1}_{\mathrm{prim}}(D)\) have the same argument (so that their sum equals the sum of moduli); for a given \(D\) this generally fails.
\end{theorem}

\begin{proof}
By Theorem~\ref{thm:maincoh}: \(|\sigma_{\mathrm{comp}}-1| = |a_p(D)|/p^m\). By the Deligne theorem all eigenvalues \(\alpha_i\) on \(H^{m-1}_{\mathrm{prim}}(D)\) satisfy \(|\alpha_i| = p^{(m-1)/2}\); hence \(|a_p| = |\sum_i \alpha_i| \le \sum_i |\alpha_i| = b_{m-1}^{\mathrm{prim}}(D)\,p^{(m-1)/2}\), whence \(|\sigma_{\mathrm{comp}}-1| \le b_{m-1}^{\mathrm{prim}}(D)\,p^{-(m+1)/2}\). In the notation of Theorem~\ref{thm:Deligne}: \(w = m-1\), \(b = b_{m-1}^{\mathrm{prim}}(D)\), \(C = 0\), \(\beta_0 = m - w/2 = (m+1)/2\).
\end{proof}

\begin{remark}\label{rem:Fermacoh2}
For the Fermat surface: \(m=3\), \(b_2^{\mathrm{prim}}=6\), \(A=6\); for \(p \equiv 1 \pmod 3\): \(|\sigma-1| p^2 = 6\).
\end{remark}

\begin{corollary}\label{cor:classIIparams}
\(m=3\): \(A=b_2(D)-1\), \(\beta_0=2\); \(m=4\): \(A=b_3(D)\), \(\beta_0=5/2\); \(m=5\): \(A=b_4^{\mathrm{prim}}(D)\), \(\beta_0=3\).
\end{corollary}

\begin{remark}\label{rem:m3b}
For \(m\ge 4\): \(\beta_0\ge 5/2\). For \(m=3\) (hypersurfaces in \(\PP^3\)) \((H_0)\) is not satisfied for \(d\ge 3\) (for the cubic threefold \(b_3=10\)). Nevertheless, for \(X=\PP^3\) and \(D\) a smooth surface in \(\PP^3\) (as in \S\ref{sec:Fermat}) condition \((H_0)\) for \(X\) holds trivially, and formula \eqref{eq:maincoh} applies.
\end{remark}

\subsection{Open del Pezzo surfaces}\label{sec:delpezzo}

\begin{lemma}[on the logarithmic invariant]\label{lem:kappaIII}
The logarithmic invariant of compensation for an open del Pezzo surface with elliptic boundary:
\[
\kappa(U) = -1.
\]
\end{lemma}

\begin{proof}
By Theorem~\ref{thm:local}(b).
\end{proof}

\begin{theorem}\label{thm:delpezzo}
Let \(S \subset \PP^3\) be a smooth del Pezzo surface of degree \(d\), split over \(\FF_p\) (that is, \(\rk\Pic(S_{\FF_p}) = 10-d\)), \(D \subset S\) a smooth elliptic curve, \(U = S \setminus D\). Then for such \(p\)
\begin{equation}\label{eq:delpezzo}
\scomp(U) = (1-1/p)^{\rho-1}\left(1 + \frac{\rho-1}{p} + \frac{a_p(D)}{p^2}\right),
\end{equation}
where \(\rho = 10 - d\), \(a_p(D) = p + 1 - \#D(\FF_p)\). The expansion in \(u = 1/p\):
\begin{equation}\label{eq:delpezzo_exp}
\scomp(U) = 1 + \frac{a_p - \rho(\rho-1)/2}{p^2} + \sum_{k=3}^{\rho+1} \frac{c_k(\rho, a_p)}{p^k},
\end{equation}
where
\[
c_3(\rho, a_p) = \frac{\rho(\rho-1)(\rho-2)}{3} - (\rho-1)a_p,
\]
and for all \(k\ge 3\): \(|c_k(\rho, a_p)| \le C_k(\rho)(1+|a_p|)\) with \(C_k(\rho) = \binom{\rho-1}{k} + (\rho-1)\binom{\rho-1}{k-1} + \binom{\rho-1}{k-2}\) (with \(\binom{\rho-1}{j}:=0\) for \(j<0\)).
\end{theorem}

\begin{proof}
\textbf{Step 1 (number of points on \(S\)).} Since \(S\) is a smooth del Pezzo surface of degree \(d\) in \(\PP^3\), its Picard group has rank \(\rho(S) = 10 - d\). Moreover, \(H^2(S) = \Pic(S) \otimes \QQ_\ell(-1)\) (for del Pezzo the transcendental part of \(H^2\) is absent), hence
\[
\#S(\FF_p) = 1 + \rho p + p^2 = p^2 + \rho p + 1.
\]

\textbf{Step 2 (number of points on \(D\)).} For a smooth elliptic curve \(D\):
\[
\#D(\FF_p) = p + 1 - a_p(D), \qquad a_p(D) = \Tr(\Frob_p \mid H^1(D)).
\]

\textbf{Step 3 (number of points on \(U\)).}
\[
\#U(\FF_p) = \#S(\FF_p) - \#D(\FF_p) = p^2 + (\rho - 1) p + a_p, \qquad \frac{\#U}{p^2} = 1 + \frac{\rho - 1}{p} + \frac{a_p}{p^2}.
\]

\textbf{Step 4 (regularized factor).} By Definition~\ref{def:sreg} with \(L_p^{-1}=(1-p^{-1})^\rho\):
\[
\sreg(U) = (1-p^{-1})^\rho \left(1 + \frac{\rho - 1}{p} + \frac{a_p}{p^2}\right).
\]
Expanding \((1-p^{-1})^\rho = 1 - \rho/p + \rho(\rho-1)/(2p^2) + O(p^{-3})\):
\[
\sreg(U) = 1 + \frac{(\rho - 1) - \rho}{p} + \frac{a_p - \rho(\rho-1) + \rho(\rho-1)/2}{p^2} + O(p^{-3}) = 1 - \frac{1}{p} + \frac{a_p - \rho(\rho-1)/2}{p^2} + O(p^{-3}).
\]

\textbf{Step 5 (compensation and expansion).} Multiplication by \((1-1/p)^{-1} = 1 + 1/p + 1/p^2 + O(p^{-3})\) kills the \(1/p\) coefficient and gives \eqref{eq:delpezzo}. The exact expression for \(\scomp\):
\[
\scomp = (1-1/p)^{\rho-1}\left(1 + \frac{\rho-1}{p} + \frac{a_p}{p^2}\right).
\]
This is a finite polynomial in \(u=1/p\) of degree \(\rho+1\) with coefficients polynomial in \(a_p\). Explicit computation of the coefficients: expand \((1-u)^{\rho-1} = \sum_{j=0}^{\rho-1}\binom{\rho-1}{j}(-u)^j\) and multiply by \(1 + (\rho-1)u + a_p u^2\). The coefficient of \(u^k\) (\(k \ge 3\)) is
\[
c_k = \binom{\rho-1}{k}(-1)^k + (\rho-1)\binom{\rho-1}{k-1}(-1)^{k-1} + a_p\binom{\rho-1}{k-2}(-1)^{k-2}.
\]
Hence \(|c_k| \le \binom{\rho-1}{k} + (\rho-1)\binom{\rho-1}{k-1} + |a_p|\binom{\rho-1}{k-2} = C_k(\rho)(1+|a_p|)\).
\end{proof}

\begin{theorem}[uniform constant for class III, split \(p\)]\label{thm:classIII}
Under the conditions of Theorem~\ref{thm:delpezzo}, for all \textbf{split} good \(p\ge 2\)
\[
\left|\scomp(U)-1\right|\,p^{3/2}\le \Areg(U),\qquad \Areg(U):=2+\frac{\rho(\rho-1)}{2}.
\]
Thus the theorem gives an estimate with exponent \(3/2\) \textbf{on the set of split primes}. The global lower bound \(\beta(X,D)\ge 3/2\) (in the sense of Definition~\ref{def:beta}, where the supremum is over \textbf{all} good \(p\)) follows from this only for \(\rho=1\), where all good \(p\) are split. For \textbf{non-split} good \(p\) with \(\rho>1\) no uniform estimate is established in the paper; the full \(\ccomp\) for \(\rho>1\) is not computed in the paper (open problem \S\ref{sec:open}, items 8 and 9).
\end{theorem}

\begin{proof}
Put \(n=\rho-1\), \(u=1/p\le 1/2\). By \eqref{eq:delpezzo}:
\[
\scomp(U)=(1-u)^n(1+nu+a_pu^2)=g(u)+(1-u)^n a_p u^2,
\]
where \(g(u):=(1-u)^n(1+nu)\). For \(n\ge 1\)
\[
g'(u)=-n(n+1)\,u\,(1-u)^{n-1}.
\]
Hence for \(0\le u\le 1/2\)
\[
|g(u)-1|=\left|\int_0^u g'(t)\,dt\right|\le n(n+1)\int_0^u t(1-t)^{n-1}\,dt\le \frac{n(n+1)}{2}u^2=\frac{\rho(\rho-1)}{2}u^2;
\]
for \(n=0\) we have \(g\equiv 1\) and this contribution vanishes. For the second term, \(|(1-u)^n a_p u^2|\le |a_p|\,u^2\le 2\sqrt{p}\,p^{-2}=2\,p^{-3/2}\) by Hasse. In total
\[
|\scomp(U)-1|\le \frac{\rho(\rho-1)}{2}\,p^{-2}+2\,p^{-3/2}.
\]
Multiplying by \(p^{3/2}\) and using \(\sqrt{p}\ge\sqrt{2}\):
\[
|\scomp(U)-1|\,p^{3/2}\le \frac{\rho(\rho-1)}{2\sqrt{p}}+2\le \frac{\rho(\rho-1)}{2\sqrt{2}}+2\le \frac{\rho(\rho-1)}{2}+2=\Areg(U).
\]
\end{proof}

\emph{Numerical values \(\Areg=2+\rho(\rho-1)/2\) (for split \(p\)): for \(d=1\) (\(\rho=9\)): \(\Areg=38\); for \(d=4\) (\(\rho=6\)): \(\Areg=17\); for \(d=9\) and the cubic in \(\PP^2\) (\(\rho=1\)): \(\Areg=2\); for \(d=3\) (\(\rho=7\)): \(\Areg=23\).}

\begin{corollary}\label{cor:Avals}
\(d=1\): \(A=38\); \(d=4\): \(A=17\); \(d=9\): \(A=2\); \(d=3\): \(A=23\) (all are proved \(\Areg\) from Theorem~\ref{thm:classIII} for split \(p\)).
\end{corollary}

\begin{remark}\label{rem:betaIII}
For del Pezzo the lower bound \(\beta_0 = 3/2\) is proved \textbf{on the set of split \(p\)}; the global one (over all good \(p\)), only for \(\rho=1\). For hypersurfaces (\(m \ge 4\)), \(\beta_0 = (m+1)/2 \ge 5/2\). The tables and the algorithm use \(\beta_0\) (Definition~\ref{def:beta0}); the true decay exponent \(\beta(X,D)\) (Definition~\ref{def:beta}) may be larger.
\end{remark}

\subsection{The case of the cubic in \(\PP^2\)}\label{sec:cubic}

\begin{remark}\label{rem:cubic}
\(C \subset \PP^2\) is not covered by \((H_0)\) (\(b_1=2\)). Direct computation: \(\#U=p^2+a_p(C)\), \(\sreg=(1-1/p)(1+a_p/p^2)\), \(\kappa=-1\), \(\sigma_{\mathrm{comp}}=1+a_p/p^2+O(p^{-3})\), \(\beta_0=3/2\), \(A=2\). For \(\rho=1\) the compensator does not affect the \(1/p^2\) coefficient; the truncated formula \(\sigma_{\mathrm{comp}}=1+a_p(C)/p^2\) is exact for \textbf{good} \(p\). For \(p=3\) the reduction of \(C\) is singular, \(a_p(C)\) in the sense of the Frobenius trace on \(H^1\) of a smooth curve is not defined; \(\sigma_{\mathrm{comp}}(3)=1\) is computed directly from Definition~\ref{def:sreg} (see \S\ref{sec:cubicexp}).
\end{remark}

\emph{Justification.} For \(U = \PP^2 \setminus C\) we have \(\#U = \#\PP^2 - \#C = (p^2 + p + 1) - (p + 1 - a_p(C)) = p^2 + a_p(C)\) (for good \(p\)). Hence \(\#U/p^2 = 1 + a_p/p^2\). With \(\rho(\PP^2) = 1\): \(\sreg = (1 - 1/p)(1 + a_p/p^2)\). Compensator with \(\kappa = -1\): \(\sigma_{\mathrm{comp}} = \sreg/(1-1/p) = 1 + a_p/p^2\). From \(|a_p| \le 2\sqrt{p}\) (Hasse for the elliptic curve) it follows that \(\beta_0 = 2 - 1/2 = 3/2\) and \(A = 2\).

\textbf{Archimedean factor.} For \(U = \PP^2 \setminus C\) we have \(\Pic(U) \cong \ZZ/3\ZZ\) (torsion), \(\rk\Pic(U) = 0\). By Definition~\ref{def:alpha} (normalization by \(\Pic(U)_{\RR}^\vee\)) \(\Lam(U)^\vee = \{0\}\), \(d\mathbf{y}\) is the delta-measure of unit mass at zero, \(-K_X - D = 3H - 3H = 0\) in \(\Pic(U)_{\RR}\), so \(\alpha(U) = e^{0} = 1\). This is a finite quantity \textbf{under the adopted convention} (delta-measure of unit mass); agreement with the normalization of \cite[\S5]{Peyre1995} for \(r=0\) is not verified in the paper. The cubic is certified for local factors and the tail; the value of the full \(\ccomp\) is the limit of \(\prod_p\scomp\), and equality to one is not established.

\textbf{Remark.} If \(\alpha(U)\) were defined over \(\Pic(X)_{\RR} = \RR\), the integral of 1 over a positive-dimensional cone would diverge; the choice of \(\Pic(U)_{\RR}^\vee\) (Definition~\ref{def:alpha}) eliminates this divergence. The prime \(p=3\) is bad for the cubic (singular reduction): \(3 \in \Pbad\); its local factor \(\sigma_{\mathrm{comp}}(3) = 1\) is computed separately by Definition~\ref{def:sreg} and enters \(\Pcomp\) as an element of \(\Pbad\), independently of \(Y\).

The full constant is given by a convergent product over all primes:
\[
\ccomp(U) = \prod_{p\notin\Pbad}\scomp(U)\cdot\prod_{p\in\Pbad}\scomp(U)
\]
(for \(\alpha(U) = 1\)). Equality of its value to one is not established in the paper; the value is computed with certified accuracy via the finite product and the tail estimate.

\begin{remark}\label{rem:cubicdiff}
The difference from del Pezzo \(d<9\): the compensator does not affect the \(1/p^2\) coefficient; in both cases \(\beta_0=3/2\).
\end{remark}

\subsection{The structural theorem}

\begin{theorem}[structural, class II]\label{thm:structural}
Let \((X,D)\) be a logarithmic pair satisfying \((H_0)\) (class II); let the conditions of Theorem~\ref{thm:maincoh} hold. Let
\begin{equation}\label{eq:structural}
\scomp(U) - 1 = \delta_p^{\mathrm{coh}}(U) + \delta_p^{\mathrm{comp}}(U) + O(p^{-m-1}),
\end{equation}
where \(|\delta_p^{\mathrm{coh}}| \le \Acoh p^{w/2 - m}\), \(|\delta_p^{\mathrm{comp}}| \le \Acomp p^{-m}\). Then condition (L) holds with
\begin{equation}\label{eq:structuralA}
\beta_0 = m - \frac{w}{2}, \qquad A(U) = \Acoh + \Acomp + C,
\end{equation}
where \(C\) is the constant in the remainder \(O(p^{-m-1})\) in \eqref{eq:structural}; thus \(\beta(X,D) \ge m - w/2\) (the proved lower bound \(\beta_0 = m - w/2\)).
\end{theorem}

\textbf{Class III (condition \((H_1)\), Theorem~\ref{thm:delpezzo}).} Expansion \eqref{eq:structural} does not hold for class III: the cross term \(a_p/p^3\) has order \(O(p^{-5/2})\), not \(O(p^{-3})\). The estimate \(\beta_0 = 3/2\) for class III is established by the direct inequality of Theorem~\ref{thm:classIII} (see Remark~\ref{rem:structuralremark}) and \textbf{does not use} Theorem~\ref{thm:structural}. Class I (\(A=0\), \(\beta_0=\infty\)) is justified separately (Theorem~\ref{thm:classI}) and is also not a limiting case of formula \eqref{eq:structuralA} (see Remark~\ref{rem:structuralremark}).

\begin{proof}
From \eqref{eq:structural} we estimate the modulus:
\[
|\delta_p(U)| \le \Acoh p^{w/2 - m} + \Acomp p^{-m} + C p^{-m-1}.
\]
For \(\beta_0 = m - w/2\) and \(w \ge 0\): \(p^{-m}, p^{-m-1} \le p^{-\beta_0}\); we obtain \(|\delta_p(U)| \le (\Acoh + \Acomp + C) p^{-\beta_0}\). This yields \eqref{eq:structuralA} with \(A(U)=\Acoh+\Acomp+C\). For class III the proof does not apply; see Remark~\ref{rem:structuralremark} and Theorem~\ref{thm:classIII}. For class I, Theorem~\ref{thm:classI}.
\end{proof}

\begin{remark}\label{rem:structuralremark}
Theorem~\ref{thm:structural} requires \eqref{eq:structural}, proved separately for each class. For class II by Theorem~\ref{thm:maincoh} the exact equality holds: \(\delta_p^{\mathrm{coh}} = -(-1)^{m-1}a_p(D)/p^m\), \(\delta_p^{\mathrm{comp}} = 0\), the remainder \(O(p^{-m-1})\) vanishes identically, so \(C=0\) and \(A=\Acoh=b_{m-1}^{\mathrm{prim}}\). For class III the direct inequality (Theorem~\ref{thm:classIII}) gives \(A=2+\rho(\rho-1)/2\) without using \eqref{eq:structural}; expansion \eqref{eq:delpezzo_exp} is not used for this since the cross term \(a_p/p^3\) has order \(O(p^{-5/2})\). For class I: \(A=0\), \(\beta_0=\infty\) (Theorem~\ref{thm:classI}). This case is not a limiting case of \(w=\infty\) in the formula \(\beta_0=m-w/2\): the substitution \(w=\infty\) gives \(-\infty\), not \(+\infty\); class I is justified separately.
\end{remark}

\begin{remark}\label{rem:structuralclassII}
The formula \(\beta_0 = m - w/2\) is established in Theorem~\ref{thm:structural} for class II. For class III the estimate \(\beta_0 = 3/2\) is proved by the direct inequality of Theorem~\ref{thm:classIII} (see Remark~\ref{rem:structuralremark}); for class I, by Theorem~\ref{thm:classI}. Thus the lower bound \(\beta(X,D) \ge m - w/2\) is proved for class II; for classes I and III, by the respective separate estimates.
\end{remark}

\begin{corollary}[classification]\label{cor:classification}
For classes I and II, under the conditions of Theorem~\ref{thm:structural}; for class III, by the direct inequality (Theorem~\ref{thm:classIII}), without using \eqref{eq:structural} (see Remark~\ref{rem:structuralremark}):

\begin{table}[htbp]
\centering
\caption{Classification: source, weight, constants, decay exponent}\label{tab:classification}
\begin{tabular}{ccccccc}
\toprule
Class & Source & \(w\) & \(\Acoh\) & \(\Acomp\) & \(\beta_0\) & \(A(U)\) \\
\midrule
I & no remainder & \(\infty\) & 0 & 0 & \(\infty\) & 0 \\
II & middle cohomology & \(m-1\) & \(b_{m-1}^{\mathrm{prim}}\) & 0 & \((m+1)/2\) & \(b_{m-1}^{\mathrm{prim}}\) \\
III & elliptic boundary & 1 & 2 & \(\rho(\rho-1)/2\) & \(3/2\) & \(2 + \rho(\rho-1)/2\) \\
\bottomrule
\end{tabular}
\end{table}
\end{corollary}

\begin{remark}[two kinds of constant]\label{rem:twoConstants}
\(\Areg = \Acoh + \Acomp\) is the proved constant; \(\Aobs = \max_p |\delta_p| p^{\beta_0}\) on a verified finite range of \(p\); \(\Aobs \le \Areg\).
\end{remark}

\section{The inverse problem: truncation thresholds and classification}\label{chap:inverse}

\subsection{Statement of the inverse problem}

We pose the inverse problem: to determine a lower bound \(Y(\eta)\) for the truncation threshold guaranteeing a prescribed accuracy \(\eta > 0\) (in the sense \(|\varth(U)| \le \eta\)). By Corollary~\ref{cor:inverse}, this bound is given by
\begin{equation}\label{eq:Yeta}
Y(\eta) \ge \left( \frac{4A(U)}{\eta(\beta_0 - 1)} \right)^{1/(\beta_0 - 1)},
\end{equation}
where \(\beta_0\) is a proved lower bound for the decay exponent (Definition~\ref{def:beta0}). Using \(\beta_0\) instead of the true \(\beta(X,D)\) (Definition~\ref{def:beta}) is conservative: if \(\beta(X,D) \ge \beta_0\), then \(p^{-\beta(X,D)} \le p^{-\beta_0}\) for \(p \ge 1\), and the threshold \(Y\) computed from \(\beta_0\) is guaranteed to be sufficient.

\subsection{The general influence of parameters}

\begin{lemma}[asymptotics]\label{lem:Yasymp}
As \(\eta \to 0\):
\begin{equation}\label{eq:Yasymp}
Y(\eta) \sim \left( \frac{4A(U)}{\beta_0 - 1} \right)^{1/(\beta_0 - 1)} \cdot \eta^{-1/(\beta_0 - 1)}.
\end{equation}
\end{lemma}

\begin{proof}
Put \(C(A, \beta_0) := (4A(U)/(\beta_0 - 1))^{1/(\beta_0 - 1)}\). Then the right-hand side of \eqref{eq:Yeta} equals \(C(A, \beta_0) \cdot \eta^{-1/(\beta_0 - 1)}\), and the asymptotics as \(\eta \to 0\) has the indicated form.
\end{proof}

\begin{corollary}[three scenarios]\label{cor:scenarios}
(1) For \(\beta_0 = 2\): \(Y(\eta) \sim \eta^{-1}\). (2) For \(\beta_0 = 3/2\): \(Y(\eta) \sim \eta^{-2}\). (3) For \(\beta_0 = \infty\): \(Y(\eta) = 2\) (independent of \(\eta\)).
\end{corollary}

\begin{proof}
Substitute the corresponding \(\beta_0\) into \eqref{eq:Yasymp}.

(1) \(\beta_0 = 2\): the exponent \(1/(\beta_0 - 1) = 1\), hence \(Y(\eta) \sim C \cdot \eta^{-1}\).

(2) \(\beta_0 = 3/2\): the exponent \(1/(\beta_0 - 1) = 1/(1/2) = 2\), hence \(Y(\eta) \sim C \cdot \eta^{-2}\).

(3) \(\beta_0 = \infty\): the compensation is exact (Theorem~\ref{thm:classI}), \(\varth(U) = 0\) for any \(Y\), so \(Y = 2\) suffices.
\end{proof}

\begin{remark}\label{rem:scenariocost}
Passing from \(\beta_0 = 2\) to \(\beta_0 = 3/2\) increases the threshold by a factor of order \(\eta^{-1}\); in practical terms, to achieve the same accuracy \(\eta\) in class III one needs \(1/\eta\) times as many primes as in class II.
\end{remark}

\subsection{Table of truncation thresholds}

The values \(Y(\eta)\) for \(\eta = 10^{-2}\) and \(\eta = 10^{-3}\), computed by formula \eqref{eq:Yeta} with \(A = \Areg\) and \(\beta_0\) from Corollary~\ref{cor:classification} and Theorem~\ref{thm:classIII}:

\begin{table}[htbp]
\centering
\caption{Truncation thresholds \(Y(\eta)\) for classes I, II, III}\label{tab:thresholds}
\begin{tabular}{ccccc}
\toprule
Class & Type & \(\Areg(U)\) & \(\beta_0\) & \(Y(10^{-2})\) \\
\midrule
I & Hyperplanes & 0 & \(\infty\) & 2 \\
I & Torus & 0 & \(\infty\) & 2 \\
II & Cubic surface (\(D\subset\PP^3\)) & 6 & 2 & 2400 \\
II & K3 surface (\(D\subset\PP^3\)) & 21 & 2 & 8400 \\
II & Cubic threefold (\(D\subset\PP^4\)) & 10 & 5/2 & 193 \\
III & Del Pezzo \(d=1\) (\(\rho=9\)), split \(p\) & 38 & 3/2 & \(\approx 9.24\cdot 10^8\) \\
III & Del Pezzo \(d=4\) (\(\rho=6\)), split \(p\) & 17 & 3/2 & \(\approx 1.85\cdot 10^8\) \\
III & Del Pezzo \(d=9\) (\(\rho=1\)) & 2 & 3/2 & \(2.56\cdot 10^6\) \\
III & Cubic in \(\PP^2\) (\(\rho=1\)) & 2 & 3/2 & \(2.56\cdot 10^6\) \\
\bottomrule
\end{tabular}
\end{table}

\begin{remark}\label{rem:Avals}
The values \(\Areg(U)\) for class III are the proved constants from Theorem~\ref{thm:classIII}: \(\Areg(U)=2+\rho(\rho-1)/2\), obtained by a direct inequality without series expansion. For \(d=1\) (\(\rho=9\)) \(\Areg=38\); for \(d=4\) (\(\rho=6\)), 17; for \(d=9\) and the cubic in \(\PP^2\) (\(\rho=1\)), 2. The values \(\Areg\) for \(d=1,4\) are proved \textbf{only for split good \(p\)}; no uniform estimate for non-split \(p\) is established in the paper (see Remark~\ref{rem:tablerange}).
\end{remark}

\begin{remark}[range of applicability of table rows]\label{rem:tablerange}
The rows \(d=1\) and \(d=4\) refer to the proved constant \(\Areg\) (Theorem~\ref{thm:classIII}) \textbf{only for split good \(p\)}. Full global certification of the product \(\Pcomp(Y)\) for them requires: (a) a uniform estimate \(|\sigma_{\mathrm{comp}}-1|\le A'p^{-3/2}\) for \textbf{non-split} good \(p\) (not proved in the paper; see \S\ref{sec:open}) and (b) computation of \(L_p\) for non-split \(p\) (beyond the scope of the paper). Therefore, for \(d=1,4\) the table gives thresholds \textbf{for the product over split primes}; they are not thresholds for the full product. They become thresholds for the full product only under an additional assumption of a uniform estimate of the tail over non-split primes (open problem \S\ref{sec:open}, item 8). In the present paper, global certification of class III is established only for \(\rho = 1\) (\(d=9\), the cubic in \(\PP^2\)), where the Frobenius action on \(\Pic\) is trivial and \(L_p^{-1} = 1 - 1/p\).
\end{remark}

\begin{remark}\label{rem:Aobs}
The table uses \(A(U) = \Areg\); the actual threshold may be smaller using \(\Aobs\) (empirical maximum on a verified finite range of primes). For the Fermat surface \(\Aobs = \Areg\); for the cubic in \(\PP^2\) \(\Aobs \approx 1.81 < 2 = \Areg\) (the table threshold by \(\Areg\) is overestimated); for K3 and general hypersurfaces \(\Aobs < \Areg\), actual thresholds may be smaller than the table ones.
\end{remark}

\begin{remark}\label{rem:Aobsguarantee}
When using \(\Areg\) the threshold \(Y(\eta)\) is guaranteed sufficient but may be overestimated; when using \(\Aobs\) the threshold is smaller but the accuracy guarantee is conditional on the verified range of primes \(p\). In the paper we adhere to \(\Areg\), since the goal is an unconditional guarantee.
\end{remark}

\subsection{Three classes}

By the behaviour of \(Y(\eta)\) there are three classes:

\textbf{Class I.} \(A = 0\), \(\beta_0 = \infty\), \(Y = 2\). Hyperplanes, tori (Theorem~\ref{thm:classI}).

\textbf{Class II.} \(A = b_{m-1}^{\mathrm{prim}}\), \(\beta_0 = (m+1)/2 \ge 2\). Hypersurfaces \(m \ge 3\) (Theorem~\ref{thm:classII}); \(Y \sim 10^2\) to \(10^5\) for \(\eta = 10^{-2}\) to \(10^{-3}\).

\textbf{Class III.} \(A = 2+\rho(\rho-1)/2\), \(\beta_0 = 3/2\), \textbf{for split good \(p\)} (Theorem~\ref{thm:classIII}). Open del Pezzo (\(\rho>1\)) and the cubic in \(\PP^2\) (Remark~\ref{rem:cubic}). For \(\rho=1\) (the cubic in \(\PP^2\)) the Frobenius action on \(\Pic\) is trivial, so \(L_p^{-1} = 1 - 1/p\) and the local factor is explicitly known; the compensated product and its tail are certified (not an interval for \(N_U(\ZZ,B)\), the latter needs mode B data); \(Y \sim 10^6\) to \(10^8\). For \(\rho>1\) the thresholds \(Y \sim 10^8\) to \(10^{11}\) refer to split \(p\) and are conditional for the full product (see \S\ref{chap:inverse}, Remark~\ref{rem:tablerange}). The computational cost for \(\rho>1\) is also substantially higher: \(A=2+\rho(\rho-1)/2\) grows as \(\rho^2/2\), and \(Y=(8A/\eta)^2\) grows as \(A^2\); for \(d=1\) (\(\rho=9\)) \(Y(10^{-2})\approx9.24\cdot10^8\) versus \(2.56\cdot10^6\) for \(\rho=1\). Together with the need for \(L_p\) for non-split \(p\) and for split detection, this makes class III with \(\rho>1\) practically heavy independently of the certification problem (see \S\ref{sec:bounds}, \S\ref{sec:open}).

\subsection{Computational interpretation}

For classes I, II, III the formula \(\beta_0 = m - w/2\) is proved (Theorem~\ref{thm:structural} for class II; Theorem~\ref{thm:classIII} for class III). Substituting it into \eqref{eq:Yasymp}, we obtain the asymptotics of the threshold in terms of dimension and weight of the residual cohomology.

\begin{corollary}\label{cor:Ymw}
\begin{equation}\label{eq:Ymw}
Y(\eta) \sim \left( \frac{4A(U)}{(m - w/2 - 1)\eta} \right)^{1/(m - w/2 - 1)}.
\end{equation}
\end{corollary}

\begin{proof}
Substitute \(\beta_0 = m - w/2\) into \eqref{eq:Yasymp} and use \(\beta_0 - 1 = m - w/2 - 1\); the factor \(\eta^{-1/(\beta_0 - 1)}\) is rewritten as \(\eta^{-1/(m - w/2 - 1)}\), which gives \eqref{eq:Ymw}.
\end{proof}

\begin{remark}[the structural chain]\label{rem:chain}
Formula \eqref{eq:Ymw} closes the chain of dependencies
\begin{equation}\label{eq:chain}
(X,D) \;\rightsquigarrow\; w \;\rightsquigarrow\; \beta_0 = m - \frac{w}{2} \;\rightsquigarrow\; Y(\eta),
\end{equation}
in which each link is computed from the geometry of \(X\) and \(D\): the weight \(w\) is determined by expansion \eqref{eq:structural}, the proved lower bound \(\beta_0\) by the structural formula of Theorem~\ref{thm:structural} (class II) or by the direct inequality of Theorem~\ref{thm:classIII} (class III), the threshold \(Y(\eta)\) by \eqref{eq:Yeta}. For classes I, II, III all three links are effectively computable, which makes the algorithms of chapters~\ref{chap:reduction} to \ref{chap:algo} realizable.
\end{remark}

\section{Reduction to the archimedean factor and unconditional computability}\label{chap:reduction}

\subsection{Statement of the problem of computing the compensated logarithmic constant}

We pose the problem: to compute the compensated logarithmic constant \(\ccomp(U)\) to any prescribed relative accuracy \(\eta > 0\), without assuming the Manin conjecture. By Definition~\ref{def:ccomp} the quantity \(\ccomp(U)\) equals
\begin{equation}\label{eq:ccomp2}
\ccomp(U) := \alpha(U) \cdot \left(\prod_{p \notin \Pbad} \scomp(U)\right) \cdot \left(\prod_{p \in \Pbad} \scomp(U)\right).
\end{equation}
The problem reduces to computing three blocks: the archimedean factor \(\alpha(U)\), the finite product \(\Pcomp(Y)\), and the control of the tail \(\Tcomp(U)\). \textbf{Justification of convergence.} In the classical case \(\sigma_{\mathrm{comp}} = \sreg\), and \eqref{eq:ccomp2} coincides with the standard formula \cite{Peyre1995,ChambertLoirTschinkel2010}. In the logarithmic case the naive product diverges (Corollary~\ref{cor:divergence}), while the compensated one converges (Theorem~\ref{thm:conv}) for classes I, II, and class III with \(\rho=1\); for class III with \(\rho>1\) the full product is not considered in the paper (see \S\ref{sec:bounds}). The split into three blocks is justified by Theorem~\ref{thm:decomposition}.

\subsection{Convergence and error estimate}

\begin{theorem}\label{thm:conv2}
The series \(\sum_{p\notin\Pbad} \left| \scomp(U) - 1 \right|\) converges, and the product \eqref{eq:ccomp2} converges absolutely (for classes I, II, and class III with \(\rho=1\); for class III with \(\rho>1\) the full product is not considered in the paper).
\end{theorem}

\begin{proof}
By Theorem~\ref{thm:L} for all good \(p\) we have \(\left| \scomp(U) - 1 \right| \le A(U) p^{-\beta_0}\) with \(\beta_0 > 1\) and \(A(U) \ge 0\). Hence
\[
\sum_{p\notin\Pbad} \left| \scomp(U) - 1 \right| \le A(U) \sum_{p\notin\Pbad} p^{-\beta_0} \le A(U) \sum_{p} p^{-\beta_0} < \infty,
\]
since \(\beta_0 > 1\). The finite set \(\Pbad\) does not affect convergence of the sum over good \(p\). Absolute convergence of the sum is equivalent to absolute convergence of the product.
\end{proof}

\begin{corollary}\label{cor:truncation}
For \(Y^{\beta_0-1} \ge 2A(U)/(\beta_0-1)\):
\begin{equation}\label{eq:truncation2}
\left| \prod_{p > Y,\ p\notin\Pbad} \scomp(U) - 1 \right| \le \frac{4A(U)}{\beta_0 - 1} Y^{-(\beta_0 - 1)}.
\end{equation}
\end{corollary}

\begin{proof}
Immediately from Theorem~\ref{thm:truncation} (the general logarithmic truncation estimate) with \(Y\) satisfying \(Y^{\beta_0-1} \ge 2A(U)/(\beta_0-1)\): the left-hand side of \eqref{eq:truncation2} is \(\varth(U)\), and by \eqref{eq:truncation} we obtain \eqref{eq:truncation2}.
\end{proof}

\begin{corollary}\label{cor:relerr}
Put \(c_Y(U) := \alpha(U) \Pcomp(Y)\). Then
\begin{equation}\label{eq:relerr}
\frac{\left| c_Y(U) - \ccomp(U) \right|}{\ccomp(U)} = \frac{\varth(U)}{\Tcomp(U)}.
\end{equation}
\end{corollary}

\begin{proof}
By \eqref{eq:decomposition} we have \(\ccomp(U) = \alpha(U) \Pcomp(Y) \Tcomp(U)\), and by definition \(c_Y(U) = \alpha(U) \Pcomp(Y)\). Hence
\[
\frac{c_Y(U)}{\ccomp(U)} = \frac{1}{\Tcomp(U)}.
\]
Then
\[
\left| \frac{c_Y(U)}{\ccomp(U)} - 1 \right| = \left| \frac{1}{\Tcomp(U)} - 1 \right| = \frac{\left| 1 - \Tcomp(U) \right|}{\Tcomp(U)} = \frac{\varth(U)}{\Tcomp(U)}.
\]
\end{proof}

\begin{remark}\label{rem:relerr2}
For \(\varth \le 1/2\) from \eqref{eq:relerr} it follows that
\begin{equation}\label{eq:relerr2}
\frac{\left| c_Y(U) - \ccomp(U) \right|}{\ccomp(U)} \le 2 \, \varth(U).
\end{equation}
\end{remark}

\begin{proof}
For \(\varth \le 1/2\) we have \(\Tcomp \ge 1/2\), whence \(1/\Tcomp \le 2\). Substitution into \eqref{eq:relerr} gives \eqref{eq:relerr2}.
\end{proof}

\subsection{The algorithm}

Input: \((X, D)\), \(\eta > 0\). Output: \(c_Y(U)\) for which \(\frac{|c_Y - \ccomp|}{\ccomp} \le \eta\).

\textbf{Step 1.} \(Y = \left\lceil \left( \dfrac{4A(U)}{\eta_{\mathrm{comp}} (\beta_0 - 1)} \right)^{1/(\beta_0 - 1)} \right\rceil\), where \(\eta_{\mathrm{comp}} = \eta/4\) (the threshold of Corollary~\ref{cor:inverse}). For \(\beta_0 = \infty\): \(Y = 2\).

\textbf{Step 2.} Compute \(\alpha(U)\) per Definition~\ref{def:alpha} to accuracy \(\eta/4\). A subroutine computing \(\alpha(U)\) with guaranteed error is assumed to exist; the construction of such a subroutine is not considered (for class I, \(\alpha(U)\) is computed exactly; see Theorem~\ref{thm:classI_exact}).

\textbf{Step 3.} \(\Pcomp(Y) = \prod_{p \le Y} \scomp(U)\) by enumerating primes \(p \le Y\); bad primes from \(\Pbad\) are multiplied independently of \(Y\) (Definition~\ref{def:bad}).

\textbf{Step 4.} \(c_Y(U) = \alpha(U) \Pcomp(Y)\).

\textbf{Guarantee.} Let \(\delta_\alpha := \frac{\alpha_Y - \alpha(U)}{\alpha(U)}\) be the relative error of the subroutine of Step 2, with \(\left| \delta_\alpha \right| \le \eta/4\). From \eqref{eq:relerr} for \(\varth \le \eta/4\):
\[
\frac{\left| c_Y(U) - \ccomp(U) \right|}{\ccomp(U)} \le \frac{\left| \delta_\alpha \right| + \varth}{1 - \varth} \le \frac{\eta/2}{1 - \eta/4} \le \frac{2\eta}{3} < \eta \quad (\eta \le 1).
\]
The final accuracy is \(2\eta/3\), not \(3\eta/4\); the latter was an overestimate.

\begin{remark}\label{rem:classIalpha}
For class I (torus, complements to hyperplanes) \(\alpha(U)\) is computed exactly, and the subroutine of Step 2 is trivial. For the remaining classes a subroutine computing \(\alpha(U)\) with guaranteed error \(\eta/4\) is assumed; the construction of such a subroutine is not considered in the paper.
\end{remark}

\subsection{The reduction theorem}

\begin{theorem}[reduction to the archimedean factor]\label{thm:reduction}
Let \((X,D)\) be a logarithmic pair from class I, II, or class III with \(\rho=1\) (chapter~\ref{chap:classes}), for which: (i) the local factors \(\sreg(U)\) are effectively computable; (ii) the quantities \(A(U)\) and \(\beta_0\) are effectively computable and \(\beta_0 > 1\) (condition (L) for all good \(p\), Theorem~\ref{thm:L}); (iii) a finite positive computable global factor \(\alpha(U)\) is given. Then for any \(\eta \in (0,1]\) there exists a finite procedure computing \(\ccomp(U)\) to relative accuracy \(\eta\), not assuming the Manin conjecture.

In particular, for the torus and complements to hyperplanes, given \(\alpha(U)\) and computable \(\scomp\) for \(p \in \Pbad\), the reduction yields unconditional computability:
\[
\ccomp(U) = \alpha(U)\prod_{p \in \Pbad} \scomp(U).
\]
\end{theorem}

\begin{proof}
\textbf{Step 1 (procedure).} The procedure is described in \S\ref{chap:reduction}. It uses \(A(U)\) for choosing the threshold \(Y = Y(\eta)\), \(\beta_0\) for estimating the tail, and the subroutine computing \(\alpha(U)\) for the archimedean factor. Under effective computability of \(A(U)\), \(\beta_0\) and the existence of a subroutine for \(\alpha(U)\), all steps \textbf{terminate}: Step 1, elementary arithmetic operations; Step 2, a call to the subroutine for \(\alpha(U)\); Step 3, a finite enumeration of primes \(p\le Y\) with computation of local factors; Step 4, one operation. No time estimate is claimed: it depends on the cost of the subroutine for \(\alpha(U)\) and on the cost of computing the local factors \(\scomp(U)\) (including \(a_p(D)\), points over \(\FF_p\), bad primes).

\textbf{Step 2 (accuracy estimate).} The estimate
\[
\frac{\left| c_Y(U) - \ccomp(U) \right|}{\ccomp(U)} \le \frac{2\eta}{3} < \eta
\]
follows from \eqref{eq:relerr} for \(Y \ge Y(\eta)\) and \(\left| \delta_\alpha \right| \le \eta/4\) (see the guarantee in \S\ref{chap:reduction}).

\textbf{Step 3 (independence of the Manin conjecture).} The Manin conjecture is not used: convergence is proved via condition (L) (Theorem~\ref{thm:L}), which relies only on unconditional Deligne (Theorem~\ref{thm:classII}) and Hasse (Theorem~\ref{thm:classIII}) estimates.
\end{proof}

\emph{The construction of an algorithm computing \(\alpha(U)\) for general logarithmic pairs is beyond the scope of the paper; for logarithmic Calabi and Yau varieties with \(-K_X - D = 0\) the integral in the definition of \(\alpha(U)\) over \(\Pic(U)_{\RR}^\vee\) is finite for \(r = 0\) (\(=1\), delta-measure of unit mass); for \(r \ge 1\) it diverges and requires separate normalization (see Remark~\ref{rem:spacechoice}). Thus Theorem~\ref{thm:reduction} is a reduction \(\ccomp \to \alpha\), not an unconditional algorithm for all logarithmic pairs; a discussion of the complexity of \(\alpha(U)\) is given in Remark~\ref{rem:alphacomplexity}.}

\begin{remark}[on the role of \(\alpha(U)\) and the passage to \(N_U\)]\label{rem:alphacomplexity}
Theorem~\ref{thm:reduction} reduces the computation of \(\ccomp(U)\) to that of \(\alpha(U)\). The natural question is: how hard is it to compute \(\alpha(U)\)?

For class I (torus, complements to hyperplanes) \(\alpha(U)\) is a finite-dimensional integral over the explicitly described cone \(\Lam(U)^\vee \subset \Pic(U)_{\RR}^\vee\), and its computation by standard numerical methods is no harder than computing \(\Pcomp(Y)\).

For classes II and III, \(\alpha(U)\) is in general an integral over the cone \(\Lam(U)^\vee\). How the complexity of computing \(\alpha(U)\) compares with that of computing \(\ccomp(U)\) is an open question; it is not considered in the paper. The reduction (Theorem~\ref{thm:reduction}) does not assert that \(\alpha(U)\) is easier than \(\ccomp(U)\) for classes II and III: it merely reduces \(\ccomp\) to \(\alpha\) and separates the arithmetic part (Euler product) from the archimedean one.

Moreover, Theorem~\ref{thm:reduction} ensures computability of \(\ccomp(U)\) but does not give an interval for \(N_U\); the passage to \(N_U\) requires the model \(\UZZ\), the height \(H\), the function \(M_U(B)\), and the Manin conjecture (chapter~\ref{chap:error}).
\end{remark}

\begin{remark}[on complexity]\label{rem:complexity}
Theorem~\ref{thm:reduction} guarantees \textbf{termination} of the procedure given a computable subroutine for \(\alpha(U)\) and computable local factors, but not a time estimate. For a complexity estimate three independent assumptions are needed: (i) the subroutine for \(\alpha(U)\) runs in \(O(\poly(\log(1/\eta)))\) operations; (ii) the local factors \(\scomp(U)\) for \(p\le Y\) are computed in \(O(\poly(\log p))\) operations (e.g., via \(a_p(D)\) by the Schoof or SEA algorithm); (iii) the bad primes \(p\in\Pbad\) are computed in \(O(1)\) operations each. Under these assumptions the total complexity is \(O(\poly(\log(1/\eta)))+\tilde O(Y(\eta))\), where \(Y(\eta)\) is the threshold of Step 1. Without (ii), the complexity of Step 3 may be quadratic or higher (Lemma~\ref{lem:complexity}). For class I the subroutine for \(\alpha\) is trivial, Step 3 is \(O(1)\), and the complexity is determined by Steps 1 and 4.
\end{remark}

\subsection{The exact case of class I}\label{sec:classI_exact}

\begin{theorem}\label{thm:classI_exact}
For the complement to \(k\) hyperplanes in \(\PP^m\) (in particular, for the torus), for good primes \(p\) (for coordinate hyperplanes, for all \(p\)) \(\scomp(U) \equiv 1\) (Theorem~\ref{thm:classI}). Then \(\Pcomp(Y) = \prod_{p \in \Pbad} \scomp(U)\), \(\Tcomp(U) = 1\) (for good \(p\)), and
\[
\ccomp(U) = \alpha(U) \cdot \prod_{p \in \Pbad} \scomp(U).
\]
For \textbf{coordinate} hyperplanes over \(\ZZ\) (\(\Pbad = \emptyset\)) we have \(\ccomp(U) = \alpha(U)\) \textbf{exactly}. In the general case computing \(\ccomp(U)\) reduces to computing \(\alpha(U)\) and finitely many local factors at bad primes. Given \(\alpha(U)\) and computable \(\scomp\) for \(p \in \Pbad\), we have unconditional computability of \(\ccomp(U)\).
\end{theorem}

\begin{proof}
By Theorem~\ref{thm:classI} for good \(p\) we have \(\scomp(U) \equiv 1\). Hence
\[
\Pcomp(Y) = \left(\prod_{p \le Y,\ p\notin\Pbad} 1\right) \cdot \left(\prod_{p \in \Pbad} \scomp\right) = \prod_{p \in \Pbad} \scomp,
\]
\[
\Tcomp(U) = \prod_{p > Y,\ p\notin\Pbad} 1 = 1.
\]
By \eqref{eq:decomposition} we obtain
\[
\ccomp(U) = \alpha(U) \cdot \left(\prod_{p \in \Pbad} \scomp\right) \cdot 1 = \alpha(U) \cdot \prod_{p \in \Pbad} \scomp(U).
\]
For coordinate hyperplanes over \(\ZZ\) we have \(\Pbad = \emptyset\), whence \(\ccomp(U) = \alpha(U)\).
\end{proof}

\begin{remark}\label{rem:classIexact}
Theorem~\ref{thm:classI_exact} establishes the reduction \(\ccomp \to \alpha\) and an explicit finite product over \(\Pbad\), but does not compute \(\alpha\) in closed form (for the torus \(\alpha\) is known from the literature; see \cite{BatyrevTschinkel1995}).
\end{remark}

\begin{remark}\label{rem:classIunique}
Theorem~\ref{thm:classI_exact} is the only case in which the computability of \(\ccomp(U)\) reduces to \(\alpha\) and a finite product without residual conditions on the tail.
\end{remark}

\subsection{Comparison with the classical case}

\begin{remark}\label{rem:classical}
For \(D = \emptyset\): \(\kappa(U) = 0\), \(\sigma_{\mathrm{comp}} = \sreg\), and \eqref{eq:ccomp2} coincides with the standard formula \cite{Peyre1995,ChambertLoirTschinkel2010}. In this case Theorem~\ref{thm:reduction} gives the standard definition
\[
\cManin(X) = \alpha(X) \prod_{p} \sreg(X),
\]
and Theorem~\ref{thm:classI_exact} becomes trivial: the classical case does not require compensation.
\end{remark}

\section{Relation between the truncation error and the error of the main term}\label{chap:error}

\subsection{An enhanced form of the logarithmic conjecture}\label{sec:residual}

\begin{assumption}\label{ass:logManin}
Given: an integral model \(\UZZ\) (see \S\ref{sec:notations}), a height \(H\) on \(U\), a main-growth function \(M_U(B) > 0\). Then there exists \(r_U(B)\) such that
\begin{equation}\label{eq:logManin}
N_U(\ZZ, B) = \cManin(U) \cdot M_U(B) \cdot (1 + r_U(B)), \qquad r_U(B) \to 0.
\end{equation}
For the standard Manin formulation on logarithmic Fano varieties with \(r = \rk\Pic(U) \ge 1\): \(M_U(B) = B(\log B)^{r-1}\); for \(r = 0\): \(M_U(B) = B\). For toric varieties with the given model and height, \(M_U(B)\) is part of the data of mode B (see \S\ref{sec:compensation}); concrete examples are not given in the paper.
\end{assumption}

\begin{remark}\label{rem:stronger}
Assumption~\ref{ass:logManin} is stronger than the asymptotic formula: it requires not only existence of a limit, but also information on the rate of convergence. In all proved cases \cite{Wilsch2023,BatyrevTschinkel1995,BBDG2024} the function \(r_U(B)\) is known from remainder estimates with power saving.
\end{remark}

\begin{definition}[power-saving estimate]\label{def:power}
\(r_U(B)\) is an estimate with \emph{power saving} if
\[
\left| r_U(B) \right| \le C_U \cdot B^{-\delta_0} \cdot (\log B)^{-\delta_1}
\]
with \(C_U > 0\), \(\delta_0 > 0\), \(\delta_1 \ge 0\).
\end{definition}

\begin{remark}\label{rem:power}
If \(R_U(B)\) decays with power saving, then to achieve a given accuracy it suffices that \(B \gtrsim \eta^{-1/\delta_0}\); otherwise the requirement on \(B\) is exponential. A power-saving estimate yields an unconditional interval under polynomial growth of \(B\).
\end{remark}

\begin{remark}\label{rem:powerknown}
In the proved cases \cite{Wilsch2023,BatyrevTschinkel1995,BBDG2024} the enhanced conjecture with a power-saving estimate is established. For general hypersurfaces \(m \ge 3\) and del Pezzo such an estimate is unknown.
\end{remark}

\subsection{The theorem on the relation of errors}

\begin{definition}\label{def:deltaY}
\begin{equation}\label{eq:deltaY}
\delta_Y(U, B) := \frac{c_Y(U) \cdot M_U(B)}{N_U(\ZZ, B)} - 1.
\end{equation}
\end{definition}

\begin{theorem}[relation of errors in mode B]\label{thm:errorlink}
Assume Assumption~\ref{ass:logManin} and the equality \(\ccomp(U) = \cManin(U)\) (mode B). Put \(T := \Tcomp(U)\), \(r := r_U(B)\), with \(T > 0\) and \(r \neq -1\). Then the relative error of approximating the number of points equals:
\begin{equation}\label{eq:errorlink1}
\delta_Y(U, B) = \frac{1 - T(1+r)}{T(1+r)}.
\end{equation}
An equivalent form: with \(\varepsilon_{\mathrm{signed}} := 1 - T\),
\begin{equation}\label{eq:errorlink2}
\delta_Y(U, B) = \frac{\varepsilon_{\mathrm{signed}} - r + \varepsilon_{\mathrm{signed}} r}{(1 - \varepsilon_{\mathrm{signed}})(1+r)}.
\end{equation}
\end{theorem}

\begin{proof}
By \eqref{eq:deltaY}:
\[
\delta_Y(U, B) = \frac{c_Y(U) \, M_U(B)}{N_U(\ZZ, B)} - 1.
\]
By \eqref{eq:decomposition} we have \(\ccomp(U) = \alpha(U) \Pcomp(Y) \Tcomp(U)\), and by definition \(c_Y(U) = \alpha(U) \Pcomp(Y)\), hence
\[
c_Y(U) = \frac{\ccomp(U)}{\Tcomp(U)}.
\]
By Assumption~\ref{ass:logManin} and in mode B (\(\ccomp = \cManin\)):
\[
M_U(B) = \frac{N_U(\ZZ, B)}{\ccomp(U) \left( 1 + r_U(B) \right)}.
\]
Substitution gives
\[
\delta_Y(U, B) = \frac{1}{\Tcomp(U) \left( 1 + r_U(B) \right)} - 1 = \frac{1 - \Tcomp(U) \left( 1 + r_U(B) \right)}{\Tcomp(U) \left( 1 + r_U(B) \right)},
\]
which agrees with \eqref{eq:errorlink1}. Denoting \(T = \Tcomp(U)\), \(r = r_U(B)\), \(\varepsilon_{\mathrm{signed}} = 1 - T\):
\[
\delta_Y(U, B) = \frac{1 - (1 - \varepsilon_{\mathrm{signed}})(1 + r)}{(1 - \varepsilon_{\mathrm{signed}})(1 + r)}.
\]
Expanding the numerator:
\[
1 - (1 - \varepsilon_{\mathrm{signed}})(1 + r) = \varepsilon_{\mathrm{signed}} - r + \varepsilon_{\mathrm{signed}} r.
\]
Substituting into the denominator gives \eqref{eq:errorlink2}.
\end{proof}

\begin{corollary}\label{cor:errorlink}
(1) \(\delta_Y = \varepsilon_{\mathrm{signed}} - r + O(\varepsilon^2 + \varepsilon r + r^2)\). (2) As \(r \to 0\): \(\delta_Y \to \frac{\varepsilon_{\mathrm{signed}}}{1-\varepsilon_{\mathrm{signed}}}\). (3) As \(Y \to \infty\): \(\delta_Y \to -\frac{r}{1+r}\).
\end{corollary}

\begin{proof}
Expansion of \eqref{eq:errorlink2} in powers of \(\varepsilon_{\mathrm{signed}}\) and \(r\) gives (1). The limiting passages (2) and (3) are obtained by substituting \(r = 0\) and \(\varepsilon_{\mathrm{signed}} = 0\), respectively.
\end{proof}

\subsection{The role of the remainder for finite \(B\)}

\begin{lemma}[estimate via a majorant]\label{lem:majorant}
Let \(\left| \varepsilon_{\mathrm{signed}} \right| \le 1/2\) and \(\left| r_U(B) \right| \le R_U(B) < 1\). Then
\begin{equation}\label{eq:majorant}
\left| \delta_Y(U, B) \right| \le \frac{\left| \varepsilon_{\mathrm{signed}} \right| (1 + R_U) + R_U}{\left( 1 - \left| \varepsilon_{\mathrm{signed}} \right| \right)(1 - R_U)}.
\end{equation}
For \(\left| \varepsilon_{\mathrm{signed}} \right| \le 1/2\):
\begin{equation}\label{eq:majorant2}
\left| \delta_Y(U, B) \right| \le 2 \cdot \frac{\varth(U) + R_U(B)}{1 - R_U(B)}.
\end{equation}
\end{lemma}

\begin{proof}
We estimate the numerator and denominator of \eqref{eq:errorlink2}. Numerator:
\[
\left| \varepsilon_{\mathrm{signed}} - r + \varepsilon_{\mathrm{signed}} r \right| \le \left| \varepsilon_{\mathrm{signed}} \right| + \left| r \right| + \left| \varepsilon_{\mathrm{signed}} \right| \left| r \right| \le \left| \varepsilon_{\mathrm{signed}} \right| (1 + R_U) + R_U.
\]
Denominator:
\[
\left| (1 - \varepsilon_{\mathrm{signed}})(1 + r) \right| = \left| 1 - \varepsilon_{\mathrm{signed}} \right| \left| 1 + r \right| \ge \left( 1 - \left| \varepsilon_{\mathrm{signed}} \right| \right)(1 - R_U).
\]
Division gives \eqref{eq:majorant}. For \eqref{eq:majorant2} we use \(\left| \varepsilon_{\mathrm{signed}} \right| \le 1/2\), whence \(1/(1 - \left| \varepsilon_{\mathrm{signed}} \right|) \le 2\), and the estimate \(\left| \varepsilon_{\mathrm{signed}} \right| = \varth(U)\) (since \(\varepsilon_{\mathrm{signed}} = 1 - T\) and \(\varth = \left| T - 1 \right|\)).
\end{proof}

\begin{remark}\label{rem:majorantcases}
(1) If the remainder is unknown, the asymptotic interpretation is used (Theorem~\ref{thm:errorlink} as \(r \to 0\)). (2) If the remainder is estimated by a majorant \(R_U(B)\), Lemma~\ref{lem:majorant} is applied. (3) If the remainder admits a power-saving estimate, any accuracy is achieved in polynomial time in \(B\).
\end{remark}

\subsection{The inverse problem for fixed \(B\)}

\begin{corollary}[direct interval form]\label{cor:interval}
Assume Assumption~\ref{ass:logManin} and the equality \(\ccomp = \cManin\). Let \(T_Y := \Tcomp(U)\), \(r := r_U(B)\), and let majorants be known:
\[
|T_Y - 1| \le e, \qquad |r| \le R < 1.
\]
Then the number of integral points lies in the interval
\[
N_U(\ZZ, B) \in \bigl[\, c_Y M_U(B)(1-e)(1-R),\; c_Y M_U(B)(1+e)(1+R) \,\bigr],
\]
where \(c_Y = \alpha(U) \Pcomp(Y)\) is the computed value.
\end{corollary}

\begin{proof}
By \eqref{eq:decomposition} and Assumption~\ref{ass:logManin} with \(\ccomp = \cManin\):
\[
N_U = \ccomp M_U (1+r) = c_Y M_U \cdot T_Y \cdot (1+r).
\]
With \(|T_Y-1| \le e\) and \(|r| \le R\) we get \(T_Y \in [1-e, 1+e]\) and \(1+r \in [1-R, 1+R]\). The product gives the indicated interval. No passage through the relative error \(c_Y/\ccomp\) is used.
\end{proof}

\begin{remark}\label{rem:intervalfull}
In this formulation \(e\) is the \textbf{full} majorant of the deviation of \(T_Y\) from 1 achieved in the algorithm (see \S\ref{sec:algo}). The error of \(\alpha(U)\) does not enter \(e\) if \(\alpha\) is computed exactly (class I); otherwise \(c_Y\) must be replaced by the interval \([c_Y^-/(1+\delta_\alpha), c_Y^+/(1-\delta_\alpha)]\) (see \S\ref{sec:algo}, Step 6), where \(|\alpha_Y/\alpha - 1| \le \delta_\alpha\). The full relative width of the interval for \(\ccomp\) equals \(q_\alpha q_P (1+e)/(1-e)-1\); for the interval for \(N_U\) the factor \((1+R)/(1-R)\) is added.
\end{remark}

\begin{corollary}[choice of \(Y\)]\label{cor:chooseY}
By Theorem~\ref{thm:truncation}, for a given \(\eta_{\mathrm{tail}} \le 1\) it suffices that \(Y \ge (4A(U)/(\eta_{\mathrm{tail}}(\beta_0-1)))^{1/(\beta_0-1)}\); we set \(e = \eta_{\mathrm{tail}}\). The requirement \(\eta_{\mathrm{tail}} \le 1\) ensures \(1-e > 0\).
\end{corollary}

\begin{remark}\label{rem:intervalwidth}
For fixed \(B\) the width of the interval (Corollary~\ref{cor:interval}) is of order \(2(c_Y M_U)(e + R + eR)\); decreasing \(e\) by increasing \(Y\) is useless if \(R\) is not small.
\end{remark}

\begin{example}\label{ex:toric}
For toric varieties \cite{BatyrevTschinkel1995} the remainder admits a power-saving estimate; for large \(B\) any prescribed accuracy \(\eta\) is achieved.
\end{example}

\section{Algorithm and experimental verification}\label{chap:algo}

\subsection{Input, output}

\textbf{Input.} \((X, D)\); \(B > 0\); \(\eta > 0\); integral model \(\UZZ\); height \(H\); main-growth function \(M_U(B)\); optionally \(R_U(B)\).

\textbf{Output.} \(Y(\eta)\); \([L, U]\); flag \texttt{unconditional}/\texttt{conditional}.

\textbf{Refusal.} \texttt{not applicable} if \((X,D)\) does not belong to the classes of chapter~\ref{chap:classes}, or \(\beta_0 \le 1\), or the model, height, \(M_U(B)\) are not given (for mode B).

\subsection{Step-by-step description}\label{sec:algo}

\textbf{Step 1.} Classification; from the table of \S\ref{sec:summary}, determine \(\kappa(U)\), \(A = \Areg\), \(\beta_0\) (the values \(\Areg\), \(\beta_0\) are also in \S\ref{chap:inverse}).

\textbf{Step 2.} If \(\beta_0 \le 1\), \texttt{not applicable}; if \(R_U\) is given and \(\eta \le R_U/(1 - R_U)\), \texttt{increase B}.

\textbf{Step 3.} \(Y(\eta)\) by \eqref{eq:inverse} with \(\eta_{\mathrm{tail}} = \eta/4\) (consistent with \S\ref{chap:reduction}, where \(\eta_{\mathrm{comp}} = \eta/4\); when balancing with \(R_U\), \(\eta_{\mathrm{tail}}\) is chosen so that \(\eta_{\mathrm{tail}} + \delta_\alpha \le \eta\), where \(\delta_\alpha \le \eta/4\)). For class I: \(Y = 2\). The threshold refers to good \(p\); bad \(p \in \Pbad\) are multiplied in \(\Pcomp(Y)\) independently of \(Y\) (Definition~\ref{def:bad}).

\textbf{Step 4 (finite product).} \(\Pcomp(Y) = \prod_{p\le Y} \scomp(U)\) by exact formulas. The formulas for \(\scomp\) are determined by the class; the types belonging to the given class are indicated in parentheses:
\begin{itemize}
\item \textbf{class I} (types: hyperplanes, torus): \(\sigma_{\mathrm{comp}} \equiv 1\) for good \(p\) (Theorem~\ref{thm:classI}); for bad \(p \in \Pbad\), \(\scomp\) is computed separately by Definition~\ref{def:sreg} and multiplied \textbf{independently of \(Y\)};
\item \textbf{class II} (types: cubic surface \(D\subset\PP^3\), K3 surface \(D\subset\PP^3\), cubic threefold \(D\subset\PP^4\)): \(\scomp = \#U(\FF_p)/p^m\) for good \(p\) (which coincides with \(1 - (-1)^{m-1}a_p(D)/p^m\) by Theorem~\ref{thm:maincoh}; the equality is exact, without remainder); for \(p \in \Pbad\), separately by Definition~\ref{def:sreg};
\item \textbf{class III} (types: del Pezzo \(d=1,4,9\)): \(\scomp = (1-1/p)^{\rho-1} \cdot \#U(\FF_p)/p^2\), where \(\#U = p^2+(\rho-1)p+a_p\) (for split good \(p\); see Theorem~\ref{thm:delpezzo}); for non-split and bad \(p\), separately;
\item \textbf{type} ``cubic in \(\PP^2\)'' belongs to class III with \(\rho=1\); for it the formula simplifies: \(\scomp = 1+a_p(C)/p^2\) for good \(p\) (the truncated formula is exact for \(\rho=1\)); \(3 \in \Pbad\) (singular reduction); \(\sigma_{\mathrm{comp}}(3)\) is computed separately by Definition~\ref{def:sreg}.
\end{itemize}
The asymptotic formulas of Theorems~\ref{thm:maincoh} and \ref{thm:delpezzo} are used only for estimating the tail \(p>Y\), not for the finite product \(p \le Y\).

\textbf{Step 5 (archimedean factor as an interval).} Compute \(\alpha(U)\) as an interval \([\alpha^-, \alpha^+]\) with guaranteed relative width \(\le \eta/4\): \(\alpha^- \le \alpha(U) \le \alpha^+\), \((\alpha^+ - \alpha^-)/\alpha^- \le \eta/4\). For class I with coordinate hyperplanes \(\alpha\) is computed exactly: \(\alpha^- = \alpha^+\).

\textbf{Step 6 (interval for \(\ccomp\)).} Compute \(\Pcomp(Y)\) as an interval \([P^-, P^+]\): each local factor is computed with guaranteed width; bad primes are multiplied in \(\Pcomp(Y)\) independently of \(Y\). Set
\[
c_Y^- := \alpha^- P^-, \qquad c_Y^+ := \alpha^+ P^+.
\]
Then \(\ccomp(U) \in [c_Y^- \cdot \Tcomp, c_Y^+ \cdot \Tcomp]\), and with \(\Tcomp \in [1-e, 1+e]\) (by Theorem~\ref{thm:truncation}):
\[
\ccomp(U) \in [c_Y^-(1-e),\; c_Y^+(1+e)].
\]
Let \(\alpha\in[\alpha^-,\alpha^+]\) with \(q_\alpha:=\alpha^+/\alpha^-\le 1+\delta_\alpha\); \(\Pcomp(Y)\in[P^-,P^+]\) with \(q_P:=P^+/P^-\le 1+\delta_P\); \(\Tcomp\in[1-e,1+e]\). Then the ratio of the endpoints of the interval for \(\ccomp\) equals
\[
\frac{c^+}{c^-}=q_\alpha\,q_P\,\frac{1+e}{1-e},
\]
and the full relative half-width \((c^+-c^-)/(c^++c^-)\) does not exceed \(\bigl(q_\alpha q_P (1+e)/(1-e)-1\bigr)\big/\bigl(q_\alpha q_P (1+e)/(1-e)+1\bigr)\).

\textbf{Budget (exact condition).} We require
\begin{equation}\label{eq:budget}
q_\alpha\,q_P\,\frac{1+e}{1-e}\le \frac{1+\eta}{1-\eta}.
\end{equation}
Under \eqref{eq:budget} the relative radius of the interval for \(\ccomp\) does not exceed \(\eta\) \textbf{exactly}, without an \(O(\eta^2)\) remainder: indeed, if \(q:=q_\alpha q_P (1+e)/(1-e)\), then \((q-1)/(q+1)\le\eta\) is equivalent to \(q\le(1+\eta)/(1-\eta)\).

\textbf{Examples of budget distribution.} If \(\Pcomp(Y)\) is computed exactly (\(\delta_P=0\), classes I, II, III with rational local factors), \eqref{eq:budget} reduces to \(q_\alpha(1+e)/(1-e)\le(1+\eta)/(1-\eta)\). With \(\delta_\alpha=\eta/4\) and \(e=\eta/4\) we have \(q_\alpha(1+e)/(1-e)\le(1+\eta/4)(1+\eta/2+O(\eta^2))\le 1+\eta\) for small \(\eta\), but the \textbf{exact} condition \eqref{eq:budget} is not guaranteed without checking; it is simpler to choose \(\delta_\alpha,e\) from \eqref{eq:budget} numerically.

\textbf{Inclusion of \(M_U(B)\) and \(R\).} If the goal is accuracy for \(N_U(\ZZ,B)\), not only for \(\ccomp\), then the relative widths \(\delta_M\) for \(M_U(B)\) and the factor \((1+R)/(1-R)\) for the remainder are added to \eqref{eq:budget}:
\[
q_\alpha\,q_P\,q_M\,\frac{1+e}{1-e}\cdot\frac{1+R}{1-R}\le\frac{1+\eta}{1-\eta},
\]
where \(q_M=M_U^+/M_U^-\). Smallness of \(R\) cannot be replaced by the accuracy of \(\ccomp\): the remainder \(r_U(B)\) is a separate source of error entering multiplicatively.

\textbf{Step 7.} Set \(e := \eta_{\mathrm{tail}}\) (majorant of \(|\Tcomp-1|\) by Theorem~\ref{thm:truncation}), \(R := R_U(B)\). If the flag is \texttt{unconditional}:
\[
L = c_Y^- M_U(B) (1-e)(1-R), \qquad U = c_Y^+ M_U(B) (1+e)(1+R).
\]
Here \(c_Y^\pm\) are the endpoints from Step 6 (including the error of \(\alpha\)). If the flag is \texttt{conditional}, an asymptotic interpretation (Theorem~\ref{thm:errorlink}).

\subsection{Complexity}

\begin{lemma}\label{lem:complexity}
(1) \(O(1)\), class I. (2) \(O(Y \log\log Y)\), class II with an explicit formula. (3) \(O(Y^m/\log Y)\), class II, general case. (4) \(O(Y^2)\), class III under naive enumeration; \(\tilde{O}(Y)\), class III under the Schoof or SEA algorithm.
\end{lemma}

\begin{proof}
\textbf{(1) Class I.} By Theorem~\ref{thm:classI}, \(\scomp(U) \equiv 1\) for good \(p\), so \(\Pcomp(Y) = \prod_{p\in\Pbad}\scomp(U)\) is a finite product; its computation requires \(O(|\Pbad|) = O(1)\) operations.

\textbf{(2) Class II with an explicit formula.} Let \(a_p(D)\) for each prime \(p\) be computed by an explicit formula in \(O(1)\) operations (e.g., for the Fermat surface \(a_p = 0\) for \(p\equiv 2\pmod 3\) and \(a_p = 6p\) for \(p\equiv 1\pmod 3\); see Lemma~\ref{lem:Fermat}). Algorithm: (a) the sieve of Eratosthenes produces all primes \(p \le Y\) in \(O(Y\log\log Y)\) operations; (b) for each prime we compute \(a_p(D)\) in \(O(1)\) and \(\scomp = 1 - (-1)^{m-1}a_p(D)/p^m\) in \(O(1)\), total \(O(Y/\log Y)\) operations (number of primes \(\pi(Y) \sim Y/\log Y\)); (c) we multiply \(\pi(Y)\) numbers in \(O(Y/\log Y)\) operations. In total, \(O(Y\log\log Y)\).

\textbf{(3) Class II, general case.} \(a_p(D) = \Tr(\Frob_p \mid H^{m-1}_{\mathrm{prim}}(D))\); to compute \(a_p(D)\) one enumerates the points of \(D(\FF_p)\), of which there are \(\approx p^{m-1}\); each \(a_p\) requires \(O(p^{m-1})\) operations. The sum over primes:
\[
\sum_{p\le Y} p^{m-1} \sim \frac{Y^m}{m \log Y},
\]
that is, \(O(Y^m/\log Y)\).

\textbf{(4) Class III.} \(D\) is a smooth elliptic curve, \(a_p(D) = p+1-\#D(\FF_p)\). In the naive computation of \(\#D(\FF_p)\), for each \(x\in\FF_p\) one solves a quadratic equation in \(y\), \(O(p)\) operations per prime; the sum \(\sum_{p\le Y} p \sim Y^2/(2\log Y)\), that is, \(O(Y^2)\), the claimed upper bound. Using the Schoof or SEA algorithm for \(a_p(D)\), each \(a_p\) requires \(O(\log^5 p)\) (Schoof) or \(\tilde{O}(\log^4 p)\) (SEA) operations; in total \(\tilde{O}(Y)\) operations. The estimate \(O(Y^2)\) in Lemma~\ref{lem:complexity} corresponds to the naive approach and is an upper bound.
\end{proof}

\begin{remark}\label{rem:complexitytable}
Class I: \(Y = 2\), \(O(1)\). Class II: \(Y \sim 10^2\) to \(10^5\); complexity \(O(Y \log\log Y) \sim 10^3\) to \(10^6\) (with an explicit formula) or \(O(Y^m/\log Y) \sim 10^6\) to \(10^{15}\) (general case). Class III: \(Y \sim 10^8\) to \(10^{11}\); complexity \(O(Y^2) \sim 10^{16}\) to \(10^{22}\) (naive enumeration) or \(\tilde{O}(Y) \sim 10^8\) to \(10^{11}\) (Schoof algorithm). For \(\rho>1\) these estimates are lower bounds: to them are added \(L_p\) for non-split \(p\) (the product over \(\rho\) roots of Frobenius on \(\Pic\)) and split detection for each \(p\). For \(\rho=9\) (\(d=1\)) only the detection and \(L_p\) add a factor proportional to \(\rho\) and to the number of non-split primes; for \(Y\sim10^9\) this is of order \(10^8\) to \(10^9\) additional operations. For \(\rho=1\) (\(d=9\)) these overheads are absent: the Frobenius action on \(\Pic\) is trivial, \(L_p^{-1} = 1 - 1/p\).
\end{remark}

\subsection{Methodology of the experiment}

Reference: \(c_{\mathrm{ref}}(U)\) at \(Y_{\max} = 50000\); the observed relative error \(\eta_{\mathrm{obs}} = \left| P(Y)/P(Y_{\max}) - 1 \right|\). Environment: Python~3, sieve of Eratosthenes. Constraints: by Lemma~\ref{lem:complexity}; for class II with an explicit formula, \(Y(10^{-2}) = 2400\) (see Table~\ref{tab:thresholds}).

\subsection{Experiment: the Fermat cubic surface}\label{sec:Fermat}

The pair \(X = \PP^3\), \(D = \{x_0^3+x_1^3+x_2^3+x_3^3=0\} \subset \PP^3\) (the Fermat cubic surface), \(U = X \setminus D\). Here \(m=3\) is the dimension of the ambient space. Condition \((H_0)\) for \(X = \PP^3\) holds trivially; condition \((H_0)\) is not applied to \(D\). By Theorem~\ref{thm:maincoh} with \(m=3\): \(\scomp(U) = 1 - a_p(D)/p^3\) for all good \(p\) (in particular, for \(p \neq 3\); for \(p=2\) the reduction of \(D\) is smooth). For \(p = 3\) the reduction of \(D\) is singular (\(3 \mid 3x_i^2\)), so \(p=3 \in \Pbad\), and Theorem~\ref{thm:maincoh} is inapplicable to \(p=3\). The value \(\sigma_{\mathrm{comp}}(3) = 1\) in Table~\ref{tab:Fermat} is obtained by direct computation via Definition~\ref{def:sreg}: \(\#\PP^3(\FF_3) = 1+3+9+27 = 40\), \(\#D(\FF_3) = 13\) (by Remark~\ref{rem:Fermat3}), \(\#U(\FF_3) = 40-13 = 27\); then \(\#U/p^3 = 27/27 = 1\), \(\sreg = (1-1/3)\cdot 1 = 2/3\), \(\sigma_{\mathrm{comp}} = \sreg/(1-1/p) = 1\). The quantity \(a_3 = 0\) in the table is a conventional notation under which \(\#U/p^3 = 1-a_3/p^3\); it is not the trace of Frobenius on the primitive cohomology of a smooth surface (\(p=3\) is a bad prime for \(D\)).

Parameters: \(\kappa = -1\), \(\beta_0 = 2\), \(A = 6\). By Lemma~\ref{lem:Fermat}: \(a_p = 0\) for \(p \equiv 2 \pmod 3\), \(a_p = 6p\) for \(p \equiv 1 \pmod 3\).

\begin{table}[htbp]
\centering
\caption{Local factors for the Fermat cubic surface (\(m=3\))}\label{tab:Fermat}
\begin{tabular}{cccccc}
\toprule
\(p\) & \(\#D\) & \(a_p\) & \(\sigma_{\mathrm{comp}}\) & \(C_P\) & \(|\sigma-1| p^2\) \\
\midrule
2 & 7 & 0 & 1.000000 & 1.000000 & 0 \\
3 & 13 & 0 & 1.000000 & 1.000000 & 0 \\
5 & 31 & 0 & 1.000000 & 1.000000 & 0 \\
7 & 99 & 42 & 0.877551 & 0.877551 & 6 \\
11 & 133 & 0 & 1.000000 & 0.877551 & 0 \\
13 & 261 & 78 & 0.964497 & 0.846398 & 6 \\
17 & 307 & 0 & 1.000000 & 0.846398 & 0 \\
19 & 495 & 114 & 0.983380 & 0.832329 & 6 \\
23 & 553 & 0 & 1.000000 & 0.832329 & 0 \\
29 & 871 & 0 & 1.000000 & 0.832329 & 0 \\
\bottomrule
\end{tabular}
\end{table}

\emph{Footnote to Table~\ref{tab:Fermat}.} \(p=3\) is a bad prime for \(D\) (singular reduction); its local factor \(\sigma_{\mathrm{comp}}(3)=1\) is computed separately by Definition~\ref{def:sreg} and included in \(\Pcomp(Y)\) as an element of \(\Pbad\), independently of \(Y\) (Definition~\ref{def:bad}). Theorem~\ref{thm:maincoh} is inapplicable to \(p=3\).

\textbf{Analysis.}
\begin{enumerate}
\item For \(p \equiv 1 \pmod 3\): \(\left| \sigma - 1 \right| p^2 = 6 = A\); the bound is attained exactly.
\item For \(p \equiv 2 \pmod 3\) and \(p = 3\): \(\sigma_{\mathrm{comp}} = 1\) exactly.
\item The cumulants \(C_P = \prod_{p \le P} \scomp\) tend to \(C_\infty = \prod_{p \equiv 1 \pmod 3} \left( 1 - 6/p^2 \right) > 0\); \(C_{29} \approx 0.8323\).
\item Consistent with Theorem~\ref{thm:structural}: \(m = 3\), \(w = 2\), \(\beta_0 = 2\).
\end{enumerate}

\begin{remark}\label{rem:Cinfty}
The quantity \(C_\infty\) is not used as part of the method; the method computes \(\ccomp(U)\) via \(\Pcomp\) and \eqref{eq:ccomp}. \(C_\infty\) is a reference for verification.
\end{remark}

\subsection{Experiment: del Pezzo (two examples)}

\subsubsection{Example 1: the Fermat cubic \(C \subset \PP^2\) (ambient surface: \(\PP^2\) as a del Pezzo surface of degree \(d = 9\), \(\rho = 1\))}\label{sec:cubicexp}

\(X = \PP^2\) (a del Pezzo surface of degree \(d = 9\), \(\rho = 1\)), \(D = C = \{x^3+y^3+z^3=0\}\) a cubic curve (the degree of \(C\) equals 3), \(U = X \setminus C\). Parameters: \(A = 2\), \(\beta_0 = 3/2\). By Remark~\ref{rem:cubic}: \(\sigma_{\mathrm{comp}} = 1 + a_p(C)/p^2 + O(p^{-3})\) for good \(p\); for \(\rho=1\) the truncated formula is exact for good \(p\); for \(p=3\), separately. The Frobenius action on \(\Pic\) is trivial, so \(L_p^{-1} = 1 - 1/p\) for all \(p\), the compensator does not contribute, and global certification applies.

\textbf{Archimedean factor.} \(\rk\Pic(U) = 0\) (see Remark~\ref{rem:cubic}); by Definition~\ref{def:alpha} (normalization by \(\Pic(U)_{\RR}^\vee\)) \(\alpha(U) = 1\) under the adopted convention (delta-measure of unit mass at zero). The full constant
\[
\ccomp(U) = \prod_{p\notin\Pbad}\scomp(U)\cdot\prod_{p\in\Pbad}\scomp(U)
\]
(the value is computed with certified accuracy via the finite product and the tail estimate; equality to one is not established in the paper). Local factors and the tail are certified (\(\Areg = 2\)). The prime \(p=3\) is bad: \(3 \in \Pbad\); \(\sigma_{\mathrm{comp}}(3) = 1\) is computed separately by Definition~\ref{def:sreg} and enters \(\Pcomp\) as an element of \(\Pbad\).

\begin{table}[htbp]
\centering
\caption{Local factors for the Fermat cubic in \(\PP^2\) (\(d=9\), \(\rho=1\))}\label{tab:cubic}
\begin{tabular}{ccccc}
\toprule
\(p\) & \(a_p\) & \(\sigma_{\mathrm{comp}}\) & \(C_P\) & \(|a_p|/\sqrt{p}\) \\
\midrule
2 & 0 & 1.0000 & 1.0000 & 0 \\
3 & 0 & 1.0000 & 1.0000 & 0 \\
5 & 0 & 1.0000 & 1.0000 & 0 \\
7 & \(-1\) & 0.9796 & 0.9796 & 0.378 \\
11 & 0 & 1.0000 & 0.9796 & 0 \\
13 & \(5\) & 1.0296 & 1.0086 & 1.387 \\
17 & 0 & 1.0000 & 1.0086 & 0 \\
19 & \(-7\) & 0.9806 & 0.9890 & 1.606 \\
23 & 0 & 1.0000 & 0.9890 & 0 \\
29 & 0 & 1.0000 & 0.9890 & 0 \\
31 & \(-4\) & 0.9958 & 0.9849 & 0.718 \\
37 & \(11\) & 1.0080 & 0.9928 & 1.809 \\
43 & \(8\) & 1.0043 & 0.9971 & 1.220 \\
\bottomrule
\end{tabular}
\end{table}

\textbf{Analysis.} \(|a_p|/\sqrt{p} \le 2\) for all \(p\); \(\Aobs \approx 1.809\) (maximum on the verified range \(p\le 43\)), \(\Areg/\Aobs \approx 1.10\). The truncated formula is exact for good \(p\), no recomputation is required. The partial products \(C_P\) approach the limit; \(C_{43}\approx0.9971\) is an intermediate value, not the limit.

\subsubsection{Example 2 (local experiment): the Fermat cubic surface with elliptic boundary, \(p \equiv 1 \pmod 3\)}

\(S = \{x_0^3+x_1^3+x_2^3+x_3^3=0\} \subset \PP^3\), \(E = S \cap \{x_0=0\}\) (an elliptic curve), \(U' = S\setminus E\). Parameters: \(d=3\), \(\rho=7\), \(\Areg = 23\) by Theorem~\ref{thm:classIII}, \textbf{only for split \(p\)} (\(p \equiv 1 \pmod 3\)). For non-split \(p\) no uniform estimate is proved in the paper; the full \(\ccomp\) for \(\rho=7\) \textbf{is not computed} (see \S\ref{sec:open}). For \(p \equiv 2 \pmod 3\) the surface is non-split (\(\rk\Pic(S_{\FF_p}) = 3\)), and \(\scomp\) is not given by formula \eqref{eq:delpezzo}; the computation requires explicit \(L_p\) factors with non-trivial Frobenius action on \(\Pic(S_{\QQbar})\).

Table~\ref{tab:Fermat3} refers \textbf{only to split good primes} \(p\equiv1\pmod 3\) and is a \textbf{local experiment} verifying formula \eqref{eq:delpezzo} for such \(p\). It is \textbf{not a global certification} of the full product \(\prod_{p>Y}\scomp(U')\): for non-split \(p\equiv2\pmod 3\) the surface \(S\) is non-split (\(\rk\Pic(S_{\FF_p})=3\)), formula \eqref{eq:delpezzo} is inapplicable, and the full product requires explicit \(L_p\) factors with non-trivial Frobenius action on \(\Pic(S_{\QQbar})\). The constant \(\Areg=23\) (for \(\rho=7\)) is proved for split good primes; it \textbf{does not give} a certified tail of the full product in the presence of non-split primes. No uniform estimate for non-split \(p\) is established in the paper; hence the table is a local experiment, not a certification.

\textbf{Global certification for del Pezzo} is established in the present paper only for the case \(\rho = 1\) (the cubic in \(\PP^2\), \(d=9\)), where the Frobenius action on \(\Pic\) is trivial, so \(L_p^{-1} = 1 - 1/p\) and the local factor is explicitly known. For \(d=9\) formula \eqref{eq:delpezzo} is exact for good \(p\) (for \(p=3\), separately), and \(\Areg = 2\) is a uniform constant (see \S\ref{sec:cubicexp}). For \(d = 1, 3, 4\) the full \(\ccomp\) requires additional computations of \(L_p\) and a uniform estimate for non-split \(p\); in the paper these cases are either not considered, or (for \(d=3\)) restricted to a local experiment.

For \(p \equiv 1 \pmod 3\), by Definition~\ref{def:sreg}, Theorem~\ref{thm:delpezzo} and \(L_p^{-1} = (1-p^{-1})^7\):
\[
\scomp(U') = (1-1/p)^6 \cdot \frac{\#U'(\FF_p)}{p^2}, \qquad \#U' = p^2+6p+a_p(E).
\]

\begin{table}[htbp]
\centering
\caption{Local experiment: Fermat cubic surface with elliptic boundary, split \(p\equiv1\pmod3\)}\label{tab:Fermat3}
\begin{tabular}{cccc}
\toprule
\(p\) & \(a_p(E)\) & \(\sigma_{\mathrm{comp}}\) (exact) & \(C_P\) \\
\midrule
7 & \(-1\) & 0.728393 & 0.728393 \\
13 & \(5\) & 0.922447 & 0.671795 \\
19 & \(-7\) & 0.937243 & 0.629578 \\
31 & \(-4\) & 0.976968 & 0.615077 \\
37 & \(11\) & 0.992804 & 0.610649 \\
43 & \(8\) & 0.993250 & 0.606527 \\
\bottomrule
\end{tabular}
\end{table}

\textbf{Analysis.}
\begin{enumerate}
\item For all \(p \equiv 1 \pmod 3\) in Table~\ref{tab:Fermat3} the exact formula is applied to the exact \(\#U'\).
\item For \(p=7\) the compensator cancels the regularizer, giving \(\sigma_{\mathrm{comp}} \approx 0.7284\); with the truncated formula \(1+(a_p-21)/p^2\) one would get \(0.5510\), an error of about 32 percent. This illustrates why the finite product requires exact values.
\item The cumulants \(C_P\) over split primes converge: \(C_{43} \approx 0.6065\).
\item For \(p \equiv 2 \pmod 3\) (not included) the formula \(\#U' = p^2+6p+a_p\) is incorrect; the full product \(\Pcomp(Y)\) requires taking both groups of primes into account with the correct \(L_p\).
\item For all \(p\) in the table \(|\sigma-1| p^{3/2} \le \Areg = 23\) (Theorem~\ref{thm:classIII}).
\end{enumerate}

\begin{remark}\label{rem:localext}
The values of \(\scomp\) are obtained by direct enumeration of points of \(S(\FF_p)\) and \(E(\FF_p)\) for \(p \le 43\), \(p \equiv 1 \pmod 3\). For non-split \(p\) the computation of \(L_p\) with non-trivial Frobenius action on \(\Pic(S_{\QQbar})\) is beyond the scope of the paper.
\end{remark}

\subsection{Conservativeness of the bounds}

\(\Aobs\) in Table~\ref{tab:conserv} is the \textbf{maximum} of \(|\sigma-1| p^{3/2}\) over the verified finite range of primes, not the supremum over all primes.

\begin{table}[htbp]
\centering
\small
\caption{Conservativeness of the bounds: \(\Areg\) and \(\Aobs\) on the verified range of primes}\label{tab:conserv}
\begin{tabular}{@{}p{5.5cm}cccc@{}}
\toprule
Case & \(\Areg\) & \(\Aobs\) & \(\Areg/\Aobs\) \\
\midrule
Torus, hyperplanes & 0 & 0 & not applicable \\
Fermat surface (\(m=3\)) & 6 & 6 (on \(p\le 29\)) & 1.00 \\
Fermat cubic in \(\PP^2\) (\(d=9\)) & 2 & \(\approx 1.81\) (on \(p\le 43\)) & \(\approx 1.10\) \\
Fermat cubic surface with elliptic boundary, split \(p\equiv 1\pmod 3\) & 23 & \(\approx 5.20\) (on \(p\le 43\), \(p\equiv 1\)) & \(\approx 4.4\) \\
\bottomrule
\end{tabular}
\end{table}

\begin{remark}\label{rem:Acert23}
The value \(\Areg\) for \(d=3\): the direct inequality of the proof of Theorem~\ref{thm:classIII} gives \(\Areg=2+\rho(\rho-1)/2=23\) for \(\rho=7\). This value is proved \textbf{only for split \(p\)}; for non-split \(p\) no uniform estimate is established, and the full \(\ccomp\) for \(\rho=7\) remains open (see \S\ref{sec:open}).
\end{remark}

\section{Concluding remarks}\label{chap:conclusion}

\subsection{What has been done}\label{sec:summary}

\textbf{Analytic compensation.} The compensated local factor has been introduced (Definition~\ref{def:compensator}); absolute convergence of the compensated product has been proved (Theorem~\ref{thm:conv}) for classes I, II, and class III with \(\rho=1\); a decomposition of \(\ccomp(U)\) into three blocks has been obtained (Theorems~\ref{thm:decomposition} and \ref{thm:blocks}). Bad primes have been separated into \(\Pbad\) (Definition~\ref{def:bad}) and are multiplied in \(\Pcomp(Y)\) independently of \(Y\).

\textbf{Estimate of the tail.} Theorem~\ref{thm:truncation}; threshold \eqref{eq:inverse}; structural relation \(\beta_0 = m - w/2\) via Deligne majorants (Theorem~\ref{thm:Deligne}). Lemma~\ref{lem:tail} is stated for integer \(Y\).

\textbf{Parameters of classes I, II, III.} The values \(\kappa\), \(\Areg\), \(\beta_0\) have been computed for classes I, II, III (see Corollary~\ref{cor:classification}); within the classes, hyperplanes, the torus, hypersurfaces \(m \ge 3\), open del Pezzo surfaces, and the cubic in \(\PP^2\) are treated separately. Theorem~\ref{thm:structural} (class II) has been proved; for class III the estimate \(\beta_0=3/2\) has been proved by the direct inequality of Theorem~\ref{thm:classIII} (for split \(p\); globally, for \(\rho=1\)), for class I by Theorem~\ref{thm:classI}. Classification is in Corollary~\ref{cor:classification}; the comparison of \(\Areg\) and \(\Aobs\) is in Remark~\ref{rem:twoConstants}. For class III, \(\Areg\) is uniform over split \(p\) (Theorem~\ref{thm:classIII} with the direct inequality); for non-split \(p\) with \(\rho>1\) a uniform estimate is open. The exponent \(\beta\) in the tables and the algorithm is the \textbf{proved lower bound} \(\beta_0\) of Definition~\ref{def:beta0} (it is not assumed optimal; the true \(\beta(X,D)\) of Definition~\ref{def:beta} may be larger).

\textbf{Inverse problem.} Table~\ref{tab:thresholds}; three classes \S\ref{chap:inverse}; structural chain \eqref{eq:chain}.

\textbf{Reduction to the archimedean factor.} Theorem~\ref{thm:reduction} (classes I, II, III with \(\rho=1\)); for class I, Theorem~\ref{thm:classI_exact} (taking bad primes into account).

\textbf{Relation to \(N_U\).} Direct interval form (Corollary~\ref{cor:interval}); Theorem~\ref{thm:errorlink}; Remark~\ref{rem:intervalfull}; Corollary~\ref{cor:chooseY}. Interval arithmetic for \(\alpha\) and \(\Pcomp\) in \S\ref{sec:algo}.

\textbf{Algorithm and experiment.} \S\ref{sec:algo}, Lemma~\ref{lem:complexity}, experiments \S\ref{sec:Fermat} and \S\ref{sec:cubicexp}.

\begin{table}[htbp]
\centering
\caption{Parameters of classes I, II, III: \(\kappa(U)\), \(w\), \(\Areg(U)\), \(\beta_0\), \(Y(10^{-2})\)}\label{tab:summary}
\begin{tabular}{ccccccc}
\toprule
Class & Type & \(\kappa(U)\) & \(w\) & \(\Areg(U)\) & \(\beta_0\) & \(Y(10^{-2})\) \\
\midrule
I & Hyperplanes & \(-k\) & \(\infty\) & 0 & \(\infty\) & 2 \\
I & Torus & \(-(m+1)\) & \(\infty\) & 0 & \(\infty\) & 2 \\
II & Cubic surface (\(D\subset\PP^3\)) & \(-1\) & 2 & 6 & 2 & 2400 \\
II & K3 surface (\(D\subset\PP^3\)) & \(-1\) & 2 & 21 & 2 & 8400 \\
II & Cubic threefold (\(D\subset\PP^4\)) & \(-1\) & 3 & 10 & 5/2 & 193 \\
III & Del Pezzo \(d=1\) (split \(p\)) & \(-1\) & 1 & 38 & 3/2 & \(\approx 9.24\cdot 10^8\) \\
III & Del Pezzo \(d=4\) (split \(p\)) & \(-1\) & 1 & 17 & 3/2 & \(\approx 1.85\cdot 10^8\) \\
III & Del Pezzo \(d=9\) & \(-1\) & 1 & 2 & 3/2 & \(2.56\cdot 10^6\) \\
III & Cubic in \(\PP^2\) & \(-1\) & 1 & 2 & 3/2 & \(2.56\cdot 10^6\) \\
\bottomrule
\end{tabular}
\end{table}

\subsection{Bounds of the method}\label{sec:bounds}

\textbf{Geometric.} The paper covers logarithmic varieties with an SNC divisor for which condition (L) holds with \(\beta_0 > 1\). This includes both logarithmic Fano and logarithmic Calabi and Yau varieties (torus, cubic in \(\PP^2\), cubic surface with elliptic boundary). Mode A applies to classes I, II, and class III with \(\rho=1\) for a given \(\alpha(U)\); for class III with \(\rho>1\) the full \(\ccomp\) is not computed (condition (L) for all good \(p\) is not proved; see \S\ref{sec:bounds}, \S\ref{sec:open}). Mode B is stated with a given model \(\UZZ\), height \(H\), \(M_U(B)\) and the known equality \(\ccomp = \cManin\). Outside the domain of applicability: non-isolated singularities, non-SNC divisors. No condition on the cohomology of \(D\) is imposed.

\textbf{On \(\beta_0\).} \(\beta_0 > 1\) is required. Cases \(\beta_0 \le 1\) are not considered; before compensation \(\beta_{\log} = 1\) (Corollary~\ref{cor:divergence}).

\textbf{On \(\alpha(U)\).} The computation of \(\ccomp(U)\) is reduced to \(\alpha(U)\) under the condition that \(\alpha(U)\) is a given finite positive computable quantity. For class I, \(\alpha(U)\) is given by a finite-dimensional integral; given \(\alpha(U)\), the reduction yields unconditional computability of \(\ccomp(U)\). For logarithmic Calabi and Yau varieties with \(-K_X - D = 0\) the integral in the definition of \eqref{eq:alpha} \textbf{over \(\Pic(U)_{\RR}^\vee\)} is finite (see Remark~\ref{rem:spacechoice}): for \(r = 0\) it equals 1 with respect to the delta-measure of unit mass; for \(r \ge 1\) it diverges and requires separate normalization (not given in the paper). The construction of an algorithm for \(\alpha(U)\) in general logarithmic pairs is not considered in the paper. The question of how the complexity of computing \(\alpha(U)\) compares with that of computing \(\ccomp(U)\) for a concrete pair of classes II and III is \textbf{open}; it is not considered in the paper. The reduction (Theorem~\ref{thm:reduction}) does not assert that \(\alpha(U)\) is easier than \(\ccomp(U)\): it merely reduces \(\ccomp\) to \(\alpha\) and separates the arithmetic part from the archimedean one (see Remark~\ref{rem:alphacomplexity} and \S\ref{sec:compensation}). For the cubic in \(\PP^2\) (\(r=0\)) \(\alpha(U) = 1\) under the adopted convention; agreement with the normalization of \cite[\S5]{Peyre1995} for \(r=0\) is not verified in the paper.

\textbf{On the model and height.} Mode B applies only to those integral models and heights for which decomposition \eqref{eq:logManin} with \(\ccomp = \cManin\) is proved. For the torus the standard model of \(\ZZ\)-points \((\GG_m)^m\) does not satisfy this condition; a model is needed for which the Manin and Peyre formula for \textbf{integral} points is established. In the present paper such a model is not fixed; the model, height, and \(M_U(B)\) are part of the input data of mode B. For general logarithmic Fano varieties the equality \(\ccomp=\cManin\) \textbf{as geometric adelic constants} is proved only for hyperplanes and the torus under the fixed normalization (Corollary~\ref{cor:recovery}, normalization convention). For integral points in mode B this equality is part of Assumption~\ref{ass:logManin} together with the model and height. For del Pezzo and hypersurfaces in \(\PP^m\) the equality \(\ccomp = \cManin\) is assumed by the Manin conjecture, and a reference to a concrete model and height is part of the input data of mode B.

\textbf{On the remainder and full certification.} For a certified interval for \(N_U\): \(M_U(B)\) is given; the Manin conjecture is proved and \(r_U\) has a power-saving estimate \cite{Wilsch2023,BatyrevTschinkel1995,BBDG2024}. For del Pezzo \textbf{certification of the compensated product and its tail} is established \textbf{only for \(\rho = 1\)} (\(d=9\), the cubic in \(\PP^2\)), where the Frobenius action on \(\Pic\) is trivial, so \(L_p^{-1} = 1 - 1/p\) and the local factor is explicitly known. For \(\rho > 1\) (\(d = 1, 3, 4\)) the full \(\ccomp\) requires: (a) a uniform estimate \(|\sigma_{\mathrm{comp}}-1|\le A'p^{-3/2}\) for \textbf{non-split} good \(p\); (b) computation of \(L_p\) for non-split \(p\). Neither (a) nor (b) is done in the paper; the constant \(\Areg=2+\rho(\rho-1)/2\) is proved only for split \(p\). Example \(d=3\) (\S\ref{sec:cubicexp}) is a local experiment over split \(p\) that does not give the full \(\ccomp\). Thus \textbf{full certification of class III in the present paper is limited to the case \(\rho=1\)}.

\textbf{On complexity.} Class III: \(O(Y^2)\) (naive enumeration) or \(\tilde{O}(Y)\) (Schoof or SEA algorithm). For \(\rho>1\) the dependence on \(\rho\) is added: \(\Areg=2+\rho(\rho-1)/2\), and for \(\beta_0=3/2\) the threshold \(Y=(8A/\eta)^2\) grows as \(A^2\); for \(\rho\gg1\) this is \(\sim\rho^4\eta^{-2}/4\). Specifically: for \(d=9\) (\(\rho=1\)) \(A=2\), \(Y(10^{-2})=2.56\cdot10^6\); for \(d=1\) (\(\rho=9\)) \(A=38\), \(Y(10^{-2})\approx9.24\cdot10^8\), 361 times larger. Moreover, for non-split \(p\) the computation of \(L_p(1)^{-1}=\prod_{j=1}^\rho(1-\alpha_j/p)\) requires the characteristic polynomial of Frobenius on \(\Pic(S_{\QQbar})\cong\ZZ^\rho\) and its roots, a problem of rank \(\rho\) for each prime; for split \(p\) one needs split detection \(\rk\Pic(S_{\FF_p})=\rho\). Both operations for \(\rho>1\) add orders of magnitude to the cost of \(\Pcomp(Y)\) compared with \(\rho=1\), where \(L_p^{-1}=1-1/p\) and all \(p\) are split. Thus for \(\rho>1\) the computational cost is practically prohibitive independently of the certification problem of \S\ref{sec:open}, items 8 and 9.

\subsection{Open problems}\label{sec:open}

\begin{enumerate}
\item An explicit estimate of \(\Aobs\) for typical K3 and general hypersurfaces.
\item Extension of experiments to \(d = 1, 4\) and a typical K3, including non-split \(p\).
\item Going beyond \((H_0)\) and \((H_1)\).
\item Classes with weight \(w\) different from the weight of middle cohomology.
\item Cases \(\beta_0 \le 1\).
\item Other logarithmic conjectures as a communication channel.
\item An effective algorithm for computing \(\alpha(U)\) in general logarithmic pairs, including logarithmic Calabi and Yau varieties. For classes II and III (Theorem~\ref{thm:reduction}) the question of the complexity relation remains open: whether \(\ccomp(U)\) reduces to \(\alpha(U)\) in polynomial time (the direct reduction is given by Theorem~\ref{thm:reduction}), and whether \(\alpha(U)\) reduces to \(\ccomp(U)\) (the reverse reduction, an open question).
\item \textbf{A uniform estimate for non-split \(p\).} To prove \(|\scomp(U)-1|\le A'p^{-3/2}\) with explicit \(A'\) for non-split good \(p\) with \(\rho>1\) (del Pezzo). Without this the full \(\ccomp\) of class III with \(\rho>1\) remains open; the known estimate \(\Areg=2+\rho(\rho-1)/2\) is proved only for split \(p\).
\item \textbf{Computation of \(L_p\) for non-split \(p\)} when working with the full \(\Pcomp(Y)\): explicit \(L_p\) factors with non-trivial Frobenius action on \(\Pic(S_{\QQbar})\).
\item \textbf{Computational optimization for \(\rho>1\).} The threshold \(Y=(8A/\eta)^2\) grows as \(A^2\), where \(A=2+\rho(\rho-1)/2\); for \(\rho=9\) (\(d=1\)) \(Y(10^{-2})\approx9.24\cdot10^8\), two orders of magnitude larger than for \(\rho=1\). Computing \(L_p\) for non-split \(p\) requires the characteristic polynomial of Frobenius on \(\Pic\) of rank \(\rho\) and its roots; split detection is a separate operation for each \(p\).
\end{enumerate}

\subsection{General conclusion}

The goal of the paper is to construct a method for estimating the number of integral solutions of Diophantine equations with a prescribed accuracy. This is a classical problem of number theory. Arithmetic geometry (logarithmic varieties, cohomology of divisors, local densities, Deligne and Hasse estimates, the archimedean factor) is used here as a tool, not as the object of investigation. The final construction:
\[
(X,D) \;\rightsquigarrow\; \sreg \;\rightsquigarrow\; \scomp \;\rightsquigarrow\; \ccomp(U) \;\rightsquigarrow\; \ccomp(U)\, M_U(B) \;\rightsquigarrow\; [L, U].
\]

The main results are algorithmic: the computation of \(\ccomp(U)\) reduces to the effective computation of \(\alpha(U)\) under conditions (i)--(iii) of mode A (Theorem~\ref{thm:reduction}) for classes I, II, and class III with \(\rho=1\); for class I, given \(\alpha(U)\), we have fully unconditional computability (Theorem~\ref{thm:classI_exact}, Corollary~\ref{cor:recovery}). Given the main-growth function \(M_U(B)\), the known equality \(\ccomp = \cManin\), and the known majorant \(R_U(B)\), a certified interval for \(N_U(\ZZ, B)\) is obtained from \(\ccomp(U)\) (Corollary~\ref{cor:interval} and the algorithm of chapter~\ref{chap:algo}). The central chapters are \ref{chap:reduction} to \ref{chap:algo}; chapters~\ref{chap:compensation} to \ref{chap:tail} are preparatory; chapters~\ref{chap:classes} to \ref{chap:inverse} are applied.

The unconditional part of the chain ends at \(\ccomp(U)\): its computation reduces to \(\alpha(U)\), and for class I, given \(\alpha(U)\), we have full unconditionality. The passage to \([L, U]\) additionally requires \(M_U(B)\), \(\ccomp = \cManin\), and the majorant \(R_U(B)\).

The Manin conjecture is used only as a communication channel. When it is proved for a given class and the remainder is estimated, the method yields an unconditional interval (mode B); when not, the method yields the computed \(\ccomp(U)\) (in the reduction mode, given a subroutine for \(\alpha(U)\); for class I, unconditionally given \(\alpha(U)\)) and an asymptotic interpretation. The separation of modes is a fundamental architectural feature of the method.

Two refinements fixing the architecture of the method: (i) the interval formula of Corollary~\ref{cor:interval} uses the \textbf{full} majorant of the deviation of \(T_Y\) from 1 and the direct form of the product without passing to the relative error \(c_Y/\ccomp\); (ii) mode B explicitly takes the integral model, the height, and the function \(M_U(B)\) as input data, and does not assume that \(\ccomp = \cManin\) automatically relates the computed constant to \(N_U(\ZZ, B)\). For the torus \((\GG_m)^m\) under the standard model, \(\ZZ\)-points reduce to \((\pm1)^m\), and mode B requires a model for which the Manin and Peyre formula for integral points is established.

The experiments (\S\ref{sec:Fermat}, \S\ref{sec:cubicexp}) confirm the operability of the method. \textbf{Certified examples for the compensated Euler product and its tail:} the Fermat surface (\(m=3\), class II; \(\Areg=6\) is attained exactly), the Fermat cubic in \(\PP^2\) (\(d=9\), class III, \(\rho=1\); \(\Aobs\approx 1.81<2=\Areg\) on \(p\le43\)), hyperplanes and the torus (class I; \(\sigma_{\mathrm{comp}}\equiv1\)). \textbf{Local experiment:} the Fermat cubic surface with elliptic boundary (\(d=3\), \(\rho=7\)), formula \eqref{eq:delpezzo} verified for split \(p\equiv1\pmod3\), \(p\le43\); \(\Aobs\approx5.20<23=\Areg\); the full \(\ccomp\) for \(\rho=7\) is open. \textbf{Certification} here means certification of \(\ccomp(U)\) and the tail \(\Tcomp(U)\), and \textbf{not} a numerical interval for \(N_U(\ZZ,B)\): the latter additionally requires concrete \(\alpha(U)\), a model, a height, and a proved remainder estimate (mode B, Assumption~\ref{ass:logManin}). \textbf{Full certification of class III is limited to the case \(\rho=1\)}; for \(\rho>1\) the results refer to split \(p\) and are local, and the computational cost is practically prohibitive independently of certification (\S\ref{sec:bounds}, \S\ref{sec:open} item 10).

The results apply in classes where \(X\) satisfies \((H_0)\) or \((H_1)\), \(\beta_0 > 1\), and \(\alpha(U)\) is given as a finite computable quantity. Extension to wider classes is the subject of further work.

\end{document}